\documentstyle[11pt]{article}

\begin{document}
\noindent
\begin{center}
  {\LARGE Higher Genus Symplectic Invariants and Sigma Model Coupled With
  Gravity}
  \end{center}

  \noindent
  \begin{center}
    {\large Yongbin Ruan}\footnote{partially supported by a NSF grant and a
    Sloan
    fellowship}\\[5pt]
      Department of Mathematics, University of Utah\\
	Salt Lake City, UT 84112\\[5pt]
	  {\large Gang Tian}\footnote{partially supported by
	  NSF grants}\\[8pt]
	    Department of Mathematics, MIT\\
	      Cambridge, MA 02139
	      \end{center}

	      \def \J{{\cal J}}
	      \def \Map{Map(S^2, V)}
	      \def \M{{\cal M}}
	      \def \A{{\cal A}}
	      \def \B{{\cal B}}
	      \def \C{{\bf C}}
	      \def \Z{{\bf Z}}
	      \def \R{{\bf R}}
	      \def \P{{\bf P}}
	      \def \I{{\bf I}}
	      \def \N{{\cal N}}
	      \def \T{{\bf T}}
	      \def \Q{{\bf Q}}
	      \def \D{{\cal D}}
	      \def \H{{\cal H}}
	      \def \S{{\cal S}}
	      \def \e{{\bf E}}
	      \def \c{{\cal C}}
	      \def \U{{\cal U}}
	      \def \E{{\cal E}}
	      \def \F{{\cal F}}
	      \def \L{{\cal L}}
	      \def \K{{\cal K}}
	      \section{ Introduction}

	       This paper is a continuation of our previous paper \cite{RT}.
In \cite{RT}, among other things, we build up the mathematical
foundation of quantum cohomology ring  on semi-positive symplectic manifolds.
      We also defined higher genus symplectic invariants without gravity
(topological sigma model) in terms of inhomogeneous
holomorphic maps from a fixed Riemann surface,
and proved the composition law they satisfy.
Topological gravity, proposed by Witten, concerns the intersection
theory of the moduli space of marked Riemann surfaces. Based
on the physical intuition, Witten suggested a relation between
those intersection numbers and the KdV hierarchy.
This relation was  clarified by Kontsevich (cf. \cite{Ko}). However,
the mathematical, as well as physical, phenomenon will become much more
interesting
if the topological sigma model
is coupled with the topological gravity.
In fact, in \cite{W2} Witten proposed an approach to
the topological sigma model coupled with gravity,
and made a very important conjecture
on the basic feature of this new model. The purpose of
this paper is to establish a mathematical foundation for the
theory of topological sigma model coupled with topological gravity over
any semi-positive symplectic manifolds. This new theory also provides
many more new geometric examples of the topological field
theory coupled with gravity.
For each semi-positive symplectic manifold $V$, we can associate a
topological sigma model with gravity, or
simply a topological field theory coupled with gravity.
This theory begins with the GW-invariants
$$ \Psi ^V_{(A,g,k)} : H_*(\overline \M_{g,k}, \Q)
\times H_*(V, \Z)^k\mapsto \Q,$$
for any $A\in H_2(V, \Z)$ and $2g +k \ge 3$. Here
$\overline \M_{g,k}$ is the Deligne-Mumford compactification
of the moduli space of  genus $g$ Riemann surfaces with $k$
marked points. The GW-invariants are multilinear and
supersymmetric on $H_*(V,\Z)^k$.

At first, we will rigorously define the GW-invariant $\Psi^V$
on semi-positive symplectic manifolds (cf. section 2).

 From the analytic point of view,
it is the most convenient to use the inhomogeneous
holomorphic maps from Riemann
surfaces in $\overline \M_{g,k}$, though
other equivalent formulations may be possible, such as using stable maps
and establishing a more sophisticated intersection theory.
An inhomogeneous holomorphic map is a solution of an
inhomogeneous Cauchy-Riemann equation (cf. Section 2).

Putting aside technical details for the time being, we can intuitively
define the GW-invariants (cf. Section 2 for details) as follows:
let $V$ be any symplectic manifold and $A\in H_2(V, \Z)$. For any
homology classes $[K] \in H_*(\overline \M_{g,k}, \Q)$ and
$\alpha _i\in H_*(V, \Z)$, represented by cycles $K$, $A_i$,
respectively, we define
$\Psi^V_{(A,g,k)}([K]; \alpha _1, \cdots, \alpha _k)$ to be
the number of tuples $(\Sigma ; x_1, \cdots, x_k; f)$
with appropriate sign, satisfying:
$\Sigma \in K$, $f: \Sigma \mapsto V$ solves a given inhomogeneous
Cauchy-Riemann equation, and $f(x_i) \in A_i$, whenever
$$\sum cod(A_i)+cod(K)= 2c_1(V)(A)+2(3-n)(g-1)+2k;
\eqno(1.1)$$
We simply put $\Psi^V_{(A,g,k)}([K]; \alpha _1, \cdots, \alpha _k)$
to be zero if (1.1) is not satisfied.

This approach towards defining new invariants has been used
before in many cases (cf. [Do], [Gr], [R], [R3], [RT],  [W1]).
For symplectic 4-manifolds, using unperturbed holomorphic maps,
the first author already defined the invariant $\Psi$ in the very important
case that
$k=0$ and $[K] = \overline \M_{g,0}$.
However, in each case, there are specified difficulties
to be overcome. Using the techniques we developed in [RT], we
will first prove
\vskip 0.1in
\noindent
{\bf Theorem A (Theorem 2.14.)} {\it If $V$ is a semi-positive symplectic
manifold,
the GW-invariant $\Psi^V_{(A,g,k)}$ can be well defined
for any $g, k \ge 0$ with $2g + k \ge 3$. Moreover,
this $\Psi^V$ depends only the symplectic structure of $V$.}
\vskip 0.1in
A symplectic manifold $V$ is semi-positive if
it is compact and there is no $J$-holomorphic
map $f: S^2 \mapsto V$ such that $3-n \le \int_{S^2} f^*c_1(V) < 0$,
where $J$ is any given compatible almost complex structure on $V$.
In particular, any algebraic manifold of dimension $\le 3$ is
semi-positive in this sense, also any algebraic manifold
$V$ with $c_1(V) \ge 0$ is semi-positive.

One new consequence of our theorem, which was not
obvious at all to physicists based on
mathematically unjustified path integrals,
is that the invariant $\Psi^V$ is a symplectic invariant.  The path
integral starts from a Lagrangian. The Lagrangian for sigma model or sigma
model coupled with gravity is valid
for any almost complex manifolds (symplectic
or not). There was a speculation that its correlation functions will be
the invariants of homotopy class of almost complex structures. This is in fact
false. Our invariants are symplectic invariants rather than the invariants of
almost complex structures. In particular, they can distinguish different
 symplectic
manifolds with the same homotopy class of almost complex structures (see
 section 5 or \cite{R}, \cite{R1}).

One of fundamental properties of a topological field theory is the
axiom on the decomposition of correlation functions. In our case,
the GW-invariants serve as the correlation functions. Therefore, in order
to make them more useful, or at least to construct a correct model for
the topological
field theory, we need to verify that our invariants satisfy
the composition law.

The composition law
governs how the GW-invariants change
during the degeneration of stable curves. Its classical cousin in enumerative
algebraic geometry is the degeneration method, which was only derived
in very special cases. The classical degeneration method never became a
general theory as neat as the composition law describes.
One reason might be that the classical counting of holomorphic
curves, particularly of higher genus,
does  not obey the  composition laws predicted by physicists, even for the
projective plane $\P^2$.  Namely the way of counting was not good.
In [RT], we found the correct counting in terms of inhomogeneous
holomorphic maps and established the composition law at least for
the mixed invariants, corresponding to the $\sigma$-models
without gravity. Based on the same techniques developed in [RT], we
are also able to prove the composition law for all GW-invariants.

Assume $g=g_1+g_2$ and $k=k_1+k_2$ with $2g_i + k_i \ge 3$.
Fix  a decomposition $S=S_1\cup S_2$ of $\{1,\cdots , k\}$ with
$|S_i|= k_i$. Then there is a canonical embedding
$\theta _S: \overline \M_{g_1,k_1+1}\times \overline \M_{g_2,k_2+1}
\mapsto \overline \M_{g,k}$, which assigns to marked
curves $(\Sigma _i; x_1^i,\cdots ,x_{k_1+1}^i)$ ($i=1,2$), their
union $\Sigma _1\cup \Sigma _2$ with $x^1_{k_1+1}$ identified to
$x^2_{k_2+1}$ and remaining points renumbered by $\{1,\cdots,k\}$ according to
$S$.

There is another natural map $\mu : \overline \M_{g-1, k+2}
\mapsto \overline \M_{g,k}$ by gluing together the last two marked
points.

Choose a homogeneous
basis $\{\beta _b\}_{1\le b\le L}$ of $H_*(V,\Z)$ modulo
torsion. Let $(\eta _{ab})$ be its intersection matrix. Note that
$\eta _{ab} = \beta _a \cdot \beta _b =0$ if the dimensions of
$\beta _a$ and $\beta _b$ are not complementary to each other.
Put $(\eta ^{ab})$ to be the inverse of $(\eta _{ab})$.
Now we can state the composition law, which consists of
two formulas.
\vskip 0.1in
\noindent
{\bf Theorem B. (Theorem 2.10)}  {\it Let $[K_i] \in H_*(\overline \M_{g_i,
k_i+1}, \Q)$ $(i=1,2)$ and $[K_0] \in H_*(\overline \M_{g-1,
k +2}, \Q)$. For any $\alpha _1,\cdots,\alpha _k$ in $H_*(V,\Z)$.
Then we have}
$$\begin{array}{rl}
&\Psi ^V_{(A,g,k)}(\theta _{S*}[K_1\times K_2];\{\alpha _i\})\\
=~& \sum \limits _{A=A_1+A_2} \sum \limits_{a,b}
\Psi ^V_{(A_1,g_1,k_1+1)}([K_1];\{\alpha _{i}\}_{i\le k_1}, \beta _a)
\eta ^{ab}
\Psi ^V_{(A_2,g_2,k_2+1)}([K_2];\beta _b,
\{\alpha _{j}\}_{j> k_1}) \\
\end{array}
\leqno (1.5)
$$
$$
\Psi ^V_{(A,g,k)}(\mu_*[K_0];\alpha _1,\cdots, \alpha _k)
=\sum _{a,b} \Psi ^V_{(A,g-1,k+2)}([ K_0];\alpha _1,\cdots, \alpha _k,
\beta _a,\beta _b) \eta ^{ab}\leqno (1.6)
$$
\vskip 0.1in

There is a natural map $\pi:
\overline{\M}_{g,k}\rightarrow \overline{\M}_{g, k-1}$ as follows: For
$(\Sigma, x_1, \cdots, x_k)\in \overline{\M}_{g,k}$, if $x_k$ is not in any
rational component of $\Sigma$ which contains only three special points,
then we define
$$\pi(\Sigma, x_1, \cdots, x_k)=(\Sigma, x_1, \cdots, x_{k-1}),$$
where a distinguished point of $\Sigma$ is either a singular point or a
marked point. If $x_k$ is in one of such rational components, we contract
this component and obtain a stable curve $(\Sigma', x_1, \cdots, x_{k-1})$ in
$\overline{\M}_{g, k-1}$, and define $\pi(\Sigma, x_1, \cdots, x_k)=(\Sigma',
x_1, \cdots, x_{k-1}).$

Clearly, $\pi$ is continuous. One should be aware that
there are two exceptional cases $(g,k)=(0,3), (1,1)$ where $\pi$ is not well
defined. Associated with $\pi$,
we have two {\em k-reduction formulas} for $\Psi^V_{(A, g, k)}$.
\vskip 0.1in
\noindent
{\bf Proposition C (Theorem 2.15). }{\it Suppose that $(g,k)
\neq (0,3),(1,1)$.
\vskip 0.1in
\noindent
(1) For any $\alpha _1, \cdots , \alpha _{k-1}$ in $H_*(V, \Z)$,
we have}
$$\Psi ^V_{(A,g,k)}([K]; \alpha _1, \cdots,\alpha _{k-1}, [V])~=~
\Psi ^V_{(A,g,k-1)}([\pi (K)]; \alpha _1, \cdots,\alpha _{k-1})
\leqno (3.3)$$
\vskip 0.1in
\noindent
(2)  Let $\alpha _k$ be in $H_{2n-2}(V, \Z)$, then
$$\Psi ^V_{(A,g,k)}([\pi^{-1}(K)]; \alpha _1,
\cdots,\alpha _{k-1}, \alpha _k)~=~\alpha^* _k (A)
\Psi ^V_{(A,g,k-1)}([K]; \alpha _1, \cdots,\alpha _{k-1})
\leqno (3.4)$$
where $ \alpha^* _k$ is the Poincare dual of $\alpha _k$.
\vskip 0.1in
\noindent
In order to formulate the generalized
Witten conjecture in terms of our invariants, we need to introduce
special cycles in $\overline \M_{g,k}$. Let
$\pi : \overline{\U_{g,k}}\rightarrow \overline{\M}_{g,k}$ be
the universal family of stable curves of genus $g$ and $k$ marked points.
Each marked point gives rise to a section $\sigma _i$ ($1\le i \le k$)
of this fibration. Following Witten,
we let $\L_i$ be the pull-back of the relative
cotangent sheave of $\pi:
\overline{\U_{g,k}}\rightarrow \overline{\M}_{g,k}$ by $\sigma_i$.
Then we put $W_{d_1, \cdots , d_k}$ to be
the Poincare dual of the cohomology class
$c_1(\L_1)^{d_1}\cup c_1(\L_2)^{d_2}\cdots\cup c_1(\L_k)^{d_k}$.
We call these $W_{d_1, \cdots , d_k}$ Witten cycles.

For convenience, as Witten did, we use
$$<\tau_{d_1, \alpha_1}, \tau_{d_2, \alpha_2}, \cdots,
\tau_{d_k, \alpha_k}>_{g,k}$$
to denote the GW-invariants $\Psi_{(A,g,k)}([W_{d_1, \cdots , d_k}];
\alpha_1, \cdots, \alpha_k)$.
Following Witten, we introduce potential functions
$$F_g=\sum_A\sum_{n_{r, \alpha}}\prod_{r, \alpha}
\frac{(t^{\alpha}_r)^{n_{r,
\alpha}}}{n_{r, \alpha}!}
<\prod_{r, \alpha}\tau^{n_{r, \alpha}}_{r, \alpha}>q^A,~~~g=0, 1,2, \cdots.$$
where $q^A$ is an element of Novikov ring \cite{MS}, \cite{RT} (section 8).
We further define
$$F^V=\sum_{g\geq 0} F_g.$$
One of fundamental problems on $F^V$, even to physicists,
is to find
the complete set of equations $F^V$ satisfies.
In Section 6, imitating the arguments of Witten in [W2], we
will prove (cf. Lemma 6.1, 6.2)
\vskip 0.1in
\noindent
{\bf Theorem C. }{\it $F^V$ satisfies the generalized string
equation
$${\partial F^V\over \partial t^1_0} = {1\over 2} \eta _{ab} t^a_0 t^b_0
+ \sum \limits_{i=0}^\infty \sum \limits _{a} t^a_{i+1}
{\partial F^V\over \partial t^a_i}.\leqno(1.7)$$
$F_g$ satisfies the dilation equation
$$\frac{\partial F_g}{\partial t^1_1}=(2g-2+\sum_{i=1}^{\infty}\sum_{a}t^a_i
\frac{\partial }{\partial t^a_i})F_g+\frac{\chi(V)}{24}\delta_{g,1},
\leqno(1.8)
$$
where $\chi(V)$ is the Euler characteristic of $V$.}
\vskip 0.1in
In general, Witten suggested
$$U = {\partial ^{2} F^V \over \partial t _{0, 1}
\partial t _{0,\sigma}},~
U' = {\partial ^{3} F^V \over \partial t _{0, 1}^2
\partial t _{0,\sigma}},~\cdots,~U^{(l)}_\sigma =
{\partial ^{l+2} F^V \over \partial t _{0, 1}^{l+1}
\partial t _{0,\sigma}},~~~~{\rm for~}~ l \ge 0
$$
We will regard $U^{(l)}$ to be of degree $l$. By a differential
function of degree $k$ we mean
a function $G(U, U', U'', \cdots)$ of degree $k$
in that sense. In particular,
any function of form $G(U)$ is of degree zero, and
$(U')^2$ has degree two.
\vskip 0.1in
\noindent
{\bf Generalized Witten Conjecture:} {\it For every $g\geq 0$, there are
differential functions

$G_{m,\alpha; n, \beta}(U_{\alpha}, {U_{\alpha}}',
{U_{\alpha}}'', \cdots)$ of degree $2g$ such that
$${\partial ^2 F_g \over \partial
\tau_{m, \alpha} \partial
\tau_{n, \beta}}= G_{m,\alpha; n, \beta}(U_{\alpha},
{U_{\alpha}}', {U_{\alpha}}'', \cdots)$$
up to terms of genus $g$.}
\vskip 0.1in
This conjecture was affirmed in case $V=pt$ by Kontsevich \cite{Ko}.
When $g=0$, it is a consequence of
the associativity equation proved in \cite{RT}. But
the general case is still open.

We call $\Psi^V_{(A,g,k)}([\overline \M_{g,k}]; \cdots)$
primitive GW-invariants of genus $g$. Those invariants
correspond to
the enumerative invariants of counting genus $g$ holomorphic
curves passing through generic $k$ cycles
in enumerative algebraic geometry.
\vskip 0.1in
\noindent
{\bf Corollary E (Proposition 6.5). }
{\it For genus $\leq 1$, the Witten invariants $<>$ can be
reduced to primitive GW-invariants.}
\vskip 0.1in
In general, we conjecture that all the Witten invariants
can be derived from primitive GW-invariants.

Our invariant can be also applied to studying topology of symplectic
manifolds. As an example, we will verify the {\it Stabilizing conjecture}
of the first author in the case of simply connected elliptic surfaces.
The conjecture claims: {\it Suppose that $X$, $Y$ are two simply connected
homeomorphic symplectic 4-manifolds. Then $X$, $Y$
are  diffeomorphic if and only if
$X\times S^2$, $Y\times S^2$ are deformation equivalent
as symplectic manifolds.} He also verified this conjecture
for certain complex surfaces (cf. [R1]). By calculating our invariants for
the product of simply connected
elliptic surfaces with $S^2$, we will prove that
\vskip 0.1in
\noindent
{\bf Theorem F: (Theorem 5.1)} {\it The stabilizing
conjecture holds for simply connected elliptic
surfaces.}
\vskip 0.1in
It has been an interesting question in symplectic topology to find how many
different deformation classes of symplectic structures with the same
 tamed almost complex structures (up to a homotopy) could exist on a fixed
smooth
manifold. In \cite{R1}, for any positive integer $n$, the first author
constructed examples admitting at least  $n$-many
different deformation classes. Using our calculation of
GW-invariant, we can produce  examples with infinitely many deformation
classes of symplectic structures.
\vskip 0.1in
\noindent
{\bf Proposition H: (Proposition 5.4) }{\it Let $X$ be the blow-up of a simply
connected elliptic
surface at one point. Then, the smooth 6-manifold $X\times S^2$ admits
infinitely many
deformation classes of symplectic structures with the same tamed almost complex
structure up to a homotopy.}
\vskip 0.1in
This paper is organized as follows. We will define the invariants and state
the basic properties (including composition law) of our invariants in section
 2. The section 3 is a technical section where we will prove the various
results about the compactification and transversality. All the results
stated in section 2 will be proved in section 4.
 We will discuss the applications to the  stabilizing conjecture and the
Witten conjecture  in section 5, 6.

Some of the results in this paper have been lectured by us in last few
years. Also, the main results of this paper were announced
in the paper [T] of the second author published in the proceeding
of the first "Current developments in Mathematics", Boston, May, 1995.
All the basic techniques were developed in [RT].

The first author wish to thank S. K. Donaldson
who suggested the example of Section 5 to him.

\section{Higher genus symplectic invariants and composition law}

In this section, we construct the higher genus symplectic invariants. Its
physical counterparts are the correlation functions of topological sigma model
coupled with gravity. Some important cases of
these invariants were
first studied
for symplectic 4-manifolds in \cite{R3} and also discussed
in [RT]. The construction here is similar
to that of \cite{RT}.

First of all, let's introduce the  inhomogeneous Cauchy-Riemann equation, which
plays a central role in \cite{RT}. Compared to that of
\cite{RT},  we
would like to define the inhomogeneous term varying continuously as we vary
the complex structures of the Riemann surfaces. This makes the construction
more complicated.
Let $(V, \omega)$ be a  symplectic manifold and $J$ be a
tamed almost complex structure. Let $\M_{g, k}$ be the moduli space
of genus $g$ Riemann surfaces with $k$-marked points and $\overline{\M}_{g, k}$
be the Deligne-Mumford compactification. Suppose that
$$\pi: \overline{\U}_{g,k}\rightarrow\overline{\M}_{g,k}$$
is the universal curve. Both $\overline{\U}_{g,k}$ and $\overline{\M}_{g,k}$
are projective varieties. Unfortunately, it is well-known that $\overline{\U}_{
g,k}$ is not a universal family. Namely, if $\Sigma\in \overline{\M}_{g,k}$ has
 a nontrivial automorphism, then $\pi^{-1}(\Sigma)=\Sigma/Aut(\Sigma)$ instead
 of $\Sigma$. We can not directly define the inhomogeneous term over
$\overline{\U}_{g,k}$. However, this problem can be overcome by constructing
some finite covers  of $\overline{\M}_{g,k}$.
\vskip 0.1in
\noindent
{\bf Definition 2.1: }{\it A finite connected  cover $p_{\mu}:\overline{\M}^{
\mu}_{g,k}\rightarrow
\overline{\M}_{g,k}$ is a good cover if $\overline{\M}^{\mu}_{g,k}$ is a normal
projective variety with
 quotient singularity such that there is universal family
$$\pi_{\M}:\overline{\U}^{\mu}_{g,k}\rightarrow \overline{\M}^{\mu}_{g,k},$$
i.e., for each $b\in \overline{\M}^{\mu}_{g,k}, \pi^{-1}_{\M}(b)$ is a stable
Riemann surface isomorphic to $p_{\mu}(b)$. Furthermore, we have following
commutative diagram
$$\begin{array}{lccc}
   p_{\mu}& \overline{\U}^{\mu}_{g,k}&\rightarrow & \overline{\U}_{g,k}\\
         & \downarrow \pi_{\M}  && \downarrow\pi\\
p_{\mu}: & \overline{\M}^{\mu}_{g,k}& \rightarrow &\overline{\M}_{g,k}
\end{array}$$
It is clear that $\overline{\U}^{\mu}_{g,k}$ is projective and  unique.}
\vskip 0.1in
To simplify the notation, we will not distinguish $b$ with $p_{\mu}(b)$ without
any confusion. As we mentioned, the problem for $\overline{\M}_{g,k}$ is that
some elements have nontrivial automorphism groups. One can resolve this problem
by taking finite cover locally. Hence, locally one always has a universal
family. Mumford proved \cite{Mu} that such local covers can be glued together
to form a global finite cover. Moreover, it can be
explicit constructed via level $m$-structure. We will not give any detail of
level $m$-structure. The idea is to fix a basis of $H^1$ with
the coefficient in $\Z_m$. Then, there is no automorphism of stable Riemann
surface preserving the fixed basis. We refer
the reader to \cite{Mu} for the detail.

Let $p_{\mu}: \overline{\M}^{\mu}_{g,k}\rightarrow \overline{\M}_{g,k}$ be a
good
cover. Suppose that
$$\phi_{\mu}:\overline{\U}^{\mu}_{g,k}\rightarrow \P^{N}$$
is a projective embedding.
 There are two
relative tangent bundles over $\P^{N} \times V$ with respect to $h_i$ ($i=1,2$)
, where $h_i$ is the projection from $\P^{N}\times V$ to
its $i$-th factor. A
section $\nu$ of $Hom(h^*_1 T\P^{N}, h^*_2 TV)$ is said to be anti-$J$-linear
if for any tangent vector $v$ in $T\P^N$,
$$ \nu(j_{\P^N}(v)) ~=~ - J (\nu(v))
\leqno (2.1)$$
where $j_{\P^N} $ is the  complex structure on $\P^N$.
Usually, we call such a $\nu$ an inhomogeneous term. We will often drop $h_i$
from the notation without any confusion.
\vskip 0.1in
\noindent
{\bf Definition 2.2.} {\it Let $\nu$ be an inhomogeneous term.
A $(J,\nu)$-perturbed holomorphic map, or simply a $(J,\nu)$-map, is
a smooth map $f:\Sigma\rightarrow V$ satisfying the inhomogeneous
Cauchy-Riemann
equation
$$(\bar{\partial}_J f)(x)=
\nu(\phi_{\mu}(x), f(x)),
\leqno (2.2)$$
where $\bar{\partial}_J$ denotes the differential operator
$d + J\cdot d\cdot j_\Sigma$.}
\vskip 0.1in
Let $\M_{g,k,\kappa}$ be the subset of $\M_{g,k}$ with automorphism group
$\kappa$. We will use $I$ to denote the trivial automorphism group.
Since $\M_{g,k,I}$ is smooth, without the loss of generality, we can assume
that $\M^{\mu}_{g,k,I}$ is smooth, where $\M^{\mu}_{g,k,
\kappa}=p^{-1}_{\mu}(\M_{g,k, \kappa})$. Let $\U^{\mu }_{
g,k, I}$ be the preimage of $\M^{\mu}_{g,k, I}$.
 We denote by $\M_A(\mu, g,k, J, \nu)_I$ the moduli space of
 $(J, \nu)$-perturbed
holomorphic maps from $(\Sigma, x_1, \cdots, x_k)\in \M^{\mu }_{g,k,I}$  into
$V$. There
are some important topological properties as follows.

 Let
$\pi: \U_A(\mu, g,k, J, \nu)_I\rightarrow \M_A(\mu, g,k,  J, \nu)_I$ be the
universal family of
curves, i.e.,
$$\pi^{-1}(f, \Sigma, \{x_i\})=\Sigma.$$
We can define the evaluation
map
$$e_A(g,k): \U_A(\mu, g,k, J, \nu)_I\rightarrow V$$
by
$$e_A(g,k)(f, \Sigma, \{x_i\}, y)=f(y).\leqno(2.3)$$
Each marked point $x_i$ defines a section
$$\sigma_i: \M_A(\mu,g,k, J,\nu)_I\rightarrow  \U_A(\mu, g,k, J, \nu)_I$$
by
$$\sigma_i((f,\Sigma, \{x_i\}))=x_i.\leqno(2.4)$$
The composition
$$e_i=e_A(g,k)\circ \sigma_i: \M_A(\mu, g,k, J,\nu)_I\rightarrow V.$$
Let
$$\Xi^A_{g,k}=\prod^{k}_{i=1}e_A(g,k)\circ \sigma_i: \M_A(g,k,  J,\nu)_I
\rightarrow V^k.$$
Evidently, we have a map $\Upsilon_A: \M_A(\mu, g,k,  J, \nu)_I\rightarrow
\M^{\mu }_{g,k,I}$
by assigning each $(J, \nu)$-map to its domain. Together, we get a smooth map
$$\Upsilon_A \times \Xi^A_{g,k}:\M_A(\mu, g,k,  J, \nu)_I\rightarrow \M^{\mu
}_{g,k,I}\times V^k.$$
In general, $\M_A(\mu, g,k,  J, \nu)_I$ is not compact. However, there
is a natural compactification $\overline{\M}_A(\mu, g,k, J, \nu)_I$, which we
call
GU-compactification (cf. Section 3). For our purpose, we also need to consider
certain quotient $\overline{\M}^r_A(\mu, g,k, J, \nu)_I$ of $\overline{\M}_A(
\mu, g,k, J, \nu)_I$.
\vskip 0.1in
\noindent
{\bf Proposition 2.3.} {\it Suppose that $(V, \omega)$ is a semi-positive
symplectic manifold. Then, there is  a Baire set of second category-$\H$ among
 all the smooth pairs $(J, \nu)$ such that for any $(J, \nu)\in \H$
\vskip 0.1in
\noindent
(1) $\M_A(\mu, g,k,J, \nu)_I$ is a smooth, oriented manifold of real dimension
$$2c_1(V)(A)+2(3-n)(g-1)+2k;$$
\vskip 0.1in
\noindent
(2) $\Upsilon_A$ and $\Xi^A_{g,k}$ extends to continuous maps, still
denoted by the same symbols, from $\overline{\M}^r_A
(\mu, g,k, J, \nu)_I$ to $\overline{\M}^{\mu}_{g,k}$ and $V^k$,
respectively;
\vskip 0.1in
\noindent
(3) The boundary $\Upsilon_A\times \Xi^A_{g,k} \left (
\overline{\M}^r_A(\mu, g,k,J,\nu)_I\backslash \M_A(\mu,g,k,J, \nu)_I\right )$
is of real
codimension at least two.

We call that $(J, \nu)$ is generic if Proposition 2.3 is satisfied.}
\vskip 0.1in
The proof of this proposition is the main topic of Section 3.

One can construct natural cohomology classes over $\overline
 \M_A(\mu, g,k,  J, \nu)_I$ by pulling
back of the cohomology classes of $\overline{\M}^{\mu}_{g,k}\times V^k$ through
 $\Upsilon _A
\times \Xi^A_{g,k}$. Then, our invariants can be defined as the paring of the
cup products of those natural cohomology classes against the fundamental class
  of $\overline{\M}_A(\mu, g,k,  J, \nu)_I$. The existence of such a
fundamental class is a much more difficult problem, which will be discussed in
\cite{RT1}. Here, we choose to avoid this problem by considering the
intersection theory as we did in \cite{RT}.

Let $\{\alpha_i\}_{1\leq i\leq k}$ be integral homology classes of $V$. Each
$\alpha_i$ can be represented by a so called pseudo-submanifold $(P,f)$. A
pseudo-submanifold is a pair $(P, f)$, where $P$ is a finite
simplicial complex of dimension
$d_i=deg(\alpha_i)$
such that $P^{top}=P-P_{d_i-2}$ ($P_{d_i-2}-(d_i-2)$
skeleton)  is a smooth manifold and $f: P\rightarrow V$ is piecewise linear
with respect to a triangulation of $V$ and  smooth over $P^{top}$ in the usual
sense. Any two such pseudo-submanifolds representing the same homology class
are the boundary of a pseudo-submanifold cobordism in the usual sense. We refer
to the Section 4 for details. We choose pseudo-submanifolds $(Y_i, F_i)$ to
represent
$\alpha_i$. Let
$$Y=\prod^k_{i=1} Y_i, F=\prod^k_{i=1} F_i.$$
We define $Y^{top}=\prod^k_{i=1} Y^{top}_i.$
Clearly, $(Y, F)$ represents $\prod^k_{i=1}\alpha_i\in H_*(V^k, \Z)$.

In general, not every
integral homology class of $\overline{\M}^{\mu}_{g,k}$ can be represented by a
pseudo-submanifold. However, some of its multiples does.
Therefore, the homology
classes represented by pseudo-submanifolds generate the rational homology
$H_*(\overline{\M}^{\mu}_{g,k}, \Q)$. We say
that  a pseudo-manifold $(G, K)$ is {\em in the general position} if
$K(G^{top})\subset \M^{\mu }_{g,k,I}$ has codimension
at least two in $K(G^{top})$.
Let $(G, K)$ be such a pseudo-submanifold in
$\M^{\mu}_{g,k,I}$. Suppose
$$\sum^k_{i=1} (2n-d_i)+(6g-6+2k-deg(G))=2c_1(V)(A)+2(3-n)(g-1)+2k.
\leqno(2.5)$$
Then, we can choose a small perturbation of $F_i, K$ such
that $K\times F$ is transverse to $\Upsilon_A\times \Xi^A_{g,k}$ as the PL-maps
with respect to some triangulation of $V$ and as the smooth maps over $Y^{top}
\times G^{top}$. By the dimension counting, we can show that
$$\begin{array}{rl}
&Im(K\times F )\cap Im(\Upsilon_A\times \Xi^A_{g,k}
(\overline{\M}_A(\mu,g,k, J, \nu)_I-
\M_A(\mu,g,k, J, \nu)_I))=\emptyset;\\
&K\times F (Y\times G -Y^{top}\times G^{top})\cap
\Upsilon_A\times \Xi^A_{g,k} =\emptyset.\\
\end{array}
\leqno(2.6)$$
Let $\Delta\subset (\M^{\mu}_{g,k,I}\times V^k) \times (\M^{\mu }_{g,k,I}\times
V^k)$ be the diagonal. Then,
$$(\Upsilon_A\times \Xi^A_{g,k}\times K\times F)^{-1}(\Delta)\subset
\M_A(\mu,g,k, J, \nu)_I \times Y^{top}\times G^{top} \leqno(2.7)$$
is a zero-dimensional smooth submanifold. Following from (2.6), it is compact
and hence
finite. Suppose that
$$(\Upsilon_A\times \Xi^A_{g,k}\times K\times F)^{-1}(\Delta)=\{(f_1, s_1),
\cdots,
(f_m, s_m)\} \leqno(2.8)$$
where $f_i\in \M_A(\mu,g,k,J,\nu)_I$ and $s_i$ represents other factors. For
each
$(f_i, s_i)$,  we define a number $\epsilon(f_i, s_i)=\pm 1$ as follows: We
define $\epsilon(f_i, s_i)=+1$ if the orientation, induced by the Jacobian of
$\Upsilon_A\times \Xi^A_{g,k}\times K\times F$ and the orientation of
$\M_A(\mu,g,k, J, \nu)_I \times Y^{top}\times G^{top}$ at $(f_i, s_i)$,
together with the orientation
of $\Delta$ matches the  orientation of $(\M^{\mu }_{g,k,I}\times V^k)
 \times (\M^{\mu}_{g,k,I}\times V^k)$. Otherwise, we define $\epsilon(f_i,
s_i)=-1$. Now, we define
$$\Psi^V_{(A, g, k,\mu)}(K; \alpha_1, \cdots, \alpha_k)=\sum^m_{i=1}
\epsilon(f_i,s_i).$$
To Justify our notation, we will show
\vskip 0.1in
\noindent
{\bf Proposition 2.4.} {\it
\vskip 0.in
\noindent
(1) $\Psi^V_{(A, g, k,\mu)}(K; \alpha_1, \cdots, \alpha_k)$
is independent of $J, \nu$, the pseudo-submanifold representatives
$(Y_i, F_i)$.
\vskip 0.1in
\noindent
(2) $\Psi^V_{(A, g, k,\mu)}(K; \alpha_1, \cdots, \alpha_k)$ is independent of
semi-positive deformations of $\omega$.
\vskip 0.1in
\noindent
(3) If $(G', K')$ is another pseudo-submanifold which is in the general
position
and represents the same homology class as that of $(G, K)$,
$$\Psi^V_{(A, g, k,\mu)}(K; \alpha_1, \cdots, \alpha_k)=\Psi^V_{(A, g, k,\mu)}
(K'; \alpha_1, \cdots, \alpha_k).$$}
\vskip 0.1in

Therefore, we can write $\Psi^V_{(A, g, k,\mu)}([K]; \alpha_1, \cdots,
\alpha_k)$
for $\Psi^V_{(A, g, k,\mu)}(K; \alpha_1, \cdots, \alpha_k)$, where
$[K]$ denotes the homology
class represented by the cycle $K$.

We postpone the proof of Proposition 2.4 to Section 4.
\vskip 0.1in
\noindent
{\bf Proposition 2.5. }{\it $\Psi^V_{(A, g, k,\mu)}(K; \alpha_1, \cdots,
\alpha_k)$ is independent of the embedding $\phi_{\mu}$. Hence, $\Psi_{(A,
g,k, \mu)}$ is a symplectic invariant}
\vskip 0.1in
\noindent
{\bf Proof: } Suppose that
$$\tilde{\phi}_{\mu}: \overline{\U}_{\mu}\rightarrow \P^{N'}$$
is a different projective embedding. Then,
$$\phi_{\mu}\times \tilde{\phi}_{\mu}:  \overline{\U}_{\mu}\rightarrow
\P^{N}\times \P^{N'}.$$
One can  consider the inhomogeneous term $\bar{\nu}\in \overline{Hom}_J(T(
\P^N\times  \P^{N'}), TV)$ where $
\overline{Hom}_J$ means anti-complex linear homomorphism. Moreover, one can
use $(J, \bar{\nu})$ to define invariant in the same fashion and prove that the
invariant is independent of $(J, \bar{\nu})$ using the same proof of
Proposition 2.4. Let $\nu$, $\tilde{\nu}$ be the inhomogeneous term defined
through the embedding
$\phi_{\mu}$ and $\tilde{\phi}_{\mu}$. Notes that
$$T(\P^N\times \P^{N})= T(\P^N)\times T(\P^{N'}).$$
Therefore, we can view both $\nu$ and $\tilde{\nu}$ as  the sections of
$\overline{Hom}_J(T(\P^N\times \P^{N'}), TV)$, where we view that $\nu$ maps
the factor $T(\P^{N'})$ to zero  and $\tilde{\nu}$ maps the first factor to
zero. Let's denote them
by $\bar{\nu}$ and $\bar{\tilde{\nu}}$. We observe that if $\nu(\tilde{
\nu})$ is generic, so is $\bar{\nu}(\bar{\tilde{\nu}})$. It follows from the
 definition that
we have the same invariant $\Psi$ using $(J, \nu)$ or $(J, \bar{\nu})$. In
the same way, we have the same invariant using $(J, \tilde{\nu})$ or
 $(J, \bar{\tilde{\nu}})$.
As we mentioned, by  repeating the proof of Proposition 2.4, we can show
that the invariant defined by $(J, \bar{\nu})$ is the same as the invariant
defined by $(J, \bar{\tilde{\nu}})$. Then, we finish the proof.
\vskip 0.1in
\noindent
{\bf Remark 2.6: }{\it We don't have to restrict ourself to projective
embedding. In
fact, we can embed $\overline{\U}^{\mu}_{g,k}$ into any smooth complex manifold
and define the inhomogeneous term in the same fashion. The proof of Proposition
2.5 shows that the resulting invariant is independent of such an embedding.}
\vskip 0.1in
For the convenience, we define $\Psi^V_{(A, g, k, \mu)}(K, \alpha_1, \cdots,
\alpha_k)=0$ if (2.5) is not satisfied.

Let's collect some properties of $\Psi^V_{(A, g, k,\mu)}.$ The
following proposition essentially
follows from the definition. We will omit its proof.
\vskip 0.1in
\noindent
{\bf Proposition 2.7.}{\it
\vskip 0.1in
\noindent
(1) $\Psi^V_{(A, g, k,\mu)}=0$, if either $\omega(A)< 0$ or $c_1(V)(A)+(3-n)(g
-1)<0$, in particular, $\Psi_{(0,g,k,\mu)}=0$ for any $g\geq 2$ and $n\geq 4$.
\vskip 0.1in
\noindent
(2) $\Psi^V_{(A, g, k,\mu)}$ is multilinear and supersymmetry on $H_*(V, \Z)^k$
with respect to the $\Z_2$-grading by even and odd degrees.}
\vskip 0.1in
We have established
the symplectic invariant
$$ \Psi ^V_{(A,g,k,\mu)} : H_*(\overline{\M}^{\mu}_{g,k} , \Q)
\times H_*(V, \Z)^k\mapsto \Q,$$
for any $A\in H_2(V, \Z)$ and $2g +k \ge 3$.

Now we discuss other more interesting properties of $\Psi^V_{(A, g, k,\mu)}$
associated with the structures of $\overline{\M}_{g,k}$.

As we mentioned in the introduction, there is a natural map
$$\pi: \overline{\M}_{g,k}\rightarrow \overline{\M}_{g, k-1} \leqno(2.9)$$
by forgetting the last
marked point and contracting the unstable rational component.
 One should be aware that
there are two exceptional cases $(g,k)=(0,3), (1,1)$ where $\pi$ is not well
defined. Suppose $\overline{\M}^{\mu}_{g,k}\rightarrow \overline{\M}_{g,k}$ is
a good
cover constructing through level-$m$ structure. Then, one can observe that
$\overline{\M}^{\mu}_{g,k+1}=\pi^*\overline{\M}^{\mu}_{g,k}$ is a good cover of
$\overline{\M}_{g, k+1}$.
Let
$$\pi_{\mu}: \overline{\M}^{\mu}_{g, k+1}\rightarrow
\overline{\M}^{\mu}_{g,k}.$$
 Then, $\pi_{\mu}$ induces a map on the universal families
(still denoted by $\pi_{\mu}$.)
$$\begin{array}{clc}
\overline{\U}^{\mu}_{g,k+1}&\stackrel{\pi_{\mu}}{\rightarrow}&\overline{\U}^{\mu}_{g,k}\\
\downarrow \pi_{\M}&                              & \downarrow \pi_{\M}\\
\overline{\M}^{\mu}_{g,k+1}&\stackrel{\pi_{\mu}}{\rightarrow} &
\overline{\M}^{\mu}_{g, k}
\end{array}$$
Let $b\in \overline{\M}^{\mu}_{g, k+1}$ and $\Sigma_b$ be the underline stable
Riemann
surface.
 Clearly,
$$\pi_{\mu}: \Sigma_b\rightarrow \Sigma_{\pi_{\mu}(b)}$$
is precisely $\pi$ defined in (2.9).

Associated with $\pi$,
we have two {\em k-reduction formulas} for $\Psi^V_{(A, g, k, \mu)}$.
\vskip 0.1in
\noindent
{\bf Proposition 2.8. }{\it Suppose that $(g,k)\neq (0,3),(1,1)$. Furthermore,
suppose that $\overline{\M}^{\mu}_{g,k+1}$ and $\overline{\M}^{\mu}_{g,k}$ are
defined as above.
\vskip 0.1in
\noindent
(1) For any $\alpha _1, \cdots , \alpha _{k-1}$ in $H_*(V, \Z)$,
we have}
$$\Psi ^V_{(A,g,k,\mu)}([K]; \alpha _1, \cdots,\alpha _{k-1}, [V])~=~
\Psi ^V_{(A,g,k-1,\mu)}((\pi_{\mu})_*[K]; \alpha _1, \cdots,\alpha _{k-1})
\leqno (2.10)$$
\vskip 0.1in
\noindent
(2)  Let $\alpha _k$ be in $H_{2n-2}(V, \Z)$, then
$$\Psi ^V_{(A,g,k,\mu)}([K]; \alpha _1,
\cdots,\alpha _{k-1}, \alpha _k)~=~\alpha^* _k (A)
\Psi ^V_{(A,g,k-1,\mu)}([(\pi_{\mu})^{-1}(K)]; \alpha _1, \cdots,\alpha _{k-1})
\leqno (2.11)$$
where $ \alpha^* _k$ is the Poincare dual of $\alpha _k$.
\vskip 0.1in
\noindent
{\bf Proof:} The proof is similar to that of Proposition 2.5. Let
$$\phi^k_{\mu}: \overline{\U}^{\mu}_{g,k}\rightarrow \P^N; \ \phi^{k-1}_{\mu}:
\overline{\U}^{\mu}_{g, k-1}\rightarrow \P^{N'}$$
be projective embedding and $\nu_k, \nu_{k-1}$ be inhomogeneous terms over
$\P^N$ or $\P^{N'}$ respectively. Consider embedding
$$\phi^k_{\mu}\times (\phi^{k-1}_{\mu}\circ \pi_{\mu}):
\overline{\U}^{\mu}_{g,k}\rightarrow \P^N\times \P^{N'}. \leqno(2.12)$$
As in Remark 2.6, we can define the invariant $\bar{\Psi}$ using
$\phi^k_{\mu}\times
(\phi^{k-1}_{\mu}\circ \pi_{\mu})$ and the inhomogeneous terms over $\P^N\times
\P^{N'}$ in the same fashion. One can show that such an invariant is
independent of inhomogeneous term. Furthermore, we can view $\nu_k, \nu_{k-1}$
(denoted by $\bar{\nu}_k, \bar{\nu}_{k-1}$) as inhomogeneous terms over
$\P^N\times \P^{N'}$ as in the proof of Proposition 2.5. Clearly, using
$(J, \bar{\nu}_k)$, we have
$$\bar{\Psi}^V_{(A,g,k,\mu)}([K]; \alpha _1, \cdots,\alpha _{k-1}, [V])~=~
\Psi^V_{(A,g,k,\mu)}([K]; \alpha _1, \cdots,\alpha _{k-1}, [V]).\leqno(2.13)$$
Using $(J, \bar{\nu}_{k-1})$, we claim that
$$\bar{\Psi}^V_{(A,g,k,\mu)}([K]; \alpha _1, \cdots,\alpha _{k-1}, [V])~=~
\Psi^V_{(A,g,k-1,\mu)}(\pi_*[K]; \alpha _1, \cdots,\alpha
_{k-1}).\leqno(2.14)$$
Suppose that $(J, \nu_{k-1})$ is generic. We claim that
$(J, \bar{\nu}_{k-1})$ satisfies the Proposition 2.2.

The Proposition 2.3 (1)
and (2) are obvious. A remark is required for (3). In the proof of (3) in the
next section, the idea is to stratify $\overline{\M}^r_A(\mu, g, k,  J, \nu)_I-
\M_A(\mu, g, k, J, \nu)_I$ and show that each strata is of real codimension
at least 2. In the proof, we use the fact that $\nu$ is generic over each
component of stable Riemann surface called the principal components. Note that
$\bar{\nu}_{k-1}$ is zero over the rational components $\pi_{\mu}$ contracts.
So it is not generic. However, we can simply treat these rational components as
bubble components (see Definition 3.6), and construct another quotient space
$$\overline{\M}^r_A\rightarrow \overline{\M}^{rr}_A$$
in the same way as we construct $\overline{\M}^r_A$. Then, the proof of
Proposition 2.3 shows that
$$\overline{\M}^{rr}_A(\mu, g, k, J, \nu)_I-\M_A(\mu, g, k,  J, \nu)_I
\leqno(2.15)$$
is of real codimension at least 2.

Once $(J, \bar{\nu}_{k-1})$ satisfies the Proposition 2.3, we can choose
$(Y_i, F_i), (G, K)$ to satisfy (2.6), (2.7). Then,
$$\bar{\Psi}_{(A, g,k,\mu)}(K; \alpha_1, \cdots, \alpha_k)=\sum^m_{i=1}
\epsilon(f_i, s_i),
$$
where $(f_i, s_i)$ is given in (2.7). Suppose that $s_i=(\Sigma, x_1, \cdots,
x_k,  y_1, \cdots, y_{k-1}, y_k).$ Then,
$$(\Sigma, x_1, \cdots,x_k)\in \M^{\mu}_{g,k,I}$$
and
$$\pi(\Sigma, x_1, \cdots, x_k)=(\Sigma, x_1, \cdots, x_{k-1}).$$
Clearly, $f\in \M_A(\mu, g,k, J, \bar{\nu}_{k-1})_I$ can also be viewed as an
element of
$\M_A(\mu, g,k-1, J, \nu)_I$.
If $Y_k=V, F_k=Id$, let
$$\bar{s}_i=(\Sigma, x_1, \cdots, x_{k-1},  y_1, \cdots, y_{k-1}).$$
Clearly
$$(\Upsilon_A \times \Xi_{g,k-1}\times \pi_{\mu}(K)\times
\prod^{k-1}_{i=1}F_i)^{-1}(\Delta)=
\{(f_1, \bar{s}_1), \cdots, (f_m, \bar{s}_m)\}. \leqno(2.16)$$
Furthermore, it is easy to check that
$$\epsilon(f_i, s_i)=\epsilon(f_i, \bar{s}_i).$$
Therefore,
$$\bar{\Psi} ^V_{(A,g,k,\mu)}([K]; \alpha _1, \cdots,\alpha _{k-1}, [V])~=~
\Psi ^V_{(A,g,k-1,\mu)}((\pi _{\mu})_*[K]; \alpha _1, \cdots,\alpha _{k-1}).$$
Here, we use the fact that $(G, \pi_{\mu}\circ K)$ represents the homology
class $(\pi _{\mu})_*[K]$.

The proof of (2) is similar.

Next, we discuss the composition law. Roughly speaking, the composition law
governs the change of $\Psi$ under the surgery of Riemann surfaces.
Compared to $k$-reduction formula, one
can view the composition law as $g$-reduction formula,
i.e., the reduction of
genus. As we mentioned in the introduction, its classical cousin in
enumerative algebraic geometry is the
degeneration formula, which was only derived individually in very special
cases. One technical reason is that it is very difficult to have a good
deformation theory in algebraic
geometry. But our $\Psi$ is just a symplectic invariant. The counterpart of
deformation theory in symplectic category can be realized by
a {\em gluing theorem}. It
has been established by the authors in \cite{RT}.

Recall that in the definition of $\Psi_{(A,g,k, \mu)}([K]; \alpha_1, \cdots,
\alpha_k)$, we require that $K$ does not
lie in the boundary of $\overline{\M}^{\mu}_{g,k}$.
We first remove this technical assumption.
\vskip 0.1in
\noindent
{\bf Proposition 2.9.} {\it $\Psi_{(A,g,k, \mu)}(K; \alpha_1,
 \cdots, \alpha_k)$ can be defined when $K$ is in the boundary
of $\overline{\M}^{\mu}_{g,k}$. Namely, $K\subset Im \theta_S$
or $K\subset Im \bar{\mu}$, where $\theta_S, \bar{\mu}$ are defined below.
 Furthermore, $\Psi_{(A,g,k, \mu)}(K; \alpha_1,
 \cdots, \alpha_k)$ depends only on the homology class represented by $(G, K)$.
Then, we extend $\Psi_{(A, g,k)}([K]; \cdots)$ for any $[K]\in H_*(\overline{
\M}^{\mu}_{g,k}, \Q)$ by the linearity.}
\vskip 0.1in
The proof follows basically from the gluing theorem in \cite{RT},
Theorem 6.1. The details will appear in Section 4.

Assume $g=g_1+g_2$ and $k=k_1+k_2$ with $2g_i + k_i \ge 3$.
Fix  a decomposition $S=S_1\cup S_2$ of $\{1,\cdots , k\}$ with
$|S_i|= k_i$. Recall that
$\theta _S: \overline \M_{g_1,k_1+1}\times \overline \M_{g_2,k_2+1}
\mapsto \overline \M_{g,k}$, which assigns to marked
curves $(\Sigma _i; x_1^i,\cdots ,x_{k_1+1}^i)$ ($i=1,2$), their
union $\Sigma _1\cup \Sigma _2$ with $x^1_{k_1+1}$ identified to
$x^2_1$ and remaining points renumbered by $\{1,\cdots,k\}$ according to $S$.
Suppose that $\overline{\M}^{\mu}_{g_1, k_1+1}, \overline{\M}^{\mu}_{g_2,
k_2+1}, \overline{\M}^{\mu}_{g,k}$ are
 good covers over $\overline \M_{g_1,k_1+1},  \overline \M_{g_2,k_2+1},
 \overline \M_{g,k}$ such that
$$\theta^*_S\overline{\M}^{\mu}_{g,k}= \overline{\M}^{\mu}_{g_1, k_1+1}\times
\overline{\M}^{\mu}_{g_2, k_2+1}.
\leqno(2.17)$$
Such good covers can be constructed using the level-$n$ structure of $\overline
{\M}_{g,k}$.
We have another natural map defined in the introduction
$\mu : \overline \M_{g-1, k+2}
\mapsto \overline \M_{g,k}$ by gluing together the last two marked
points.
Let
$$\overline{\M}^{m_{\mu}}_{g-1, k+2}=\mu^*\overline{\M}^{\mu}_{g,k}, \
\bar{\mu}: \overline{\M}^{m_{\mu}}_{g-1, k+2} \rightarrow
\overline{\M}^{\mu}_{g,k}.\leqno(2.18)$$

Choose a homogeneous
basis $\{\beta _b\}_{1\le b\le L}$ of $H_*(V,\Z)$ modulo
torsion. Let $(\eta _{ab})$ be its intersection matrix. Note that
$\eta _{ab} = \beta _a \cdot \beta _b =0$ if the dimensions of
$\beta _a$ and $\beta _b$ are not complementary to each other.
Put $(\eta ^{ab})$ to be the inverse of $(\eta _{ab})$.
Now we can state the composition law, which consists of two formulas.
\vskip 0.1in
\noindent
{\bf Theorem 2.10 }. {\it Let $[K_i] \in H_*(\overline \M^{\mu}_{g_i,
k_i+1}, \Q)$ $(i=1,2)$ and $[K_0] \in H_*(\overline \M^{m_{\mu}}_{g-1,
k +2}, \Q)$. Suppose that $\overline{\M}^{\mu}_{g_1, k_1+1}, \overline{\M}^{
\mu}_{g_2, k_2+1},
\overline{\M}^{\mu}_{g,k}, \overline{\M}^{m_{\mu}}_{g-1, k+2}$ are defined as
(2.17), (2.18). For any $\alpha _1,
\cdots,\alpha _k$ in $H_*(V,\Z)$.
Then we have}
$$\begin{array}{rl}
&\Psi ^V_{(A,g,k, \mu)}(\theta _{S*}[K_1\times K_2];\{\alpha _i\})\\
=~& \sum \limits _{A=A_1+A_2} \sum \limits_{a,b}
\Psi ^V_{(A_1,g_1,k_1+1, \mu)}([K_1];\{\alpha _{i}\}_{i\le k}, \beta _a)
\eta ^{ab}
\Psi ^V_{(A_2,g_2,k_2+1, \mu)}([K_2];\beta _b,
\{\alpha _{j}\}_{j> k}) \\
\end{array}
\leqno (2.19)
$$
$$
\Psi ^V_{(A,g,k, \mu)}(\bar{\mu}_*[K_0];\alpha _1,\cdots, \alpha _k)
=\sum _{a,b} \Psi ^V_{(A,g-1,k+2, m_{\mu})}([ K_0];\alpha _1,\cdots, \alpha _k,
\beta _a,\beta _b) \eta ^{ab}\leqno (2.20)
$$
\vskip 0.1in

The proof of composition law essentially follows from Proposition 2.9 by some
topological arguments. We
postpone it to Section 4.

So far, we are working on the covers $\overline{\M}^{\mu}_{g,k}$. To define
invariants over $\overline{\M}_{g,k}$, we introduce following classical
notion in algebraic topology to relate homology of $\overline{\M}^{\mu}_{g,k}$
to homology of $\overline{\M}_{g,k}$.
\vskip 0.1in
\noindent
{\bf Definition 2.11. }{\it Suppose that $f: M\rightarrow N$ be a continuous
map such that both $M$ and $N$ have Poincare duality. Then we define transfer
map
$$f_{!}: H_*(M)\rightarrow H_*(N)$$
by $f_{!}(\alpha)=PD_{M}\circ f^*\circ PD^{-1}_N(\alpha)$, where $PD_{M},
PD_{N}$
are Poincare duality maps. Furthermore, transfer maps compose functorially. }
\vskip 0.1in
For our case, we use $\Q$ as coefficients since Poincare duality only holds
over rational coefficient. If $M$ is a finite cover of $N$, it is easy to
observe that
$$f_*f_{!}(\alpha)=\lambda\alpha,$$
where $\lambda$ is the order of covers.

  Suppose that $\lambda^{\mu}_{g,k}$ is the order of cover for
$p_{\mu}: \overline{\M}^{\mu}_{g,k}\rightarrow \overline{\M}_{g,k}$.
\vskip 0.1in
\noindent
{\bf Definition 2.12. }{\it For any $[K]\in H_*(\overline{M}_{g,k}, \Q)$ and
$\{\alpha_i\}\in H_*(V, \Z)$, define
$$\Psi_{(A,g,k)}([K];
\{\alpha_i\})=\frac{1}{\lambda^{\mu}_{g,k}}\Psi_{(A,g,k,\mu)}((p_{\mu})_{!}([K]); \{\alpha_i\}). \leqno(2.21)$$}
\vskip 0.1in
\noindent
{\bf Lemma 2.13.}{\it $\Psi_{(A,g,k)}([K]; \{\alpha_i\})$ is independent of
$\overline{\M}^{\mu}_{g,k}$.}
\vskip 0.1in
\noindent
{\bf Proof: }  Consider the fiber product
$$\overline{\M}=\overline{\M}^{\mu}_{g,k}\times_{\overline{\M}_{g,k}}\overline{\M}^{\mu'}_{
g,k}.$$
Let
$$p^1: \overline{\M}\rightarrow \overline{\M}^{\mu}_{g,k}; \  p^2: \overline{
\M}\rightarrow \overline{\M}^{\mu'}_{g,k}$$
be projections. Then, we can pull back the universal family $\overline{\U}^{
\mu}_{g,k},  \overline{\U}^{\mu'}_{g,k}$ by $p^1, p^2$. Both are obviously the
universal family over
$\overline{\M}$. By the uniqueness of universal family, they must be the same.
In other
words, let $\overline{\U}$ be the universal family over $\overline{\M}$.
$$\overline{\U}=(p^1)^*\overline{\U}^{\mu}_{g,k}=(p^2)^*\overline{\U}^{\mu'}_{g,k}. \leqno(2.22).$$
Let
$$\phi_{\mu}: \overline{\U}^{\mu}_{g,k}\rightarrow \P^N; \phi_{\mu'}:
\overline{\U}^{\mu'}_{g,k}\rightarrow \P^{N'}$$
be projective embedding.  We have a natural
embedding
$$\phi=\phi_{\mu}\times \phi_{\mu'}: \overline{\U}\rightarrow
\P^N\times\P^{N'}.
\leqno(2.23)$$
By Remark 2.6, we can use this embedding to define inhomogeneous term and
define
the analogue of $\Psi$. Furthermore, we can show that such invariant is
independent of the choice of a generic choice of inhomogeneous term. Let's
denote this invariant as $\tilde{\Psi}_{(A, g,k)}$.
We claim that
$$\tilde{\Psi}_{(A, g,k)}((p^1p_{\mu})_{!}([K]), \{\alpha_i\})=\lambda^{\mu'}_{
g,k}\Psi_{(A,g,k,\mu)}((p_{\mu})_{!}([K]),  \{\alpha_i\})\leqno(2.24)$$
and
$$\tilde{\Psi}_{(A, g,k)}((p^2p_{\mu'})_{!}([K]), \{\alpha_i\})=\lambda^{\mu}_{
g,k}\Psi_{(A,g,k,\mu')}((p_{\mu'})_{!}([K]),  \{\alpha_i\}).\leqno(2.25)$$
For any inhomogeneous term $\nu$ over $\P^N$, we can again view it as an
inhomogeneous term over $\P^N\times \P^{N'}$ and denote it by $\bar{\nu}$.
If $(J, \nu)$ is generic, so is $(J, \bar{\nu})$. Choose $(H, L)$ to represent
$(p^1p_{\mu})_{!}([K])$ and together with $Y$ satisfying (2.6), (2.7). Then,
$(H, p^1\circ L)$ represents $p^1_*((p^1p_{\mu})_{!}([K]))=\lambda^{\mu'}_{g,k}
(p_{\mu})_{!}(
[K])$ since the order of cover $p^1: \overline{\M}\rightarrow \overline{\M}^{
\mu}_{g,k}$ is
$\lambda^{\mu'}_{g,k}$. It is straightforward to check that
$$(\Upsilon_A\times \Xi^A_{g,k}\times F \times L)^{-1}(\Delta)=(\Upsilon_A
\times \Xi^A_{g,k}\times F \times p^1\circ L)^{-1}(\Delta).$$
Furthermore, the orientation matches. Then,
$$\begin{array}{rcl}
\tilde{\Psi}_{(A, g,k)}((p^1p_{\mu})_{!}([K]),
\{\alpha_i\})&=&\Psi_{(A,g,k,\mu)}(p_*(p^1p_{\mu})_{!}([K]),  \{\alpha_i\})\\
&=&\lambda^{\mu'}_{g,k}\Psi_{(A,g,k,\mu)}((p_{\mu})_{!}([K]),  \{\alpha_i\})
\end{array}.\leqno(2.26)$$
It is the same argument to show that
$$\tilde{\Psi}_{(A, g,k)}((p^2p_{\mu'})_{!}([K]), \{\alpha_i\})=\lambda^{\mu}_{
g,k}
\Psi_{(A,g,k,\mu')}((p_{\mu'})_{!}([K]),  \{\alpha_i\}).\leqno(2.27)$$
On the other hand, $p^1p_{\mu}=p^2p_{\mu'}$. Therefore,
$$\frac{1}{\lambda^{\mu}_{g,k}}\Psi_{(A,g,k,\mu)}((p_{\mu})_{!}([K]),
\{\alpha_i\})
=\frac{1}{\lambda^{\mu'}_{g,k}}\Psi_{(A,g,k,\mu')}((p_{\mu'})_{!}([K]),
\{\alpha_i\}).\leqno(2.28)$$
This finishes the proof.
\vskip 0.1in
\noindent
{\bf Proposition 2.14. }{\it
\vskip 0.1in
\noindent
(1) $\Psi^V_{(A, g, k)}(K, \alpha_1, \dots, \alpha_k)$ is
a symplectic invariant.
\vskip 0.1in
\noindent
(2) $\Psi^V_{(A, g, k)}(K, \alpha_1, \dots, \alpha_k)$ is independent of
semi-positive deformations of $\omega$.
\vskip 0.1in
\noindent
(3) $\Psi^V_{(A, g, k)}=0$, if either $\omega(A)\leq 0$ or $C_1(V)(A)+(3-n)(g
-1)<0$, in particular, $\Psi_{(0,g,k)}=0$ for any $g\geq 2$ and $n\geq 4$.
\vskip 0.1in
\noindent
(4) $\Psi^V_{(A, g, k)}=0$ is multilinear and supersymmetry on $H_*(V, \Z)^k$
with respect to the $\Z_2$-grading by even and odd degrees.}
\vskip 0.1in
The proof follows from the Proposition 2.4,
2.7 and the Definition 2.12.
\vskip 0.1in
\noindent
{\bf Proposition 2.15. }{it Suppose that  $(g,k)
\neq (0,3),(1,1)$.
\vskip 0.1in
\noindent
(1) For any $\alpha _1, \cdots , \alpha _{k-1}$ in $H_*(V, \Z)$,
we have}
$$\Psi ^V_{(A,g,k)}(K; \alpha _1, \cdots,\alpha _{k-1}, [V])~=~
\Psi ^V_{(A,g,k-1)}([\pi (K)]; \alpha _1, \cdots,\alpha _{k-1})
\leqno (2.29)$$
\vskip 0.1in
(2)  Let $\alpha _k$ be in $H_{2n-2}(V, \Z)$, then
$$\Psi ^V_{(A,g,k)}([\pi^{-1}(K)]; \alpha _1,
\cdots,\alpha _{k-1}, \alpha _k)~=~\alpha^* _k (A)
\Psi ^V_{(A,g,k-1)}(K; \alpha _1, \cdots,\alpha _{k-1})
\leqno (2.30)$$
where $ \alpha^* _k$ is the Poincare dual of $\alpha _k$.
\vskip 0.1in
\noindent
{\bf Proof: } Notes that $\overline{\M}^{\mu}_{g,
k+1}=\pi^*\overline{\M}^{\mu}_{g, k}$.
Geometrically, $(p_{\mu})_{!}([K])$ is represented by $(p_{\mu})^{-1}(K)$.
By the construction of $\M^{\mu}_{g, k+1}$,
$$\pi_{\mu}(p_{\mu})^{-1}(K)=(p_{\mu})^{-1}(\pi_{\mu}(K)). \leqno(2.31)$$
Hence,
$$(\pi_{\mu})_*(p_{\mu})_{!}([K])=(p_{\mu})_{!}(\pi_{\mu})_*([K]).
\leqno(2.32)$$
Then, (1) follows from Proposition 2.8 and (2.32).

 For (2), $[\pi^{-1}((p_{\mu})_{!}([K]))]=\pi_{!}((p_{\mu})_{!}([K]))$. By the
natality of transfer map,
$$(p_{\mu})_{!}\pi_{!}([K])=(p_{\mu}\pi)_{!}([K])=(\pi_{\mu}p_{\mu})_{!}([K])=
 (\pi_{\mu})_{!}(p_{\mu})_{!}([K]). \leqno(2.33)$$
Then, (2) follows from the Proposition 2.8(2) and (2.33).
\vskip 0.1in
\noindent
{\bf Theorem 2.16 (Composition Law).} {\it Let $[K_i] \in H_*(\overline
\M_{g_i,
k_i+1}, \Q)$ $(i=1,2)$ and $[K_0] \in H_*(\overline \M_{g-1,
k +2}, \Q)$. For any $\alpha _1,\cdots,\alpha _k$ in $H_*(V,\Z)$.
Then we have}
$$\begin{array}{rl}
&\Psi ^V_{(A,g,k)}(\theta _{S*}[K_1\times K_2];\{\alpha _i\})\\
=~& \sum \limits _{A=A_1+A_2} \sum \limits_{a,b}
\Psi ^V_{(A_1,g_1,k_1+1)}([K_1];\{\alpha _{i}\}_{i\le k}, \beta _a)
\eta ^{ab}
\Psi ^V_{(A_2,g_2,k_2+1)}([K_2];\beta _b,
\{\alpha _{j}\}_{j> k}) \\
\end{array}
\leqno (2.34)
$$
$$
\Psi ^V_{(A,g,k)}(\mu_*[K_0];\alpha _1,\cdots, \alpha _k)
=\sum _{a,b} \Psi ^V_{(A,g-1,k+2)}([ K_0];\alpha _1,\cdots, \alpha _k,
\beta _a,\beta _b) \eta ^{ab}\leqno (2.35)
$$
\vskip 0.1in
\noindent
{\bf Proof: } Let
$$p^{\mu}_{g, k}: \overline{\M}^{\mu}_{g,k}\rightarrow \overline{\M}_{g,k}, \
p^{m_{\mu}}_{g-1, k+2}: \overline{\M}^{m_{\mu}}_{g-1. k+2}\rightarrow
\overline{\M}_{g-1,
k+2}. \leqno(2.36)$$
By the Proposition 2.10,
$$\begin{array}{rl}
&\Psi ^V_{(A,g,k)}(\theta _{S*}[K_1\times K_2];\{\alpha _i\})\\
=~&\frac{1}{\lambda^{\mu}_{g,k}}\Psi ^V_{(A,g,k, \mu)}(\theta _{S*}(p^{\mu}_{g,
k})_{!}[K_1\times K_2];\{\alpha _i\})\\
=~&\frac{1}{\lambda^{\mu}_{g_1, k_1+1}\lambda^{\mu}_{g_2, k_2+1}}\Psi
^V_{(A,g,k, \mu)}(\theta _{S*}(p^{\mu}_{
g_1, k_1+1})_{!}([K_1])\times (p^{\mu}_{g_2, k_2+1})_{!}([K_2]);\{\alpha
_i\})\\
=~&\sum \limits _{A=A_1+A_2} \sum \limits_{a,b}
\frac{1}{\lambda^{\mu}_{g_1, k_1+1}}\Psi ^V_{(A_1,g_1,k_1+1,
\mu)}((p^{\mu}_{g_1,k_1+1})_{!}[K_1];\{\alpha _{i}\}_{i\le k}, \beta _a) \\
&\eta ^{ab} \frac{1}{\lambda^{\mu}_{g_2, k_2+1}}
\Psi ^V_{(A_2,g_2,k_2+1, \mu)}(p^{\mu}_{g_2, k_2+1})_{!}[K_2];\beta _b,
\{\alpha _{j}\}_{j> k}) \\
=~&\sum \limits _{A=A_1+A_2} \sum \limits_{a,b}
\Psi ^V_{(A_1,g_1,k_1+1)}([K_1];\{\alpha _{i}\}_{i\le k}, \beta _a)
\eta ^{ab}
\Psi ^V_{(A_2,g_2,k_2+1)}([K_2];\beta _b,
\{\alpha _{j}\}_{j> k})
\end{array}
\leqno(2.37)
$$
It is the same argument for (2).
\vskip 0.1in
\noindent
{\bf Remark 2.17: } The mixed invariant $\Phi _{(A,\omega ,g)}$ in [RT] can
be identified with certain $\Psi$ by choosing appropriate cycles $[K]$.
More precisely, for any $k,l\ge 0$ with $2 g + k \ge 3$, put
$K_{k,l}$ to be the closure of the cycle
$K^0_{k,l}$ in $\overline \M_{g,k+l}$, where
$K^0_{k,l}$ is the set of all
$(\Sigma , x_1,\cdots, x_{k+l})$ in $\M_{g,k+l}$
with $(\Sigma , x_1,\cdots,x_k)$ being a fixed point in $\M_{g,k}$.
Then
$$\Psi ^V_{(A,g,k+l)}(K_{k,l}; \alpha_1,\cdots ,\alpha _{k+l})
=\Phi_{(A,\omega ,g)}(\alpha _1,\cdots, \alpha _k |
\alpha _{k+1},\cdots ,\alpha _{k+l})\leqno(2.38)$$
It follows from the Proposition 2.5 that
$\Phi_{(A,\omega ,g)}(\alpha _1,\cdots |\cdots ,\alpha _{k+l})$ is $0$ if
$\dim (\alpha _{k+l}) > 2n -2$ (cf. [RT]).
\vskip 0.1in
Recall that $\Psi$ is only defined for so call ``generic $(J, \nu)\in \H$''.
Following from \cite{RT}, we can relax this genericity condition as follows:
\vskip 0.1in
\noindent
{\bf Definition 2.18. }{\it We call $(J, \nu)$ to be $A$-good if the following
two conditions are satisfied.
\vskip 0.1in
\noindent
(1). The set $\{ f\in \M_A(g,k, J,\nu)_I| Coker(Jdf \oplus D_f)\neq 0\}$ is
of real codimension 2, where $D_f$ is the linearization of the inhomogeneous
Cauchy-Riemann equation for $(J, \nu)$-maps at $f$.
\vskip 0.1in
\noindent
(2). $\overline{\M}^r_A(g,k,J,\nu)_I-\M_A(g,k, J, \nu)_I$ is of real
codimension
 2.}
\vskip 0.in
One can see from the construction that $\Psi$ is well defined if $(J,\nu)$ is
$A$-good.
\vskip 0.1in
\noindent
{\bf Remark 2.19 (Relation of $\Psi$ to enumerative function):} One of main
applications of the composition law of
the genus-0 GW-invariant $\Psi^V_{(A,0,k)}$ is
to compute the enumerative invariants of rational curves in
complex homogeneous
spaces. In \cite{RT}, Lemma 10.1, the authors proved the equivalence
of $\Psi^V_{(A,0,k)}$ with the enumerative function of
rational curves on complex Grassmannian manifolds.
More precisely, we showed
that for Grassmannian manifolds,
the integrable complex structure $J_0$ with zero perturbation term satisfies
the condition of Proposition 2.3, i.e., $(J_0, 0)\in \H$.
Therefore, one doesn't
need to deform the integrable complex structure or add any perturbation terms.
The integrable complex structure is already $A$-good for calculating the
invariant $\Psi^V_{(A,0,k)}$. Then, by the definition, $\Psi^V_{(A,0,k)}$ is
an enumerative function.

The same has been proved by
Jun Li and the second author in [LT] for any complex
homogeneous manifolds (also see \cite{CF} for an alternative proof
for Flag manifolds).

In contrast to the genus 0 case, $\Psi^V_{(A,g,k)}$($g\geq 1$) is not the same
as the enumerative function even for the projective plane $\P^2$.
For example, a
simple computation through the composition law will yield the mixed
invariant
$$\Psi^{\P^2}_{([L], 1, 3)}([K_{3,0}]; [L], [L], [L])=3, \Psi^{\P^2}_{(3[L],
1, 9)}([K_{1,8}]; [L], [pt], \cdots, [pt])=27,$$
where $L$ is a line in $\P^2$ and $pt$ is a point.
The first number is supposed to represent the
number of degree 1 elliptic curves with fixed $j$-invariant mapping three
marked points to three distinct lines. But it is well-known that there is no
degree 1 elliptic curves at all in $\P^2$. The second number is supposed to
represent the  number of degree three elliptic curves with fixed
$j$-invariant passing through 8 points. It is known in classical algebraic
geometry that such a number should be 12. What happen was that for the
integrable complex structure $J_0$ on $\P^2$, the boundary of the Gromov-
Uhlenbeck compactification $\overline{\M}_A(1,k, J_0, 0)-\M_A(1,k, J_0, 0)_I$
has a  component whose dimension is larger than the dimension of
$\M_A(1,k,J_0,0)$ itself.
Such a component consists of the maps from the union of
 an elliptic curve and
a rational curve to $\P^2$, which map the elliptic curve to a point.
The effect of considering the inhomogeneous Cauchy Riemann equation
$\bar{\partial}_J f=\nu$ instead of the
homogeneous Cauchy Riemann equation is to
perturb away all those maps. In the process, we creat finitely many
$(J, \nu)$-maps, which provides the correct account of the contribution of
the component described above. Only adding those contributions, our invariants
will satisfy the composition law, while
the classical enumerative invariants
do not.
In fact, the composition law of the mixed invariant we proved in \cite{RT}
computes
all the mixed invariants of any genus for $\P^{N}$ and many other Fano
manifolds.
 It is still a
major problem how to use it to compute the enumerative invariants. Recall that
to define the mixed invariants, we fix the complex structure on the Riemann
surfaces. If we allow the complex structure of Riemann surfaces to vary, it is
not clear how the composition law
will even help to compute the enumerative invariants.

\section{Compactification and Transversality}

In this section, we discuss the structures of  the
moduli space $\M_A(\mu, g,k, J, \nu)$ and the certain quotient of its
Gromov-Uhlenbeck compactifications
$\overline{\M}^r_A(\mu, g,k, J, \nu)$. The problems
we will discuss here are smoothness,  the orientation of
$\M_A(\mu, g,k, J,
\nu)$ and the stratification of $\overline{\M}^r_A(\mu, g,k, J, \nu)$.
By the natural
of those questions, this section is rather technical. For the readers who only
wish to get a sense of big picture, one can skip over this section. On the
other hand, if reader wish to have a good understanding about the results in
this paper, the properties we discuss in this section are crucial or the
proof of all the results in Section 2, which will be provided in the Section 4.
This section is roughly divided as two parts. In the first part, we prove the
smoothness of $\M_A(\mu, g,k, J, \nu)$ and construct its canonical orientation.
The idea of their proof is quite standard.
For the smoothness, the basic
tool is the Sard-Small Transversality Theorem.
We refer the readers to \cite{M1},
\cite{R3}, \cite{M2} in the case of the Cauchy-Riemann equation.
In both \cite{R3} and \cite{M2}, the argument relies on
a cumbersome norm on the space of tamed
almost complex structures defined by
Floer. Because the case (for the
inhomogeneous Cauchy-Riemann equation) here didn't appear
in the literature, we will include an outlined proof. For the
orientation, the genus 0 case was also due to McDuff \cite{M1}. But the
treatment we follow is that of the first author \cite{R}
(see also 4.12 of \cite{RT}),
which was in the spirit of Donaldson's treatment
of the orientation problem
in the gauge theory.  The
second part will be devoted to prove Proposition 2.3,
which is similar to
Section 3, 4, \cite{RT}.

The idea of applying the Sard-Smale Transversality Theorem
to the moduli problem was due to Freed-Uhlenbeck \cite{FU}. Recall that the
Sards-Smale Theorem says that if $X, Y$ are Banach manifolds and
$\F': X\rightarrow Y$ is
a Fredholm map of index $k$, then the set $Y_{reg}$ of regular values of $\F$
is of the Baire second category, provided that $\F$ is sufficiently
differentiable. Recall that
$y\in Y$ is called a regular value, if the derivative
$D\F(x): T_x X\rightarrow T_y Y$ is surjective at any $x$ with
$\F(x)=y$. It
then follows from the Implicit Function Theorem that $\F^{-1}(y)$ is
$k$-dimensional manifold for every $y\in Y_{reg}$.

One obvious problem is that $\M^{\mu}_{g,k}$ may not be smooth.
But it can be stratified by smooth manifolds. Notes that $\M_{g,k}$
has a stratification parameterized by
the automorphism group of Riemann surfaces. Namely, one can write
$$\M_{g,k}=\sum_{\alpha\in I} \T^{\kappa}_{g,k}, \leqno(3.1)$$
where each smooth strata $\T^{\kappa}_{g,k}$ consists of
the stable  Riemann surfaces of a fixed
automorphism group $\kappa$. Without the loss of generality, we can assume that
$$\M^{\mu}_{g,k}=\sum_{\alpha\in I} \T^{\mu, \kappa}_{g,k}, \leqno(3.2)$$
where $\T^{\mu, \kappa}_{g,k}=p^{-1}_{\mu}(\T^{\kappa}_{g,k})$ is smooth.

Using the same arguments as in the smooth case, we will
establish the transversality theory for each stratum
$\T^{\mu,\kappa}_{g,k}$. It is rather straightforward.
The precise structure of $\T^{\mu,\kappa}_{g,k}$ is not needed. One
only has to know that each $\T^{\mu, \kappa}_{g,k}$ is smooth.

Let $\M_A(\mu, g,k, J,\nu)_\kappa$ consist of $(J, \nu)$-maps $f$ such
that the domain of $f$ has automorphism group $\kappa$.
We shall prove that
\vskip 0.1in
\noindent
{\bf Theorem 3.1. }{\it There is a set $\H_{reg}$ of
Baire second category among
all the smooth pairs $(J, \nu)$ such that for any $(J, \nu)$,
$\M_A(\mu, g,k, J, \nu)_\kappa$ is a smooth manifold of dimension $2c_1(V)(A)
-2n ( g-1)+\dim \T^{\mu,\kappa}_{g,k}$.}
\vskip 0.1in
Fix a smooth topological surface $\Sigma_g$ of genus $g$. Our basic
topological object is
$$Map_A(\Sigma_g, V)=\{f:\Sigma_g\rightarrow V \mbox{ such that $f$ is smooth
and } f_*[\Sigma_g]=A\}. $$
To apply the Sard-Smale theorem, we  need to put some Sobolev norm on
$Map_A(\Sigma_g, V)$, so that it has a structure of Banach manifold.

To specify a Sobolev norm, we choose a smooth family of metrics on $\Sigma_g$,
parameterized by the elements of $\T^{\mu, \kappa}_{g,k}$. For example, one can
choose a projective embedding
$$\phi_{\mu}: \overline{\U}^{\mu}_{g,k}\rightarrow \P^N$$
in section 2 and consider the restriction of Fubini-Study metric on
$\phi_{\mu}\pi^{-1}_{\M}(j)$ for each $j\in \T^{\mu, \kappa}_{g,k}$.
 Then for each $j\in \T^{\mu, \kappa}_{g,k}$, $
j$ defines a Sobolev
$L^p_{m, j}$-norm on $Map_A(\Sigma_g, V).$ Its completion under this norm is a
smooth
Banach manifold if $pm>2$. We shall also use $j$ to denote the complex
structure of the underlying Riemann surface. When $j$ varies in $\T^{\mu,
\kappa}_{g,k}$,
$$\chi^{p,m}_{(A,\kappa,g,k)}=\bigcup_{j\in \T^{\mu,\kappa}_{g,k}}L^p_{m, j}
(Map_A(\Sigma_g, V))\times \{j\} \leqno(3.3)$$
is a smooth Banach manifold as well since $\T^{\mu, \kappa}$ is smooth.
Obviously,
there
is a map $\chi^{p,m}_{(A,\kappa,g,k)}\rightarrow \T^{\mu, \kappa}_{g,k}$. Let
$\H^l$ be the completion of $\H$
the space of all smooth pairs $(J, \nu)$ under
$C^{l}$-topology. Then, $\H^l$ is
a smooth Banach manifold. Consider the
{\em universal moduli space}
$$\M^l_A(\kappa,g,k)=\{(f,j, J, \nu)\in \chi^{m,p}_{(A,\kappa,g,k)}\times \H^l;
\bar{\partial}_Jf(x)=\nu(
\phi(x), f(x))\}.\leqno(3.4)$$
When $p>2, 1\leq m\leq l$, by the elliptic regularity, $\M_A^l(\kappa,g,k)$ is
independent of $m, p$.
\vskip 0.1in
\noindent
{\bf Proposition 3.2. }{\it For every $A\in H_2(V, \Z)$ and $g\geq 0, l\geq 1$,
the universal moduli space $\M_A^l(\kappa,g,k)$ is a smooth Banach manifold.}
\vskip 0.1in
\noindent
{\bf Proof:} There is an infinite dimensional vector bundle
$$\E^{m-1,p}_{(f,j,J, \nu)}
\rightarrow \chi^{m,p}_{(A,\kappa,g,k)}\times \H^l, \leqno(3.5)$$
where the fiber $\E^{m-1,p}_{(f,j,J, \nu)}=W^{m-1,p}(\Lambda^{0,1}_{
j}T^*\Sigma_g\otimes_J f^*TV).$ The perturbed holomorphic equation
defines a section of this bundle by
$$\F: \chi^{m,p}_{(A,\kappa,g,k)}\times
\H^l\rightarrow \E^{m-1,p}_{(A,\kappa ,g,k)},\
F(f,j, J, \nu)(x)=\bar{
\partial}_Jf(x)-\nu(\phi(x), f(x)).\leqno(3.6)$$
Notice that the definition of $\bar{\partial}_J$ depends on the complex
structure $j$ on $\Sigma_g$. Then, it is enough to show that
$\F$ is transverse to the zero section. Suppose $\Sigma_j=(\Sigma, x_1, x_2,
\dots, x_k).$ Let $T\Sigma_j=T\Sigma\otimes^k_{i=1} {\cal O}(-x_i)$. Then,
$T_j \T^{\mu, \kappa}_{g,k}= H^{0,1}_{(\kappa,j)}(T\Sigma_j)$-space of $\kappa$
invariant (0,1)-forms. Notice that one can also identity
$$H^{0,1}_{(\kappa,j)}(T\Sigma_j)=(H^{0,1}_{(\kappa,j)}(T\Sigma)\oplus^k_{i=1}
 T_{x_i}\Sigma)^{\kappa},$$
where $(.)^{\kappa}$ means the $\kappa$-invariant subspace.

Let $\F(f, j, J, \nu)=0$. We have
$$\begin{array}{rl}
&T_{(f, j)}\chi^{m,p}_{(A,\kappa,g,k)}=W^{m,p}(\Lambda^0f^*TV)\oplus H^{0,1}_{(
\kappa,j)}(T\Sigma_j);\\
& T_{(J, \nu)}\H^l=C^l(End(TV, J))\oplus C^l(\overline{
Hom}_J(T\P^N, TV)),\\
\end{array}
\leqno(3.7)$$
where $End(TV, J)=\{Y: TV\rightarrow TV; YJ+JY=0\}$.
Furthermore, $\overline{Hom}_J(T\P^N, TV)$ is
the space of anti-complex linear
homomorphism with respect to the complex structure.
There is a natural identification
$$\overline{Hom}_J(T\P^N, TV)|_{\Gamma_f}=\Omega^{0,1}_J(f^*TV),\leqno(3.8)$$
where $\Gamma_f\subset \Sigma_g\times V\subset \P^N\times V$ is the graph of
$f$.

Now, we shall show the surjectivity of the differential
$$\begin{array}{rl}
&D\F(f, j, J, \nu):W^{m,p}(\Lambda^0f^*TV)\oplus H^{0,1}_{(\kappa, j)}
(T\Sigma_j)\oplus C^l(End(TV, J))\oplus C^l(\overline{
Hom}_J(T \P^{N}, TV))\\
&~~~~\longrightarrow W^{m-1,p}(\Lambda^{0,1}_{j}T^*\Sigma\otimes_J f^*TV).\\
\end{array}
\leqno(3.9)$$

An easy computation yields
$$D\F(f, j, J, \nu)(\xi, s,Y,X)=D_{f}\xi+J\circ df\circ
s+ f^*Y\circ df
\circ j-X|_{\Gamma_f}, \leqno(3.10)$$
where
$$D_f: W^{m,p}(\Lambda^0f^*TV)\rightarrow W^{m-1,p}(\Lambda^{0,1}_{j, J}
(f^*TV)) \leqno(3.11) $$
is the linearization of the Cauchy-Riemann equation at $f$.

By the elliptic regularity theory, a $(J, \nu)$-map $f$ is in $W^{l+1, q}$ for
any $q > 0$ if
$(J, \nu)$ is in $C^l$. It follows that $D_f$ in (3.11) is well-defined.
Moreover, it follows that the cokernel of $D_f$ is contained
in $C^l(\Lambda^{0,1}_{j, J}(f^*TV))$. Since $D_f$ is elliptic,
its cokernel is of finite dimension.
However, the map
$$X|_{\Gamma_f}: C^l(\overline{Hom}_J(T\P^N, TV)\rightarrow C^l(
\Lambda^{0,1}_{j, J}(f^*TV))\leqno(3.12)$$
is surjective (here we also use the fact that $f$ is in the space
$C^{l}$). Therefore, by (3.10), $D\F(f,j, J, \nu)$ is surjective.
\vskip 0.1in
By the Implicit Function Theorem, we conclude that the
universal moduli space $\M^l_A(\kappa, g)$ is a smooth Banach manifold.
\vskip 0.1in
\noindent
{\bf Proof of Theorem 3.1:} Based on Proposition 3.2, the proof is just
a standard application of the Sard-Smale Theorem.
For the reader's convenience, we outline the
arguments here.

Consider the projection
$$\pi: \M^l_A(\kappa, g, k)\rightarrow \H^l\leqno(3.13)$$
as a map between the Banach manifolds. The tangent space $T_{(f,j, J, \nu)}
\M^l_A(\kappa, g,k)$ consists of $(\xi, s, Y, X)$ such that
$$D_{f}\xi+J\circ df\circ s+f^*Y\circ
df\circ j-X|_{\Gamma_f}=0. \leqno(3.14)$$
The derivative
$$d\pi: T_{(f,j, J, \nu)}\M^l_A(\mu, g)\rightarrow T_{(J,\nu)}\H^l$$
is just the projection to $(Y, X)$ factors. One can show that
$d\pi$ is a Fredholm operator whose kernel is isomorphic to the kernel of
$D_f\oplus J\circ df$ and has the same index as that of
$D_f\oplus J\circ df$, where
$$J\circ df: H^{0,1}_{(\kappa, j)}( T\Sigma_j)\rightarrow W^{m-1,
p}(\Lambda^{0,1}_{j, J}(f^*TV)). \leqno(3.15)$$
Hence, the operator $d\pi$ is onto precisely when
$D_f\oplus J\circ df$ is onto for any $(J, \nu)$-map $(f, j)$
in $\M^l_A(\mu, g,k,J, \nu)_{\kappa}=
\pi^{-1}((J, \nu)$. In other words,
$$\H^l_{reg}=\{(J, \nu)\in\H^l; D_f\oplus J\circ df
\mbox{ is onto for all } (f, j)
\in \M^l_A(\mu, g,k,J, \nu)_\kappa\}\leqno(3.16)$$
is precisely the set of the regular values of $\pi$. By the
Sard-Smale
Theorem, this set is of the second category. Thus we have proved that
$\H^l_{reg}$ is dense in $\H^l$ with respect $C^l$-topology.
Then one can easily deduce that $\H_{reg}$ is of the second category in
$\H$ with respect to $C^{\infty}$ topology.

Let $\T^{\mu,\kappa}_{g,k}=\bigcup_{K=1}^{\infty} N_K$, where
$N_K$ is compact and
$N_K\subset N_{K+1}$.
Consider the set $\H_{reg, K}\subset \H$ of all smooth $(J, \nu)$ with the
property that the operator $D_f\oplus J\circ df$ is onto for
any $(f, j)$ satisfying: $||df||_
{L^{\infty}}< K$ and $j\in N_K$. Clearly,
$$\H_{reg}=\bigcap_{K>0}\H_{reg, K}.\leqno(3.17)$$
Similarly, we can define $\H_{reg, K}^l$.

We claim that $\H_{reg, K}$ is open and dense in $\H$ with respect to the
$C^{\infty}-topology$.

The openness is clear. It remains
to prove that $\H_{reg, K}$ is dense in $\H$ with respect to
$C^{\infty}$-topology. Note that $\H_{reg, K}=\H^l_{reg, K}\cap \H$.
Then $\H^l_{reg, K}$ is open in $\H^l$ with respect to
$C^l$-topology. Since $\H^l_{reg}\subset \H^l_{
reg, K}$ and $\H^l_{reg}$ is dense in $\H^l$
with respect to
$C^l$-topology, so is $\H^l_{reg, K}$.
It follows that $\H_{reg, K}$ is dense in $\H$ with respect
to $C^{\infty}$-topology. Notice that
$\H_{reg}$ is an intersection of countable
open dense subsets, so it is of
second category.

The dimension formula follows from the Riemann-Roch Theorem.
\vskip 0.1in
\noindent
{\bf Theorem 3.3. }{\it For any $(J, \nu), (J',\nu')\in \H_{reg}$, there is
a second category set of paths $\H_{((J,\nu),(J', \nu'))}$ connecting
$(J, \nu), (J',\nu')$ among the set of all such smooth paths such that for
any path $(J_t, \nu_t)\in \H_{((J,\nu),(J', \nu'))}$
$$\M_A(\mu, g,k,(J_t,\nu_t))_\kappa=\bigcup_{t\in[0,1]}
\M_A(\mu, g,k,J_t,\nu_t)_\kappa\times
\{t\}\leqno(3.18)$$
is a smooth cobordism.}
\vskip 0.1in
The proof is identical to that of 3.1. We omit it.
\vskip 0.1in
\noindent
{\bf Remark 3.4:} In the case of homogeneous Cauchy-Riemann equation ($\nu=0$),
one can use the Teichmuller space $\T_{g,k}$ ( which is smooth) in the place of
$\T^{\mu, \kappa}_{g,k}$. Moreover, there is no need to consider finite covers.
Let $\T^*_A( g,k,J,0)$ be the set of
$(f, \lambda)\in Map_A(\Sigma_g, V)\times \T_{g,k} $
such that $f$ is a $(J,0)$ map for the complex
structure induced by $\lambda$ but not a multiple cover of another $(J,0)$
map. The same argument implies that there is $\H_{reg}$ of second
category  among all the tamed almost complex structure such that for any $J\in
\H_{reg}$, $\T^*_A( g,k,J,0)$ is smooth. The mapping class group $G_g$ acts
freely on $\T^*_A(g,k,J,0)$. Hence,
$$\T^*_A(g,k,J,0)/G_g\subset \M_A( g,k,J,0) \leqno(3.19)$$
is smooth. In our case, the inhomogeneous Cauchy-Riemann equation is not
preserved under the action of mapping class group. Therefore, we have to
consider the smoothness for each strata of $\M_g$. On the other hand, our
inhomogeneous equation can handle the multiple covered map, which can not be
handled by the homogeneous equation except dimension 4 \cite{R3}.
\vskip 0.1in
Next, we construct the canonical orientation of $\M_A(\mu, g,k, J,\nu)_\kappa$
and $\M_A(\mu, g,k, (J_t,\nu_t))_\kappa$.
The construction is identical to that
in \cite{R} (3.3.1) and \cite{RT} (4.12).
\vskip 0.1in
\noindent
{\bf Theorem 3.5. }{\it There is a canonical orientation over $\M_A(\mu, g,k,J,
\nu)_\kappa$ and $\M_A(\mu, g,k,(J_t,\nu_t))_\kappa$.}
\vskip 0.1in
\noindent
{\bf Proof:}
Recall that the linearization of $\F$ at $(f,j)\in \M_A(\mu, g,k,
J,\nu)_\kappa$ is
$$D_f\oplus J\cdot df: \Omega^0(f^*TV)\times H^{0,1}_{\kappa,j}(T\Sigma_j)
\rightarrow \Omega^{0,1}(f^*TV).$$
The tangent space $T_{f,j} \M_A(\mu, g,k, J,\nu)_\kappa=Ker (D_f\oplus J\cdot
df).$ Its determinant is
$$det(T\M_A(\mu, g,k, J,\nu)_\kappa)=det (D_f\oplus J\cdot df),$$
 which is defined over
$Map_A(\Sigma, V)\times \T^{\mu, \kappa}_{g,k}$. As usual, an orientation of
$\M_A(\mu, g,k, J,\nu)_\kappa$ is just a nowhere
vanishing section of $det(T\M_A(\mu, g,k, J,\nu)_\kappa)$ up to
multiplication by positive
functions. We shall omit `` up to multiplication by positive functions'' if
there is no
confusion. Therefore, to construct a
canonical orientation of $\M_A(\mu, g,k, J,\nu)_\kappa$,
it is enough to construct a
canonical section  of $det (D_f\oplus J\cdot df)$ over the
whole $Map_A(\Sigma, V)\times
\T_{\kappa, g}$. We define
$$D_f^J = {1\over 2} ( D_f - J\cdot D_f \cdot J ) \leqno(3.20)$$
Clearly, it is $J$-linear. Moreover, we have
$$D_f=D_f^J + Z_f, \leqno(3.21)$$
where $Z_f$ is
the zero order term. Let
$$D_{f,t}=D^J_f +tZ_f.\leqno(3.22)$$
Then, $det(D_{f, t}\oplus J \cdot df)$ is isomorphic to $det(D_{f, 0}\oplus J
\cdot df)$.
Hence, $det(D_f\oplus J\cdot df)$ is
isomorphic to $det(D^J_f\oplus J\cdot df)$.

On the other hand, both $ker (
D^J_f\oplus J\cdot df)$ and
$coker (D^J_f\oplus J\cdot df)$ are complex vector spaces.
Therefore, there is a canonical
section of the determinant line bundle $det(D_f\oplus id)$
corresponding this
complex structure.

Similarly, one can construct a canonical orientation on
$\M_A(\mu, g,k, (J_t,\nu_t))_{\kappa}$.
\vskip 0.1in
As the
oriented manifolds, we have
$$\M_A(\mu, g,k, J_0, \nu_0)_\kappa \times \{0\}\bigcup \M_A(\mu, g,k,
J_1,\nu_1)_\kappa
\times \{1\}=\M_A(\mu, g,k, J_0, \nu_0)_\kappa\bigcup -\M_A(\mu,
g,k,J_1,\nu_1)_\kappa,$$
where ``-'' means the opposite orientation. Let
$$\M_A(\mu, g,k, J,\nu)=\bigcup_\kappa \M_A(\mu, g,k, J,\nu)_\kappa.$$

In the second half of this section, we focus on
the compactification of ${\M}_A(\mu, g,k,J,\nu)$.
\vskip 0.1in
\noindent
{\bf Definition 3.6 ([PW], [Ye], [Ko]). }{\it Let $(\Sigma, \{x_i\})$ be a
stable Riemann surface.
A stable map (associated with $(\Sigma, \{x_i\})$) is an equivalence class
of continuous maps $f$ from $\Sigma'$ to
$V$ which are smooth at smooth points of $\Sigma'$,
where the domain $\Sigma'$ is obtained by joining chains of $\P^1$'s
at some double points of $\Sigma$ to separate the two components, and then
attaching some
trees of $\P^1$'s. We call components of $\Sigma$ {\em principal components}
and others {\em bubble components}. Furthermore,
\begin{description}
\item[(1)] If we attach a tree of $\P^1$ at a marked point $x_i$, then $x_i$
will
be replaced by a point different from intersection points on some component of
the tree. Otherwise, the marked
points do not change.
\item[(2)] The singularities of $\Sigma'$ are normal crossing and there are at
most two components intersecting at one point.
\item[(3)]  If the restriction of $f$ on a bubble component is constant, then
it has
at least three special points (intersection points or marked points). We call
this component  {\em a ghost bubble} \cite{PW}.
\item[(4)]  For each principal component, the restriction of
$f$ is a $(J,\nu)$-map.
\item[(5)]  The restriction of $f$ to
each bubble is $J$-holomorphic.
\end{description}
Two such maps are equivalent if one is the composition of the other with
an automorphism of bubble components fixing
the special points.}
\vskip 0.1in
Evidently, the equivalence relation is trivial unless some bubble component
has one or two special points. The restriction of $f$ to each component
carries a homology class. We shall use $[f]$ to denote the summation of
all those homology classes.
\vskip 0.1in
\noindent
{\bf Remark 3.7: } The terminology {\em stable maps} was first used by
Kontsevich and Manin in
\cite{KM}.
It had appeared before
in Parker-Wolfson-Ye's proof of Gromov-Uhlenbeck compactness
theorem under the name {\it cusp curves} \cite{PW}, \cite{Ye}. Later, it was
introduced to algebraic geometry by Kontsevich and Manin and known more
commonly as stable maps. Here, we follow their terminology.
\vskip 0.1in
Then Theorem 3.1 of \cite{RT} (see also \cite{PW} and \cite{Ye}) can be
restated as follows:
\vskip 0.1in
\noindent
{\bf Theorem 3.8. }{\it Let  $f_m\in \M_A(\mu, g,k, J, \nu)$. Suppose that the
domain $(\Sigma_m, \{x^m_i\})$ of $f$ converges to a stable Riemann surface
$(\Sigma, \{x_i\})$ in the sense of Deligne-Mumford. Then, there is a
subsequence $\{f_{m_t}\}$  ``weakly converging'' to a stable map $f$(
associated with $(\Sigma, \{x_i\})$) such that $[f]=A$.
Here, by the weak convergence, we mean that the
image of $f_{m_t}$
converges  to the image of $f$ in the Hausdorff topology.}
\vskip 0.1in
Strictly speaking, Proposition 3.1 of \cite{RT} only proves the version of
Theorem 3.7 without marked points. But one can easily keep track of marked
points in the proof and deduce Theorem 3.7 as we stated.

We denote the space of  stable maps with fundamental class $A$
by $\overline{\M}_A
(\mu, g,k,J,\nu)$. Clearly,
$$\overline{\M}_A(\mu, g,k,J,\nu)\supset \M_A(\mu, g,k,J,\nu). \leqno(3.23)$$
One can easily deduce from Theorem 3.7 the following:
\vskip 0.1in
\noindent
{\bf Corollary 3.9. }{\it $\overline{\M}_A(\mu, g,k,J,\nu)$ is compact
in the Hausdorff topology. Moreover,
the evaluation map $e_i$ extends to a continuous map from
$\overline {\M} _A(\mu, g,k,J,\nu)$.}
\vskip 0.1in
We shall call $\overline {\M}_A(\mu, g,k,J,\nu)$ GU-compactification of
$\M_A(\mu, g,k,J,\nu)$, since Gromov and Uhlenbeck first studied
the compactness problem for harmonic maps and pseudo-holomorphic
curves.
\vskip 0.1in
\noindent
{\bf Definition 3.10. }{\it We call $f$ a reduced GU-map if $f$ satisfies (1),
(4),(5) of Definition 3.6. Furthermore, it satisfies:
\begin{description}
\item[(2')] The singularities of $\Sigma'$ are of normal crossing, but it
could have three or more components intersecting at one point;
\item[(3')] There are no ghost bubbles;
\item[(6)] There are no bubble components which are  multiple covering maps;
\item[(7)] There are no subtrees of the bubbles whose components have
the same image.
\end{description}}
\vskip 0.1in
For any stable map, we can define a  GU-map (with possibly different
fundamental class) as follows: (i) we
collapse the ghost bubbles; (ii) we replace each multiple covered bubble
component by its
image; (iii) we collapse each subtree of the bubbles whose components have the
same image.

Clearly, this reduction may destroy the property
(2) of the Definition 3.6, but still preserve the property (2').
Also in this reduction, the fundamental class may change.

Define
$\overline{\M}^r_A(\mu, g,k,J,\nu)$ to be the quotient of $\overline{\M}_A(\mu,
g,k,J,
\nu)$ by the above reduction. Furthermore, we define the topology
on $\overline{\M}^r_A(\mu, g,k,J,\nu)$ as the quotient topology.

By the definition, $\overline{\M}^r_A(\mu, g,k,J,\nu)$ is a union of
GU-maps
with possibly different fundamental classes. We will prove the following
structure theorem.
\vskip 0.1in
\noindent
{\bf Theorem 3.11. }{\it Let $(V, \omega)$ be a semi-positive symplectic
manifold. There is a set $\H_{reg}$ of Baire second category among all the
smooth pairs $(J, \nu)$ such that for any $(J,\nu)\in \H_{reg}$,
$\overline{\M}^r_A(\mu, g,k, J, \nu)-\M_A(\mu, g,k,
 J, \nu)_I$ consists of finitely many strata, such that
each stratum  is a smooth manifold of
real codimension at least 2.}
\vskip 0.1in
\noindent
{\bf Proof of Proposition 2.3.} (1) follows from Theorem
3.1. (2) is obvious (cf. Corollary 3.9).
By the construction of $\overline{
\M}^r_A(\mu, g,k,J,\nu)$, both $\Upsilon$ and $\Xi^A_{g,k} $ descend to
$\overline{\M}^r_A(\mu, g,k,J,\nu)$. Then (3) follows from Theorem 3.11.
\vskip 0.1in
\noindent
{\bf Remark 3.12:} One may ask whether or not
$\overline{\M}^r_A(\mu, g,k, J,\nu)$ carries a fundamental class.
This is indeed the case if $\overline{\M}^r_A(\mu, g,k, J,\nu)$ admits a real
analytic structure. These will be established in
\cite{RT1} by
more delicate analysis.
Then one can directly use the GU-compactification
to prove Proposition 2.3.
\vskip 0.1in
In the rest of this section, we outline the proof of Theorem 3.11. The
proof is identical to Section 4 of \cite{RT}. We refer the readers to
\cite{RT} for certain details.

First we shall decompose $\overline{\M}^r_A(\mu, g,k, J, \nu)-\M_A(\mu, g,k,J,
\nu)$
into strata. A stratum is the set of GU-maps (possibly with
total homology class different from $A$) satisfying:
(1) their domains with
marked points  are of
the same homeomorphic
type; (2) Each connected component carries a fixed homology class.
Furthermore, for technical
reasons, we need to specify those bubble components,  which have the
same image even though
they may not be adjacent to each other, and their intersection points having
the
same image.  Therefore,
the strata of $\overline{\M}^r_A(\mu, g,k, J, \nu)$ are indexed by  data:
(i) homeomorphism type of the domain of GU-maps with marked points; (ii) a
homology class associated to each
component; (iii) a specification of components with the same image and their
intersection points with the same image. We
denote by $D$ a
set of those three data. Let $\D_{g,k}$ be the collection of
such $D$'s. Note that when we drop
the multiplicity from a multiple covering
map, we change the homology class.
However it is still $A$-admissible in the following sense:
\vskip 0.1in
\noindent
{\bf Definition 3.13.} {\it Let $D$ be given as above. We
define $[D]$ to be the sum of homology classes of components in (ii).
Let $P_1, \cdots, P_o$ be principal components and $B_1,
\cdots, B_p$ be bubble components of $D$. We say that
$D$ is called $A$-admissible if there
are positive integers $b_1, \cdots, b_k$ such that
$$A=\sum^o_1 [P_i]+\sum^p_1 b_j [B_j]\leqno(3.24)$$
where $[P_i]$, $[B_j]$ are the homology classes of $P_i$, $B_j$.
We say that $D$ is $(J, \nu)$-effective if every principal
component can be represented
by a $(J, \nu)$-map and every bubble component can be represented by a
$J$-holomorphic map. }
\vskip 0.1in
We will always denote by $\Sigma _i$ the domain of the $(J,\nu)$-map
representing $P_i$. Let $\D^{J, \nu}_{g,k}\subset \D_{g,k}$
be the set of $A$-admissible,
$(J, \nu)$-effective $D$.
\vskip 0.1in
\noindent
{\bf Lemma 3.14.} {\it The set $\D^{J, \nu}_{g,k}$ is finite.}
\vskip 0.1in
\noindent
This is the analogue of Lemmas 4.5 of \cite{RT} and a simple
corollary of the Gromov-Uhlenbeck compactness theorem (cf. Theorem 3.8).
The presence of marked points
doesn't affect the proof at all. We omit it.

One can consider $\D^{J, \nu}_{g,k}$ as the set of indices of strata.
For each $D\in \D^{J, \nu}_{g,k}$, let $\M_{D}(\mu, g,k, J, \nu)$ be the space
of
GU-maps such that the homeomorphism type of its domain with marked
points,
homology class of each component, and components and their intersection points
which have the same image  are specified by $D$.
The following lemma can be deduced from
the definition.
\vskip 0.1in
\noindent
{\bf Lemma 3.15. }{\it
$$\overline{\M}^r_A(\mu, g,k,J,\nu)=\bigcup_{D\in \D^{J,\nu}_{g,k}} \M_D(\mu,
g,k,J,\nu).\leqno(3.25)$$
}
\vskip 0.1in
By the definition, each $D$ is associated with a stable Riemann surface
$\Sigma_{D, r}$, which can be obtained by contracting all the bubble
components. Recall that for each principal component we have to fix the
automorphism
group preserving the special points to make the transversality arguments work
(Theorem 3.1). Here we use
$\bar {\kappa}$ to denote the $o$-tuple $(\kappa _1, \cdots, \kappa_o)$,
where $\kappa_i$ is the
automorphism group of the principal component $P_i$ fixing the special points
(marked points and intersection points) of $\Sigma_{D, r}$.
Then,
$$\M_D(\mu, g,k,J,\nu)=\bigcup_{\bar {\kappa}}\M_D(\mu,
g,k,J,\nu)_{\bar{\kappa}},\leqno(3.26)$$
where $\M_D(\mu, g,k,J,\nu)_{\bar{\kappa}}$ consists of all maps in
$\M_D(\mu, g,k,J,\nu)$ whose $i^{\rm th}$-principal component has
the automorphism group $\kappa_i$. Next we prove the smoothness of
$\M_D(\mu, g,k,J,\nu)_{\bar{\kappa}}$.
First we make another reduction by identifying the domains of those bubble
components which have the same image,
and change the homology class accordingly. Furthermore, we identify the
corresponding intersection points with the same image. Suppose that the
resulting new domain and homology class of each component are specified
by $\bar{D}$. This process may destroy the tree structure and
creat some cycles in the domain.
The total homology class
may also change. However, it remains to be $A$-admissible.

Given such $D$ and $\bar{D}$, we can identify
$\M_D(\mu, g,k,J,\nu)_{\bar{\kappa}}$ with
the space of $(J, \nu)$-maps whose
domain, homology class of each component are specified by  $\bar{D}$ and the
automorphism group of its principal components are specified by $\bar
{\kappa}$.
Let us denote this space by $\M_{\bar D}(\mu, g,k,J,\nu)_{\bar{\kappa}}$. Then,
$$\M_D(\mu, g,k,J,\nu)_{\bar{\kappa}}=\\M_{\bar D}(\mu,
g,k,J,\nu)_{\bar{\kappa}}.\leqno(3.27)$$
For each $f$ in $\M_{\bar D}(\mu, g,k,J,\nu)_{\bar{\kappa}}$,
the bubble components have different images.  As before,
let $P_1, \cdots, P_o$
be the
principal components and let $B_1, \cdots, B_p$ be
the bubble components.
Let $\Sigma _{\bar D}$ be the domain of maps in the stratum
$\M_{\bar{D}}(\mu, g,k, J, \nu)_{\bar {\kappa}}$. This is a union
of $\Sigma _i$ (genus $g_i$) and some $S^2$'s intersecting each other
according to the intersection pattern given by $D$. Let $h_i$ be the number of
intersection points on the component
$P_i$. Note that we count a self-intersection point twice. Here, the
intersection points between the components are the points in
their domain, not in their image.
Similarly,
let $h^j$ be the number of intersection points on the bubble component
$B_j$. Let $k_i(k^j)$ be the number of marked points on $P_i$-component (bubble
component $B_j$), which are different from intersection points. Notice that
$k^j=0$ or $1$. Moreover,
$$\sum k_i +\sum k^j\leq 2k.\leqno(3.28)$$
It may not be equal to $2k$ because the collapsing of the ghost bubbles
containing a marked point forces the marked point to lie on the intersection.
Let $\M^*_{[B_j]}(S^2, J,0)\subset
\M_{[B_j]}(S^2, J,0)$ be the space of non-multiple covering maps and
$$\M^*_{[B_j]}(S^2, h^j+k^j, J, 0)=\M^*_{[B_j]}(S^2, J,0)\times \overline{
S^2}^{h^j+k^j}$$
 where $\overline{S^2}^{h^j+k^j}$ is the set of distinct $h^j+k^j$-tuple points
of $S^2$.
Consider
$$
\tilde{\M}_{\bar{D}}(\mu, g,k, J, \nu)_{ \bar {\kappa}}
= \{f: \Sigma _{\bar D}
\rightarrow V ~|~
f_{P_i}\in \M_{[P_i]}(\mu, g_i, h_i+k_i, J, \nu_i)_{\kappa_i},
Im(f_{B_j})\neq Im(f_{B_{j'}}) ~~\mbox{ if }
j\neq j' \}
\leqno(3.29)
$$
For each bubble component, there is a parameterization group $G=PSL_2$.
Therefore, $G^{p_{\bar{D}}}$ acts on $\tilde{\M}_{\bar{D}}(\mu, g, J, \nu)_{
\bar {\kappa}}$, where $p_{\bar{D}}$ is the number of bubble components.
Then $\M_{\bar D}(\mu, g,J,\nu)_{\bar{\kappa}}=
\tilde{\M}_{\bar{D}}(\mu, g,k, J,\nu)_{ \bar {\kappa}}/G^{p_{\bar{D}}}$.
Clearly,
$$\tilde{\M}_{\bar D}(\mu, g,k,J,\nu)_{\bar{\kappa}} \subset \prod \M_{[P_i]}
(\mu, g_i, h_i+k_i, J, \nu_i)_{\kappa_i}
\times  \prod \M^*_{[B_j]}(S^2, h^j+k^j, J,0),\leqno(3.30)$$
whose components intersecting each other according to the intersection pattern
given by $\bar{D}$.
Consider the evaluation map
$$e_{\bar D}:
\prod \M_{[P_i]}(\mu, g,h_i+k_i, J,\nu)_{{\kappa _i}}\times \prod
\M^*_{[B_j]}(S^2, h^j+k^j, J,0)\mapsto \prod
V^{h_i}\times \prod V^{h^j}=V^{h_{\bar{D}}},\leqno(3.31)$$
where $h_{\bar{D}}=\sum h_i +\sum h^j$, and $e_{\bar{D}}$ is defined as
follows:
We first define
$$\begin{array}{rl}
&e_{P_i}:\M_{[P_i]}(\mu, g, h_i+k_i, J,\nu)_{{\kappa _i}}\rightarrow V^{h_i}\\
&e_{P_i}(f, x_1, \cdots, x_{h_i}, x_{h_i+1}\cdots x_{h_i+k_i})=(f(x_1), \cdots,
f(x_{h_i}))
\end{array}
\leqno(3.32)$$
For each $B_j$, we define $e_{B_j}: \M^*_{[B_j]}(S^2, h^j+k^j, J,0)
\rightarrow
 V^{h^j}$ by
$$e_{B_j}(f, y_1, \cdots, y_{h^j}, y_{h^j+1}\cdot y_{h^j+k^j})=(f(y_1), \cdots,
f(y_{h^j})).
\leqno(3.33)$$
Then we define $e_{\bar{D}}=\prod e_{P_i}\times \prod e_{B_j}$.
Recall that if $M, N$ are
submanifolds of $X$, $M\cap N$ can be interpreted as $M\times N\cap \Delta$,
where $\Delta$ is the diagonal of $X\times X$. This means that we can realize
any intersection pattern by constructing a ``diagonal'' in the product.
Let us construct a submanifold $\Delta_{\bar{D}}\subset V^h$ which
plays the role of the diagonal. Let $z_1,\cdots,z_{t_{\bar{D}}}$ be all
the intersection
points. For each $z_s$, let
$$I_s=\{P_{i_1}, \cdots, P_{i_q}, B_{j_1}, \cdots,
B_{j_r}\}$$
be the set of components which intersect at $z_s$.
Now we will construct a product $V_s$ of $V$ such that
its diagonal describes the intersection at
$z_s$. This is done as follows: we allocate one or two factors from
each of $V^{h_{i_1}}, \cdots,
V^{h_{i_q}}$, according to whether or not $z_s$
is a self-intersection point of the corresponding principal component.
We allocate one
factor from each of
$V^{h^{j_1}}, \cdots, V^{h^{j_r}}$.
Here $V^{h_i}$ or $V^{h^j}$
are the image of $e_{P_i}$ or $e_{B_j}$.
Then, we take the product of
those factors and denote it by $V_s$.
Let $\Delta_s$ be the diagonal of $V_s$. Then the product
$\Delta_{\bar{D}}=\Delta_1\times \cdots \times
\Delta_{t_{\bar{D}}}\subset V^{h_{\bar{D}}}$ is the diagonal to realize the
intersection pattern between the components of $\bar{D}$.
Then
$e^{-1}_{\bar{D}}(\Delta_{\bar{D}})\supset \tilde{\M}_{\bar D}
(\mu, g,k, J, \nu)$. But
they may not be equal because we require that bubble
components have different image. But $\tilde{\M}_{\bar D}
(\mu, g,k, J, \nu)$ is an open subset.
Moreover, the group $G^{p_{\bar{D}}}$ acts
on $e^{-1}_{\bar D}(\Delta_{\bar{D}})$.
\vskip 0.1in
\noindent
{\bf Theorem 3.16.} {\it There is a set $\H_{reg}$ of Baire second category
among all the smooth pairs $(J,\nu)$ such that for any  $(J, \nu)\in
\H_{reg}$, $\M_{\bar D}(\mu, g,k,J,\nu)_{\bar{\kappa}}$ is a
smooth manifold of dimension
$$\sum(2c_1(V)(P_i)+2n(1-g_i))+\sum \dim \T^{\mu, \kappa}_{g_i,0}+\sum
(2c_1(V)(B_j)+2n-6)+2h_{\bar{D}}+\sum k_i+\sum
k^j-2n(h_{\bar{D}}-t_{\bar{D}}),$$
where $g_i$ is the genus of $\Sigma_i$ and $t_{\bar{D}}$ is the number of
intersection points of $
\bar{D}$.

Moreover, for any $(J,\nu)$ and $(J', \nu')$ of $\H_{reg}$, there is
a second category set of paths $\H_{((J,\nu), (J', \nu'))}$  connecting
$(J,\nu)$ and $(J', \nu')$ among all the smooth paths such that for any path
$(J_t, \nu_t)\in \H_{((J,\nu), (J', \nu'))}$,
$$\bigcup_{t\in [0,1]}\M_{\bar  D}(\mu, g,k,J,\nu)_{\bar{\kappa}}$$
 is a smooth cobordism of one dimension higher.}
\vskip 0.1in
By the construction of $\bar{D}$, it is evident that $t_{\bar{D}}\leq t_D$ and
$h_{\bar{D}}\leq h_D$. But, $h_{\bar{D}}-t_{\bar{D}}=h_D-t_D$. Therefore,
\vskip 0.1in
\noindent
{\bf Corollary 3.17.} {\it Under the conditions of Theorem 3.16,
the dimension of $ \M_{\bar D}(\mu, g,k,J,\nu)_{\bar{\kappa}}$ is less than or
equal to
$$\sum(2c_1(V)(P_i)+2n(1-g_i))+\sum \dim \T^{\mu, \kappa}_{g_i,0}+2k+\sum
(2c_1(V)(B_j)+
2n-6)+2h_{D}-2n(h_{D}-t_{D}).$$
}
\vskip 0.1in
\noindent
{\bf Proof of Theorem 3.16:} The idea of the proof
is similar to that in the proof of
Lemma 4.8-4.11 in \cite{M1}.
Also many arguments are the same as those in the proof of Theorem
4.7 in \cite{RT}. But we will avoid the Floer's
norm on the space of almost complex structures as we did before.

First of all, $\tilde{\M}_{
\bar{D}}(\mu, g,k,J,\nu)_{\bar{\kappa }}$ is an open subset of
$e^{-1}_{\bar{D}}(\Delta_{
\bar{D}})$. Hence for the purpose of proving smoothness, we can assume that
$$\tilde{\M}_{\bar{D}}(\mu, g,k,J,\nu)_{\bar{\kappa }}=
e^{-1}_{\bar{D}}(\Delta_{
\bar{D}})$$
to simplify the argument.

Suppose $p>2,m\geq 1$. Following (3.3),
$$\chi^{p,m}_{\bar {\kappa}, \bar{D}}=\prod_{i=1}\bigcup_{j\in
\T^{\mu,\kappa}_{g_i,
h_i+k_i}} L^p_{m,j}
(Map_{[P_i]}(\Sigma_i, V))\times \{j\} \times \prod_{s=1}L^p_k (Map_{[B_s]}
(\P^1, V))\times (\overline{S^2})^{h^j+k^j}\leqno(3.34)$$
is a smooth Banach manifold. Then, we define
$$\M^l(\bar{D}, \bar {\kappa})=\{(f,j,J,\nu)\in \chi^{p,m}_{\bar
{\kappa},\bar{D}}\times
\H^l ; \bar{\partial}_J f(x)=\nu(\phi(x), f(x))\},\leqno(3.35)$$
where $m\leq l$ and the equation
$$\bar{\partial}_J f(x)=\nu(\phi(x), f(x))$$
means the $(J,\nu)$-holomorphic equation for each $P_i$ component and
$J$-holomorphic equation for each $[B_s]$ component.
It is not hard to observe that
$$\M^l(\bar{D}, \bar {\kappa})=\bigcup_{(J,\nu)\in \H^l}(\prod_1 \M^l_{[P_i]}(
\mu, g_i,h_i+k_i, J,\nu)_{{\kappa _i}} \times \prod_1 (\M^*)^l_{[B_s]}(S^2,h^j
+k^j, J,0).
\leqno(3.36)$$
\vskip 0.1in
\noindent
{\bf Lemma 3.18. }{\it $\M^l(\bar{D}, \bar {\kappa})$
is smooth Banach manifold.}
\vskip 0.1in
\noindent
{\bf Proof of Lemma 3.18: } $\M^l(\bar{D}, \bar {\kappa})$
is just the analogue of
$\Theta(\Delta_{\J})$ in Lemma 4.9 of \cite{RT}.
The proof of Lemma 3.18
is identical to that of Lemma 4.9 in \cite{RT}. We omit it.
\vskip 0.2in
Note that the  evaluation map extends
$$e_{\bar{D}}: \M^l(\bar{D}, \bar {\kappa})\rightarrow
V^{h_{\bar{D}}}.\leqno(3.37)$$
We define
$$\tilde{\M}^l_{\bar{D}}(\mu, g, \bar
{\kappa})=e^{-1}_{\bar{D}}(\Delta_{\bar{D}})\leqno(
3.38)$$
\vskip 0.1in
\noindent
{\bf Lemma 3.19. }{\it $\tilde{\M}^l_{\bar{D}}(\mu, g, \bar {\kappa})$
is a smooth Banach manifold.}
\vskip 0.1in
Its proof is identical to the proof of Theorem 4.7 of \cite{RT}. We omit it.

Then, the rest of proof of Theorem 3.16 is similar to that of Theorem 3.1.
We sketch the argument.

Consider the projection
$$\pi: \tilde{\M}^l_{\bar{D}}(\mu, g, \bar {\kappa})\rightarrow \H^l.\leqno(3.
39)$$
Then,
$$\tilde{\M}^l_{\bar{D}}(\mu, g,k,J,\nu)_{\bar
{\kappa}}=(\pi)^{-1}((J,\nu)).\leqno(3.40)$$
Then $\pi$ is a Fredholm map between two Banach manifolds. It follows from
the Sard-Smale
Transversality Theorem that the set
$$\H^l_{reg}=\{(J, \nu)\in\H^l\,|\, d\pi \mbox{ is onto for all } f\in
\M^l_{\bar{D}}
(\mu, g,k,J,\nu)_{\bar {\kappa}}\}$$
is of Baire second category. Let
$$\H_{reg}=\bigcap_l \H^l_{reg}.$$
Then,
$$\tilde{\M}_{\bar{D}}(\mu, g,k,J,\nu)_{\bar {\kappa }}=\bigcap_l
\M^l_{\bar{D}}(\mu, g,k,J,\nu)_{\bar {\kappa }}.
\leqno(3.41)$$
As in the proof of Theorem 3.1, we can deduce that $\H_{reg}$ is
of Baire second category. We leave it to the readers as an exercise.
Then, for any
$(J,\nu)\in \H_{reg}$,
$$\tilde{\M}_{\bar{D}}(\mu, g,k,J,\nu)_{\bar {\kappa }}=\pi^{-1}((J, \nu))$$
is a smooth manifold. Since $G^{p_{\bar{D}}}$ acts freely on
$\tilde{\M}_{\bar{D}}(\mu, g,k,J,\nu)_{\bar {\kappa }}$,
$$\M_{\bar D}(\mu, g,k,J,\nu)_{\bar {\kappa }}=\tilde{\M}_{\bar{D}}(\mu, g,k,
J,\nu)_{\bar{\kappa}}/G^{p_{\bar{D}}}$$
is smooth.

An routine counting argument yields the dimension formula.

The proof of the second part of Theorem 3.16 is identical.
\vskip 0.1in
Recall that  if we contract all the bubble components of $D$, we obtain the
stable Riemann surface $\Sigma _{D, r}$ in the sense of Deligne-Mumford.

An interesting special case of Theorem 3.18 is
\vskip 0.1in
\noindent
{\bf Corollary 3.20.} {\it If $\Sigma _D$ has no bubble components at all,
i.e., $\Sigma _{D}=\Sigma_{D,r}$,
$\M_{\Sigma_D}(\mu, g,k,  J, \nu)_{\bar {\kappa}}$
is smooth for a generic
$(J, \nu)$. Moreover, for a generic $(J_t, \nu_t)$,
$\bigcup_t \M_{D}{(\mu, g,k,J_t,\nu _t)_{\bar {\kappa }}}
\times \{t\}=\bigcup_t \M_{\Sigma_D}(\mu, g,k,J_t,\nu _t)_{\bar {\kappa
}}\times \{t\}$ is a smooth cobordism.
Here, the word ``generic'' means that it is in a set of Baire second
category.}
\vskip 0.1in
Next, we compute the codimension of $\M_{\bar D}(\mu, g,k,J,\nu)_{{\kappa
_i}}$. First
\vskip 0.1in
\noindent
{\bf Proposition 3.21.} {\it Suppose that $(V, \omega)$ is a semi-positive
symplectic
manifold. Let $\M_{\Sigma_{D,r}}\subset \overline{\M}_{g,k}$
be the set of stable Riemann surfaces such that their homeomorphism types are
specified by $\Sigma_{D, r}$. Then,
$$\dim \M_{\bar{D}}(\mu, g,k,J,\nu )_{\bar {\kappa }}\leq
2c_1(V)(A)+2n(1-g)+\dim \M_{
\Sigma_{D,r}}-2p_D,\leqno(3.42)$$
where $p_D$ is the number of bubble components of $D$(not $\bar{D}$)}
\vskip 0.1in
\noindent
{\bf Proof:} By Corollary 3.17,
the dimension of $\M_{\bar{D}}(\mu, g,k,J,\nu)_{\bar {\kappa }}$ is less than
or equal to
$$\sum_i (2c_1(V)([P_i])+2n(1- g_i))+2k+\sum_i \dim \T^{\mu, \kappa_i}_{g_i,0}+
\sum_j (
2c_1(V)(B_j)+(2n-6))+2h_{D}-2n(h_{D}-t_{D})$$
$$=2c_1(V)([\bar{D}])+2n\sum_i(1-g_i)+2k+\sum \dim_i \T^{\mu, \kappa_i}_{g_i,0}
+(2n-6)p_{\bar{D}}+2h_D-2n(h_{D}-t_{D})$$
For a generic $J$,
$$2c_1(V)(B_j)+2n-6=\dim \M^*_{[B_j]}(S^2, J,0)/PSL_2\geq 0. \leqno(3.43)$$
If some bubble component $B_j$ happens
to be the image of two or more bubble
components of $D$, by adding $2c_1(V)(B_j)+2n-6$ to the dimension formula,
$$\dim \M_{\bar{D}}(\mu, g, k,J, \nu)\leq  2c_1(V)([D])+2n\sum_i(1- g_i)+2k+
\sum_i \dim \T^{\mu, \kappa}_{g_i,0}+(2n-6)p_{D}+2h_{D}-2n(h_{D}-t_{D}).$$
Since $(V, \omega)$ is semi-positive, $c_1(V)(B_j)\geq 0$ for a generic
$J$. Since $D$ is $A$-admissible,   $c_1(V)([D])\leq c_1(V)(A)$.  Let
$$\lambda_D=(2n-6)p_D+2h_D-2n(h_D-t_D).\leqno(3.44)$$
Then, Proposition 4.13 of $\cite{RT}$ implies that
$$\lambda_D\leq \lambda_{\Sigma_D}-2p_D. \leqno(3.45)$$
Note that in \cite{RT}, we use $k_D$ to
denote the number of bubble component
instead of $p_D$ we used in this paper
(here $k$ was used to denote the number
of marked points). Therefore,
$$dim \M_{\bar{D}}(\mu, g,k,J,\nu )_{\bar {\kappa }}
\leq \sum_i 2c_1(V)(A)+2n\sum_i(1- g_i)+2k+
\sum_i \dim \T^{\mu, \kappa_i}_{g_i,0}+
2h_{\Sigma_D}-2n(h_{\Sigma_{D}}-t_{\Sigma_{D}})-2p_D.
$$
Since $\Sigma_{D,r}$ is homeomorphic to a stable curve, $t_{\Sigma_D}$
is the number of double points and $h_{\Sigma_D}=2t_{\Sigma_D}$. An easy
inductive argument (Proposition of 4.13, \cite{RT}) shows that
$$o-t_{\Sigma_D}-\sum_i g_i=1-g.\leqno(3.46)$$
It is easy to observe that
$$dim \M_{\Sigma_{D,r}}=2k+\sum_i \dim
\T^{\mu,\kappa_i}_{g_i,0}+2h_{\Sigma_D}\leq \dim \M_{\Sigma_D} \leqno(3.47)$$
and equal iff all the $\kappa_i$ are trivial.
\vskip 0.1in
\noindent
{\bf Proof of Theorem 3.11 (Structure Theorem):} It follows from Lemma 3.14,
3.15, Theorem 3.16, Proposition 3.21.

\section{Proof of Proposition 2.4, 2.9 and  2.10}

In this section, we first establish the existence
of the GW-invariants $\Psi^V$ (Proposition 2.4) and
its independence from various parameters.
Hence, $\Psi^V$ is a symplectic invariant. Then, we will prove Proposition
2.9 and 2.10.
Technically speaking, this section is the analogue of
section 5 and 7 of \cite{RT}.
We shall repeatly use the word ``generic'' to
mean something belonging to a set of Baire second category.

First of all, we extend the definition of
pseudo-submanifolds to the singular
space with  quotient singularities such has $\overline{\M}^{\mu}_{g,k}$.
Furthermore, we also need to consider the transversality theory of such
pseudo-submanifolds. It is well-known that over the rational coefficient $\Q$,
the usual theory for the smooth manifolds extends to the singular space with
quotient singularities, where the Poincare duality holds over $\Q$.
\vskip .1in
\noindent
{\bf Definition 4.1: }{\it An n-dimensional finite simplicial complex $P$ is
called an abstract pseudo-manifold if $P^{top}=P-P_{n-2}$ ($n-2$ skeleton)
is an open smooth oriented $n$-dimensional manifold.
$P$ is called an abstract
pseudo-manifold with boundary if $P^{top}$ is a
$n$-dimensional oriented smooth
manifold with boundary $\partial P^{top}$.
Let $\partial P=\overline{\partial
P^{top}}$. Then $\partial P \cap P_{n-2}$ is a subcomplex of dimension less
than or equal to $n-3$. Let $V$ be a stratified PL-manifold such that each
stratum is even dimensional and its triangulation is compatible with the
stratification. A pseudo-submanifold is a pair $(P, f)$, where $P$ is
an abstract pseudo-manifold and $f: P\rightarrow V$ is a piece-wise linear
map (PL) with respect to  some triangulation of $V$. Furthermore, we require
that $f$ maps $P^{top}$ into one stratum and smooth.
A pseudo-submanifold
cobordism between pseudo-submanifolds
$(P, f), (Q, h)$ is a pair $(K, H)$ such
that $K$ is an abstract pseudo-manifold
with boundary with $\partial K= P\cup
-Q$ and $H$ is PL with respect to some
triangulation of $V$ and  smooth over
$K^{top}$ in the sense that $H$ maps the $K^{top}$
smoothly in one stratum or maps
the interior $K^o$ of $K^{top}$ to
one stratum and $\partial K^{top}$ to the
lower strata.
Moreover,  $H|_{P\cup -Q}=f\cup - h$,
where $-$ means the opposite orientation.}
\vskip .1in
Furthermore, we have the following lemma on transversality.
\vskip 0.1in
\noindent
{\bf Lemma 4.2. } {\it Let $V$ be a stratified PL-manifold with fundamental
class $[V]$ such that the Poincare duality holds over the rational coefficient.
 Then, for each homology class $\alpha$, there exists $p$ and
 a pseudo-submanifold representative $(P,f)$
of a homology class $p(\alpha)$. Furthermore,
if $h_i: X_i \rightarrow V$ be
smooth maps from
smooth manifolds $X_i$ to smooth strata of $V$, then there is a small
perturbation $\tilde{f}: P
\rightarrow V$ such that $\tilde{f}$ is transverse to each $h_i$, i.e.,
$\tilde{f}$ is
transverse to $h_i$ as a PL map and transverse over $P^{top}$ as a smooth map.}
\vskip 0.1in
\noindent
{\bf Proof:} This lemma is the consequence of
standard transversality results
in PL topology  \cite{Mc}(Theorem 5.2). One remark is that if $V$ is
a smooth manifold, this lemma holds  for any $\alpha$. Otherwise, the lemma
holds for those class of the form $[V]\cap \beta^*$ for a cohomology class
$\beta^*$ (Theorem 5.2 of \cite{Mc}). Then, the lemma follows from the
assumption
that the Poincare Duality holds over the rational coefficient $\Q$.
\vskip 0.1in
Recall that in (2.4), we have defined the evaluation map
$$e_i: \M_A(\mu, g,k, J,\nu)_\kappa\rightarrow V.\leqno(4.1)$$
$e_i$ extends obviously to each
$\M_{\bar D}(\mu, g,k,J,\nu )_{\bar {\kappa }}$ (Lemma 3.15).
We shall still denote it by $e_i$. Furthermore,
$$\Upsilon: \M_A(\mu, g,k, J,\nu)_\kappa\rightarrow
\M^{\mu}_{g,k,\kappa},\leqno(4.2)$$
which takes the domain of maps in $\M_A(\mu, g,k, J,\nu)_\kappa$ (section 2)
extends over
$$\Upsilon: \M_{\bar D}(\mu, g,k,J,\nu )_{\bar {\kappa }}\rightarrow \M^{\mu}_{
\Sigma_{D,r}, \bar{\kappa}}$$
as well. There is an obvious
version of the maps $e_i$ and $\Upsilon$ for the corbordisms $\bigcup_t
\M_A(\mu, g,k,J_t,\nu_t)_\kappa \times
\{t\}$ and
$\bigcup_t \M_{\bar D}(\mu, g,k,J_t,\nu_t )_{\bar {\kappa }}\times \{t\}$. We
denote them by $e_i^{(t)}$ and $\Upsilon^{(t)}$. Then, we
define
$$(\Xi^A_{g,k})^{(t)}=\prod_i e_i^{(t)}.\leqno(4.3)$$
\vskip 0.1in
\noindent
{\bf Definition 4.3. }{\it Let $(P_i, f_i)$ be a pseudo-submanifold of $V$.
We say that $(P_i, f_i)$ is transverse to $e_i$ (hence $\Xi^A_{g,k}$) if
$(P_i,f_i)$ is transverse to
$e_i$ as the maps from $\M_A(\mu, g,k, J,\nu)_\kappa$ and their extensions over
$\M_{\bar D}(\mu, g,k,J,\nu )_{\bar {\kappa }}$ for each $D\in
\D^{J,\nu}_{g,k}$ in the sense of
Lemma 4.2. We say that $(P_i, f_i)$ is transverse to $e_i^{(t)}$ if it is
transverse to $e_i^{(t)}$ (hence $(\Xi^A_{g,k})^{(t)}$) as the maps from
$\bigcup_t \M_A(\mu, g,k,J_t,\nu_t )_{\bar {\kappa }}\times
\{t\}$ and their
extensions over $\bigcup_t
\M_{\bar D}(\mu, g,k,J_t,\nu_t )_{\bar {\kappa }}\times \{t\}$ for each
$D\in \D^{J,\nu}_{g,k}$ in the sense of Lemma 4.2. Similarly, we say
that a pseudo-submanifold $(G, K)$ of $\overline{\M}^{\mu}_{g,k}$ is
transverse to $\Upsilon$ (or $\Upsilon^{(t)}$) if
it transverse to them as the maps from $\M_A(\mu, g,k, J,\nu)_\kappa$ and its
extensions
over $\M_{\bar D}(\mu, g,k,J,\nu )_{\bar {\kappa }}$ (or  from $\bigcup_t
\M_A(\mu, g,k,J_t,\nu_t )_{\bar {\kappa }}
\times \{t\}$ and its extensions over $\bigcup_t
\M_{\bar D}(\mu, g,k,J_t,\nu_t )_{\bar {\kappa }}
\times \{t\}$) respectively.}
\vskip 0.1in
Let us recall the construction of section 2.
Let $\alpha_i$ be integral homology
classes of $V$. We choose pseudo-submanifolds $(Y_i, F_i)$ to represent
$\alpha_i$. Let
$$Y=\prod^k_{i=1} Y_i, F=\prod^k_{i=1} F_i.$$
Then, $Y^{top}=\prod^k_{i=1} Y^{top}_i.$
Clearly, $(Y, F)$ represents $\prod^k_{i=1}\alpha_i\in H_*(V^k, \Z)$.

Note that
$\overline{\M}^{\mu}_{g,k}$ may not be smooth, but the Poincare Duality holds
over
rational coefficients. Hence, we can assume that each
homology class can be represented by a pseudo-submanifold and the
corresponding
transversality holds as long as we work over $H_*(\overline{\M}^{\mu}_{g,k},
\Q)$.
Let $(G, K)$ be a pseudo-submanifold in
$\overline{\M}^{\mu}_{g,k}$ and first we assume that  $(G, K)$ is {\em  in
general position } i.e., $K(G^{top})\subset \M^{\mu}_{g,k, I}$, where $I$
represents trivial automorphism group.
\vskip 0.1in
\noindent
{\bf Lemma 4.4. }{\it By choosing small perturbations
if necessary, we have that
$K\times F$ is transverse to $\Upsilon\times \Xi^A_{g,k}$ as the PL-maps
with respect to some triangulation of $V$ and as the smooth maps over $Y^{top}
\times G^{top}$.}
\vskip 0.1in
\noindent
{\bf Proof: } This follows obviously from Lemma 4.2.
\vskip 0.1in
\noindent
{\bf Corollary 4.5. }{\it Suppose that
$$\sum^k_{i=1} (2n-d_i)+(6g-6+2k-deg(\mu, g))=2c_1(V)(A)+2(3-n)(\mu,
g-1)+2k.\leqno(4.4)$$
and $F_i, K$ satisfy the statement of Lemma 4.4. Then,
$$\begin{array}{rl}
&K\times F \cap \Upsilon\times \Xi^A_{g,k}(\overline{\M}^r_A(\mu, g,k, J, \nu)-
\M_A(\mu, g,k, J, \nu))_I=\emptyset;\\
& K\times F (Y\times G -Y^{top}\times G^{top})\cap
\Upsilon\times \Xi^A_{g,k} =\emptyset.\\
\end{array}
 \leqno(4.5)$$}
{\bf Proof: } By Lemma 4.4, $K\times F$ is transverse to $\Upsilon\times
\Xi^A_{g,k}$. To prove
$$K\times F \cap \Upsilon\times \Xi^A_{g,k}
(\overline{\M}^r_A(\mu, g,k, J, \nu)-
\M_A(\mu, g,k, J, \nu)_I)=\emptyset,$$
by (3.25), (3.26),
$$\overline{\M}^r_A(\mu, g,k,J,\nu)=\bigcup_{D\in \D^{J,\nu}_{g,k}}\bigcup_{
\bar {\kappa}}\M_{\bar{D}}(\mu, g,k,J,\nu)_{\bar {\kappa}}.$$
By the Proposition 3.20, except the main component $\M_A(\mu, g,k,  J, \nu)_I$,
all other components $\M_{\bar{D}}(\mu, g,k,J,\nu)_{\bar {\kappa}}$ are
 smooth manifolds  of dimension
$$\leq 2c_1(V)(A)+2(3-n)(g-1)+2k-2 \leqno(4.6)$$
Since $K\times F$ is transverse to
$\M_{\bar{D}}(\mu, g,k,J,\nu)_{\bar {\kappa}}$, then
$$\dim (K\times F \cap \Upsilon\times \Xi^A_{g,k}
(\M_{\bar{D}}(\mu, g,k,J,\nu)_{\bar {\kappa}}))\leq -2.\leqno(4.7)$$
Hence, it is empty. Therefore,
$$K\times F\cap \Upsilon\times \Xi^A_{g,k}(\M_{D}(\mu, g,k,J,\nu)_{\bar
{\kappa}}))=K\times F \cap \Upsilon\times \Xi^A_{g,k}
(\M_{\bar{D}}(\mu, g,k,J,\nu)_{\bar {\kappa}}))=\emptyset.\leqno(4.8)$$
$$K\times F (Y\times G -Y^{top}\times G^{top})\cap \Upsilon\times \Xi^A_{g,k}
=\emptyset\leqno(4.9)$$
follows from a similar dimension counting argument. We leave it to the readers.

Now we adopt the notations from section 2. By (2.7),
$$(\Upsilon\times \Xi^A_{g,k}\times K\times F)^{-1}(\Delta)\subset
\M_A(\mu, g,k,  J, \nu)_I \times Y^{top}\times G^{top} \leqno(4.10)$$
is a zero-dimensional smooth submanifold.
\vskip 0.1in
\noindent
{\bf Lemma 4.6. }{\it $(\Upsilon\times \Xi^A_{g,k}\times K\times F)^{-1}
(\Delta)$ is compact and hence finite.}
\vskip 0.1in
\noindent
{\bf Proof: }{ Suppose that there is a sequence of distinct elements
$$(f_s, X_s, x_s)\in (\Upsilon\times \Xi^A_{g,k}\times K\times F)^{-1}
(\Delta).$$
By taking a subsequence, we can assume that
$$f_s\rightarrow f\in \overline{\M}^r_A(\mu, g,k,J,\nu)$$
and
$$(X_s, x_s)\rightarrow (X, x)\in Y\times G.$$
However, $(\Upsilon\times \Xi^A_{g,k}\times K\times F)^{-1}(\Delta)$ is smooth.
Thus, either
$$(f,X,x)\in K\times F \cap \Upsilon\times \Xi^A_{g,k}(\overline{\M}^r_A(\mu,
g,
k, J, \nu)-\M_A(\mu, g,k, J, \nu)_I)=\emptyset\leqno(4.11)$$
or
$$(f, X, x)\in K\times F (Y\times G -Y^{top}\times G^{top})\cap \Upsilon\times
\Xi^A_{g,k}
=\emptyset.\leqno(4.12)$$
In both cases, we have a contradiction.

Once the Lemma 4.6 is proved, as in section 2, we can define
$$\Psi^V_{(A,g,k, \mu)}(K, \alpha_1, \cdots, \alpha_k)$$
as the algebraic sum of $(\Upsilon\times \Xi^A_{g,k}\times K\times F)^{-1}
(\Delta)$. To emphasis the dependence  on $(J, \nu)$ at this
moment, we define
$$Z(A,\mu, g,k, J,\nu,F,K)=(\Upsilon\times \Xi^A_{g,k}\times K\times F)^{-1}
(\Delta)\leqno(4.13)$$
Sometimes (for example (4.29)), we also use $Y, G$ in the place of $F, K$
in $Z(\cdots)$, if there is no confusion.
To abuse the notation, we  denote its algebraic sum by
$|Z(A,\mu,g,k, J,\nu,F,K,)|$.

Next we have to prove Proposition 2.4. The proof
will be divided into a series of Lemmas:
\vskip 0.1in
\noindent
{\bf Lemma 4.7. }{\it $|Z(A,\mu,g,k,J,\nu,F,K)|$ is independent of the
representative
$(Y,F)$ and $(G,K)$, whenever $(G,K)$ is  in general position.}
\vskip 0.1in
\noindent
{\bf Proof: } Suppose that $(Y',F'),(G' ,K')$ are other representatives such
that $(G', K')$ is  in general position. There are corbordisms $(Q,H)$ and
$(L, P)$ such that
$$\partial(Q)=Y\cup -Y', \ H|_{\partial(Q)}=F\cup -F'; \partial(L)=G\cup -G',
\ P|_{\partial(L)}=K\cup -K'\leqno(4.14)$$
Let's first work on $(Q, H)$.
By choosing a small perturbation of $H$ relative to $\partial(Q)$ if
necessary, we can assume that $H\times K$ is transverse to
$\Upsilon \times \Xi^A_{g,k}$. Then, by counting dimensions  as
one did in the proof of Lemma 4.4, one can show that
$$\begin{array}{rl}
&H\times K\cap \Upsilon \times \Xi^A_{g,k}(\overline{\M}^r_A(\mu, g,k,J,\nu)-
\M_A(\mu, g,k,J,\nu)_I)=\emptyset;\\
&\ H\times K(Q\times K-Q^{top}\times K^{top})\cap
\Upsilon \times \Xi^A_{g,k}=\emptyset.\\
\end{array} \leqno(4.15)$$
Then, it follows from the same argument as in
the proof of Corollary 4.5 that
$$Z(A,\mu,g,k,J,\nu, H,K)=(\Upsilon \times \Xi^A_{g,k}\times H\times
P)^{-1}(\Delta)\subset \M_A(\mu, g,k,
J,\nu)_I\times Q^{top}\times K^{top}\leqno(4.16)$$
is a compact, smooth, oriented  $1$-manifold with boundary
$$\partial(Z(A,\mu,g,k,J,\nu, H,K))=Z(A,\mu,g,k,J,\nu, F,K)\bigcup
-Z(A,\mu,g,k,J,\nu, F',K).\leqno(4.17)$$
Hence,
$$|Z(A,\mu,g,k,J,\nu,F,K)|=|Z(A,\mu,g,k,J,\nu,F',K)|.\leqno(4.18)$$
Next, we fix a $(Y,F)$ and consider $(L, P)$. By choosing a small perturbation
of $P$ relative to $\partial(L)$,
we can assume that $(L,P)$ is  in general position and $F\times P$ is
transverse
to $\Upsilon\times \Xi^A_{g,k}$. Repeating the previous argument, we can
show that
$$|Z(A,\mu,g,k,J,\nu,F,K)|=|Z(A,\mu,g,k,J,\nu,F,K')|.\leqno(4.19)$$
\vskip 0.1in
\noindent
{\bf Lemma 4.8. }{\it $|Z(A,\mu,g,k,J,\nu,F,K)|$ is independent of generic
$(J,\nu)$.}
\vskip 0.1in
\noindent
{\bf Proof: }Let $(J', \nu')$ be another generic  pair.   Choose a generic path
$(J_t, {\nu}_t)$ connecting $(J,\nu)$ to $(J', \nu')$, such
that
$$\M_{\bar{D}}(\mu, g,k, (J_t),(\nu_t))_{\bar {\kappa}}=\bigcup_t
\M_{\bar{D}}(\mu, g,k, J_t, \nu_t)_{\bar {\kappa}} \times \{t\}$$
is a smooth, oriented cobordism between
$\M_{\bar{D}}(\mu, g,k, J, \nu)_{\bar {\kappa}}$
and $\M_{\bar{D}}(\mu, g,k, J', \nu')_{\bar {\kappa}}$.
By Lemma 4.2 and choosing a small perturbation if necessary,
we can assume that $F\times
K$ is transverse to $\Upsilon^{(t)}\times (\Xi^A_{g,k})^{(t)}$, where the
choices of $F, K$ do not affect our result by the Lemma 4.7. Then, a dimension
accounting argument shows that
$$K\times F\cap \Upsilon^{(t)} \times (\Xi^A_{g,k})^{(t)}
(\overline{\M}^r(\mu, g,k,
(J_t),(\nu_t))-\M_A(\mu, g,k,(J_t), (\nu_t))_I)=\emptyset,$$
$$ K\times F(Y\times G-Y^{top}\times G^{top})\cap \Upsilon^{(t)} \times
(\Xi^A_{g,k})^{(t)}=\emptyset. \leqno(4.20)$$
Then,
$$\begin{array}{rl}
&Z(A,\mu,g,k, (J_t), (\nu_t), F,K)=(\Upsilon^{(t)}
\times (\Xi^A_{g,k})^{(t)}\times K\times F)^{-1}(\Delta)\\
\subset &\M_A(\mu, g,k, (J_t), (\nu_t))_I
\times Y^{top}\times G^{top}\\
\end{array}
\leqno(4.21)$$
is a compact, smooth, oriented $1$-manifold with boundary. Moreover,
$$\partial(Z(A,\mu,g,k, (J_t), (\nu_t), F,K))=Z(A,\mu,g,k, J, \nu, F,K)\cup
-Z(A,\mu,g,
k, J, \nu', F,K).$$
Hence,
$$|Z(A,\mu,g,k,J',\nu',F,K)|=|Z(A,\mu,g,k,J,\nu,F,K)|.\leqno(4.22)$$
\vskip 0.1in
\noindent
{\bf Lemma 4.9. } {\it $|Z(A,\mu,g,k,J,\nu,F,K)|$ is independent of
semi-positive
symplectic deformations.}
\vskip 0.1in
\noindent
{\bf Proof:} Let $\omega_t$ be a family of semi-positive symplectic deformation
of $\omega_0=\omega$. Then, we can choose a generic path   $(J_t, \nu_t)$
such that $J_t$ is $\omega_t$-tamed. Then, this lemma follows from the
same arguments as  the proof of Lemma 4.7.
\vskip 0.1in
Next, we prove Proposition 2.9, which gives the analytic foundation
for the composition law (Theorem 2.10). Then,
the composition law can be easily
derived from Proposition 2.9 by using an observation of Witten.

First, we extend the definition of $\Psi$ to the case
that $(G,K)$ is in the
boundary, say $K(G)\subset \overline{\M}^{\mu}_{\Sigma}$, where
$\overline{\M}^{\mu}_{\Sigma}$ is one stratum of $\overline{\M}^{\mu}_{g,k}$.
If $(G,K)$ has more than
two components lying in different strata, the corresponding GW-invariant is
just the sum of the GW-invariants of each component. Therefore, without
loss of generality, we can assume that $(G,K)$ is in the one stratum of
$\overline{\M}^{\mu}_{g,k}$. Note that $\overline{\M}^{\mu}_{\Sigma}$ is also a
PL-manifold with
local  quotient singularities where the
Poincare Duality holds over the
rational coefficient. Therefore, we can replace $\overline{\M}^{\mu}_{g,k}$ by
$\overline{\M}^{\mu}_{\Sigma}$ in our construction of $\Psi$. In this case, the
proper moduli space
is $\overline{\M}^r_{\Sigma}(\mu, g,k,J,\nu)$. Note that we still have
the restriction of $\Xi^A_{g,k}$, $\Upsilon$
to $\overline{\M}^r_{\Sigma}(\mu,
g,k,J,\nu)$. We shall denote it by $\Xi^A_{\Sigma}$ and $\Upsilon_{\Sigma}$.
Then, we can choose a generic $J,\nu$ and $F$ and $K$ in general position
(inside $\overline{\M}^{\mu}_{\Sigma}$),i.e., $K(G^{top})\subset \M^{\mu}_{
\Sigma, \bar{I}}$, where $\bar{I}$ means the trivial automorphism group for
each component. Repeating the
construction of Lemma 4.4, 4.5,
we can define
$$\Psi^V_{(A,g,k, \mu)}(K, \alpha_1, \cdots, \alpha_k)$$
as the algebraic sum of
$$Z(A, \overline{\M}^{\mu}_{\Sigma}, J,\nu, F,K)=
(\Upsilon_{\Sigma}\times \Xi^A_{\Sigma}\times K\times F)^{-1}(\Delta)\subset
\M_{\Sigma}(\mu, g,k,J,\nu)_{\bar{I}}\times Y^{top}\times
G^{top}\leqno(4.23).$$

Before proving Proposition 2.9, we need the following family version of
the gluing theorem in \cite{RT}.

Let $U\subset \bar{U}\subset\M^{\mu}_{\Sigma, \bar{I}}$ be
an oriented submanifold (not necessarily closed) and
$$\phi: U\times [0, \epsilon) \rightarrow \overline{\M}^{\mu}_{g,k}$$
be a diffeomorphism such that  $\phi(U\times \{0\})=U$ and if $t\neq 0$,
 $\phi(U\times \{t\})\subset
\M^{\mu}_{g,k, I}$.
Let
$$U_t=\phi(U\times \{t\}).\leqno(4.24)$$
Define
$$\M_{A}(U_0, J, \nu)=(\Upsilon)^{-1}(U_0)\cap \M_{A,\Sigma}(\mu, g,k,J,\nu)_{
\bar{I}}.$$
For each $t\neq 0$, we define
$$\M_{A}(U_t, J, \nu)=(\Upsilon)^{-1}(U_t)\cap \M_{A}(\mu, g,k,J,\nu)_I.
\leqno(4.25)$$
Fix a generic $(J,\nu)$ and $\phi$ such that for each t, $\M_{A}(U_t, J, \nu)$
 is a smooth   manifold of dimension
$$2c_1(V)(A)+2n(1-g)+\dim U_t.\leqno(4.26)$$
Then, the following is a slight modification of Theorem 6.1, [RT].
We leave the details to the readers.
\vskip 0.1in
\noindent
{\bf Gluing Theorem: }{\it \it Let $f_0$ be any map in
${\cal M} _A ( U _0 , J , \nu )$.
Then there is
a continuous family of injective maps $T _t$ from
${\cal W}$ into ${\cal M} _A ( U _t , J , \nu   )$, where $t$
is small and ${\cal W}$ is a neighborhood of $f_0$ in
${\cal M} _A ( U _0 , J , \nu  )$,
such that (1) for any $f$ in ${\cal W}$, as $t$ goes to zero,
$T _t(f)$ converges to $f$ in $C^0$-topology on $\Sigma _0$ and
in $C^3$-topology outside
the singular set of $\Sigma _0$; (2) there are $\epsilon , \delta > 0$
satisfying : if $f^\prime$ is in ${\cal M} _A ( U_t , J , \nu  )$
and
$$d_V(f^\prime(x), f_0(y)) \leq \epsilon, ~~~{\rm whenever}~x \in \Sigma _t
\in U_t,$$
where $d_V$ is the distance functions
of a $J$-invariant metric $h_V$ on $V$ and then $f^\prime$ is
in $T_t({\cal W})$. Moreover, for each $t$, $T_t$ is an orientation-preserving
smooth map from
${\cal W}$ into ${\cal M} _A ( U _t , J , \nu   )$.}
\vskip 0.1in
\noindent
{\bf Proof of Proposition 2.9: } If $(G',K')$ is another pseudo-submanifold
representing the same homology class, we can assume
that there is a pseudo-manifold cobordism $(L, P)$ between $(G, K)$ and
$(G', K')$.

Without loss of the generality, we may assume that $(G', K')$ is
 in general position. By choosing a small perturbation relative to the
boundary if necessary, we can further assume that
$$P(L^o)\subset \M^{\mu}_{g,k, I}. \leqno(4.27)$$
Again, by counting dimensions, we can show that
$$F\times P \bigcap \Upsilon\times \Xi^A_{g,k}
(\overline{\M}^r_A(\mu, g,k,J,\nu)
-\M_A(\mu, g,k,  J, \nu)_I\cup \M_{\Sigma}(\mu,
g,k,J,\nu)_{\bar{I}})=\emptyset,$$
$$F\times P(Y\times L-Y^{top}\times L^{top})=\emptyset. \leqno(4.28)$$
Then
$$\begin{array}{rl}
&Z(A,\mu,g,k,J,\nu, F, L)=
(\Upsilon\times \Xi^A_{g,k}\times F\times P)^{-1}(\Delta)\\
\subset &\M_A(\mu, g,k,  J, \nu)_I\cup \M_{\Sigma}(\mu, g,k,J,\nu)_{\bar{I}}
\times Y^{top}\times L^{top}.\\
\end{array}
\leqno(4.29)$$
Since $L^{top}$ is a manifold with boundary, we can choose an open subset
$$\tilde{L}\subset \overline{\tilde{L}}\subset L^{top}\leqno(4.30)$$
such that
$$Z(A,\mu,g,k,J,\nu, F, L)\subset \M_A(\mu, g,k,  J, \nu)_I\cup
\M_{\Sigma}(\mu,g,k,J,\nu)_{\bar{I}}\times
Y^{top}\times \tilde{L},\leqno(4.31)$$
$$\partial{\tilde{L}}=U\cup U'$$
where $U\subset G^{top}, U'\subset (G')^{top}.$ Furthermore, there is a
uniform constant $\epsilon$ such that there are collars
$$U\times [0, \epsilon)\subset \tilde{L};~~~ U'\times [0, \epsilon)\subset
\tilde{L}.\leqno(4.32)$$
Thus,
$$Z(A,\mu,g,k,J, \nu, F, \tilde{L})=Z(A,\mu,g,k,J,\nu, F, L).\leqno(4.33)$$
Let $U_t=U\times\{t\}$ and $U'_s=U'\times \{s\}.$
It follows from
an ordinary cobordism argument, which is identical to that of Lemma 4.6,
that
$Z(A, g,k,J,\nu, \tilde{L}-U\times [0,\epsilon)\cup U'\times [0,\epsilon))$
is a smooth, compact, oriented 1-manifold with boundary
$$\begin{array}{rl}
&\partial(Z(A, g,k,J,\nu, \tilde{L}-U\times [0,\epsilon)\cup U'\times [0,
\epsilon)))\\
=&Z(A,\mu,g,k,J,\nu, F,
U_{\epsilon})\cup - Z(A,\mu,g,k, J,\nu, F, U'_{\epsilon}).\\
\end{array}
\leqno(4.34)$$
Therefore,
$$|Z(A,\mu,g,k,J,\nu, F, U_{\epsilon})|=|Z(A,\mu,g,k,J,\nu, F,
U'_{\epsilon})|.\leqno(
4.35)$$
Obviously, we can use any $U_t, U'_s$ in the place of $U_{\epsilon}, U'_{
\epsilon}$. Then, we  show
$$|Z(A,\mu,g,k,J,\nu, F, U_t)|=|Z(A,\mu,g,k,J,\nu, F, U'_{s})|\leqno(4.36).$$
Next, we study the behavior of $Z(A,\mu,g,k,J,\nu,U_t)$ as $t\rightarrow 0$.
By the Gromov-Uhlenbeck
compactness theorem, and by taking a subsequence if necessary,
we may assume that any sequence $f_t\in Z(A,\mu,g,k,J,\nu,U_t)$ converges to a
$$f\in \overline{\M}_{\Sigma}(\mu, g,k,J,\nu)_{\bar{I}}\times Y\times K.
\leqno(4.37)$$
By (4.28),
$$f\in \M_{\Sigma}(\mu, g,k,J, \nu)_{\bar{I}}\times Y^{top}\times
K^{top}.\leqno(4.38
)$$
Hence, $f\in Z(A,\mu,g,k,J,\nu, U_0)$ by (4.31). On the other hand, by
the gluing
theorem, for $t$ small, there is a bijective map
$$T_t: Z(A,\mu,g,k,J,\nu, U_0)\rightarrow Z(A,\mu,g,k,J,\nu, U_t)\leqno(4.39)$$
preserving the orientation. Therefore,
$$|Z(A,\mu,g,k,J,\nu, U_0)|=|Z(A,\mu,g,k,J,\nu, U_t)|. \leqno(4.40)$$
Since $(G',K')$ is  in general position, we can allow $s=0$ in the (4.36).
Therefore,
$$|Z(A, g,k, J,\nu, U_0)|=|Z(A,\mu,g,k,J,\nu,U'_0)|. \leqno(4.41)$$
Thus, we finish the proof of Proposition 2.9.
\vskip 0.1in
\noindent
{\bf Proof of Theorem 2.10. } Let's consider the first case where we have the
embedding
$$\theta_S: \overline{\M}_{g_1, k_1+1}\times \overline{\M}_{g_2, k_2+1}$$
by identifying the $k_1+1$-th marked point of the first component with the
$1$-st marked point of the second component. Suppose that the marked points
in the first component are $\{x_1, \cdots, x_{k_1}, y_1\}$ and the marked
points in the second component are
$\{y_2, x_{k_1+1}, \cdots, x_{k_1+k_2}\}$.
The map $\theta_S$ identifies $y_1$ with $ y_2$ and gives
rise to a stable curve of
genus $g_1+g_2$ with marked points $\{x_1, \cdots x_{k_1+k_2}\}$.
By (2.17), for their finite covers,
$$\theta_S^*\overline{\M}^{\mu}_{g,k}=\overline{\M}^{\mu}_{g_1, k_1+1}\times
\overline{\M}^{\mu}_{g_2, k_2+1}.$$
Let $[K_i]\in H_*(\overline{\M}^{\mu}_{g_i,k_i+1}, \Q)$ represented by the
pseudo-submanifold $K_i$. By Proposition 2.9,
$$\Psi^V_{A,g,k, \mu}([\theta_S(K_1\times K_2)]; \{\alpha_i\})$$
doesn't depend on the particular representative of $[\theta_S(K_1\times K_2)]$.
In particular, it could have a representative,
which is  in general position. This is very
important in the later applications. Here, we choose a representative
$\theta_S(K_1\times K_2)$. Let
$$\overline{\M}^{\mu}_{\Sigma}=\theta_S
(\overline{\M}^{\mu}_{g_1, k_1+1}\times \overline{
\M}^{\mu}_{g_2, k_2+1})\leqno(4.42)$$
Then,
$$Z(A, \mu, g,k,J,\nu, F, \theta_S(K_1\times K_2))=Z(A, \theta_S(\overline{\M}^
{\mu}_{g_1, k_1+1}\times \overline{\M}^{\mu}_{g_2, k_2+1}), J, \nu, F,
\theta_S(K_1\times K_2)).
\leqno(4.43)$$
Recall (4.23) that
$$Z(A,\overline{\M}^{\mu}_{\Sigma},J,\nu, F, \theta_S(K_1\times
K_2))\leqno(4.44)$$
$$\mbox{\hskip 0.5in} =(\Upsilon_{\Sigma}\times \Xi^A_{\Sigma}\times F\times
\theta_S(K_1\times K_2))^{-1}(\Delta)\subset \M_{\Sigma}(\mu, g,k,J,\nu)_{
\bar{I}}
\times Y^{top}\times \theta_S(K_1\times K_2)^{top}.$$
But any element of $\M_{\Sigma}(\mu, g,k,J,\nu)_{\bar{I}}$ is a pair
$$(f_1, f_2)\in \M_{A_1}^{\mu}
(\mu, g_1,k_1+1,J,\nu)_I\times \M^{\mu}_{A_2}(\mu,g_2,k_2+1,
J,\nu)_I; ~~f_1(y_1)=f_2(y_2)\leqno(4.45)$$
with $A_1+A_2=A$. Let $e^{A_1}_{y_1}$ be the evaluation map of $y_1$'s and
$e^{A_2}_{y_2}$ be  the evaluation map of $y_2$'s. Then,
$$\M_{\Sigma}(\mu, g,k,J,\nu)_{\bar{I}}=\bigcup_{A_1+A_2=A}
(e^{A_1}_{y_1}\times
e^{A_2}_{y_2})^{-1}(\delta), \leqno(4.46)$$
where $\delta$ is the ordinary diagonal in $V\times V$.
Using this decomposition and switching the order of the components
appropriately, we have
$$\bigcup_{A_1+A_2=A}\Upsilon_{A_1}\times \Upsilon_{A_2}\times \Xi^{A_1}_{g_1,
k_1}\times \Xi^{A_2}_{g_2, k_2}\times K_1\times K_2\times F \times
e^{A_1}_{y_1}
\times e^{A_2}_{y_2},\leqno(4.47)$$
which maps the space
$$\bigcup_{A_1+A_2=A}\M_{A_1} (\mu, g_1, k_1+1, J, \nu)_I\times
\M_{A_2}(\mu, g_2, k_2+1,  J, \nu)_I\times K^{top}_1\times K^{top}_2\times
Y^{top}
$$
into
$$\mbox{\hskip 1in} \M^{\mu}_{\Sigma_1,I}\times \M^{\mu}_{\Sigma_2,I}\times
V^{k_1}\times
V^{k_2}\times \M^{\mu}_{\Sigma_1,I}\times
\M^{\mu}_{\Sigma_2,I}\times V^{k_1+k_2}\times V \times V.$$
We will denote this map
by $\Upsilon_{\Sigma}\times \Xi^A_{\Sigma}\times F\times \theta_S(K_1\times
K_2)$. Let
$$\Delta_{A_1, A_2}\subset (\M^{\mu}_{\Sigma_1,I}\times \M^{\mu}_{\Sigma_2,I}
\times V^{k_1}\times V^{k_2})\times ( \M^{\mu}_{\Sigma_1,I}\times
\M^{\mu}_{\Sigma_2,I}\times V^{k_1+k_2})\leqno(4.48)$$
be the diagonal.
Define
$$\begin{array}{rl}
&Z(A_1, A_2, g_1, g_2, k_1, k_2, \delta, J,\nu)\\
=&(\Upsilon_{A_1}\times
\Upsilon_{A_2}\times \Xi^{A_1}_{g_1, k_1}\times \Xi^{A_2}_{g_2, k_2}\times
K_1\times K_2\times F \times e^{A_1}_{y_1}\times e^{A_2}_{y_2})^{-1}(
\Delta_{A_1, A_2}\times \delta).
\end{array}
\leqno (4.49) $$
Then,
$$Z(A,\overline{\M}^{\mu}_{\Sigma},J,\nu, F, \theta_S(K_1\times
K_2))=\bigcup_{A_1+
A_2=A}Z(A_1, A_2, g_1, g_2, k_1, k_2, \delta, J,\nu)\leqno(4.50)$$
$$=\bigcup_{A_1+A_2=A}(\Upsilon_{A_1}\times
\Upsilon_{A_2}\times \Xi^{A_1}_{g_1, k_1}\times \Xi^{A_2}_{g_2, k_2}\times
K_1\times K_2\times F \times e^{A_1}_{k_1+1}\times e^{A_2}_{k_2+1})^{-1}(
\Delta_{A_1, A_2}\times \delta). $$
A corbordism argument as that of Lemma 4.6 shows that
$|Z(A_1, A_2, g_1, g_2,
k_1, k_2, \delta, J,\nu)|$ only depends on the homology class of
$[\delta]\in H_*(V\times V, \Z)$. Choose a homogeneous basis $\{\beta_b\}_{1
\leq b\leq L}$ of $H_*(V, \Z)$ modulo torsion. Let $(\eta_{ab})$ be its
intersection matrix. Note that $\eta_{ab}=\beta_a.\beta_b=0$ if the dimension
of $\beta_a$ and $\beta_b$ are not complementary to each other. Let
$(\eta^{ab})$ be the inverse of $(\eta_{ab})$. Then,
$$[\delta]=\sum_{a,b} \eta^{ab}\beta_a\otimes \beta_b. \leqno(4.51)$$
Choose a pseudo-submanifold representing $\beta_a$ (still denoted it by
$\beta_a$). Then, one observes that
$$|Z(A_1, A_2, g_1, g_2, k_1, k_2, \delta, J,\nu)|\leqno(4.52)$$
$$=\sum_{a,b} \eta^{ab}|Z(A_1, \mu, g_1,  k_1, \beta_a, J,\nu, \prod_{i\leq
k_1}
Y_i, K_1)||Z(A_2, \mu, g_2, k_2, \beta_2, J,\nu, \prod_{j>k_1} Y_j, K_2)|.$$
Together with (4.50), it yields the first formula of Theorem 2.10
$$\begin{array}{rl}
&\Psi ^V_{(A,g,k,\mu)}(\theta _{S*}[K_1\times K_2];\{\alpha _i\})\\
=~& \sum \limits _{A=A_1+A_2} \sum \limits_{a,b}
\Psi ^V_{(A_1,g_1,k_1+1, \mu)}([K_1];\{\alpha _{i}\}_{i\le k}, \beta _a)
\eta ^{ab}
\Psi ^V_{(A_2,g_2,k_2+1,\mu)}([K_2];\beta _b,
\{\alpha _{j}\}_{j> k}) \\
\end{array}
\leqno (4.53)
$$
The second formula can be derived in the similar fashion. Here, we have an
embedding
$$\mu: \overline{\M}_{g-1, k+2}\rightarrow \overline{\M}_{g, k}\leqno(4.54)$$
by gluing last two marked points $x_{k+1}, x_{k+2}$. By (2.18), $\mu$ induces
a map on the corresponding finite covers
$$\bar{\mu}: \overline{\M}^{m_{\mu}}_{g-1, k+2}\rightarrow
\overline{\M}^{\mu}_{
g, k}.$$
We choose
a representative
$\bar{\mu}(K)$. Let
$$\overline{\M}^{\mu}_{\Sigma}=\bar{\mu}(\overline{\M}^{m_{\mu}}_{g,
k+2}).\leqno(4.55)$$
Notes that
$$\M_{\Sigma}(\mu, g,k,I,J,\nu)=(e_{k+1}\times e_{k+2})^{-1}(\delta).$$
It implies that
$$Z(A,\mu, g,k,J,\nu, F, \mu(K))=Z(A, \overline{\M}^{\mu}_{\Sigma}, J,\nu, F,
\mu(K))\leqno(4.56)$$
$$=Z(A, m_{\mu}, g-1, k+2, \delta, J, \nu, F, K).$$
By (4.52),
$$|Z(A,\mu, g,k,J,\nu, F, \mu(K))|=\sum_{a,b} \eta^{ab}|Z(A, m_{\mu}, g-1,
k+2,  J,\nu,
Y\times \beta_a\times \beta_b, K)|.\leqno(4.57)$$
It yields the second formula
$$\Psi ^V_{(A,g,k, \mu)}(\bar{\mu}_*[K_0];\alpha _1,\cdots, \alpha _k)
=\sum _{a,b} \Psi ^V_{(A,g-1,k+2, m_{\mu})}([ K_0];\alpha _1,\cdots, \alpha _k,
\beta _a,\beta _b) \eta ^{ab}\leqno (4.58)
$$
We finish the proof of Theorem 2.10.

\section{Stabilizing Conjecture}

One of the initial motivations for studying the
GW-invariants is to use it to distinguish nondeformation equivalent symplectic
manifolds. For example, the first author had  successfully calculated the
genus-$0$ GW-invariant in \cite{R1} to produce a large number of diffeomorphic,
non-symplectic-deformation equivalent  symplectic manifolds, whose
existence were unknown in symplectic topology before. During the course of work
in \cite{R1}, the first author observed a correspondence
 between the
differentiable category of symplectic 4-manifolds $X$ and the symplectic
deformation category of its stabilized manifold $X\times \P^1$. It can be
summarized as the following stabilizing conjecture:
\vskip 0.1in
\noindent
{\bf Stabilizing Conjecture \cite{R1}: }{\it Suppose that $X,Y$ are two simply
connected
homeomorphic symplectic 4-manifolds. $X, Y$ are diffeomorphic iff the
stabilized manifolds $X\times \P^1,
Y\times \P^1$ with the product symplectic structures are deformation
equivalent.}
\vskip 0.1in
It has been demonstrated by the
Donaldson theory that the differentiable structure of smooth
4-manifolds is a delicate problem. However, by the results of Freedman, the
two simply connected 4-manifolds $X, Y$ are
homeomorphic if and only if $X\times \P^1, Y\times \P^1$ are diffeomorphic.
Therefore, the delicate problem  about the differentiable structures of smooth
4-manifolds will disappear after the stabilizing process. On the other hand,
the stabilizing conjecture can be viewed as an analogy of Freedman's theorem
between
the smooth and the symplectic category.
The first pair of examples supporting the
conjecture  were constructed in
\cite{R1}, where $X$ is the blow-up of $\P^2$ at 8-points and $Y$ is a Barlow
surface. Furthermore, the first author also verified the conjecture for the
cases: (1). $X$ is rational, $Y$ is irrational; (2). $X, Y$ are irrational but
have different number of $(-1)$ curves. Since then, a lot of more evidences
supporting the stabilizing Conjecture have been discovered.
Note that the first Chern class is an obvious symplectic deformation
invariant. The stabilizing conjecture implies that the first Chern class
of a simply connected symplectic 4-manifold is a differentiable invariant,
which was only proved
recently by Taubes \cite{T1}. Recently, the first author was informed by
Donaldson that his results on the existence of
symplectic submanifolds imply that if
$X\times \P^1, Y\times \P^1$ are symplectic deformation equivalent, then some
branched covers of $X, Y$ are diffeomorphic. In this section, we will compute
some higher genus GW-invariants to prove the stabilizing conjecture for the
case of simply connected elliptic surfaces $E^n_{p,q}$. The examples
$E^n_{p,q}\times \P^1$ were suggested to the first author by Donaldson.

Let's recall the construction of $E^n_{p,q}$.  Let $E^1$ be the blow-up
of $\P^2$ at generic 9 points, and let $E^n$ be the fiber connected
sum of $n$ copies of $E^1$. Then $E^n_{p,q}$ can be obtained from $E^n$ by
logarithmic transformations alone two smooth fibers with multiplicity $p$ and
$q$. Note that
$E^n_{p,q}$ is simply connected if and only if $p,q$ are coprime.
Moreover, the Euler number $\chi(E^n_{p,q})=12n$, and hence $n$ is a
topological number.
\vskip 0.1in
\noindent
{\bf Theorem 5.1. }{\it $E^n_{p,q}\times \P^1, E^n_{p',q'}\times \P^1$ with
product symplectic structures are symplectic deformation equivalent if and only
if $(p,q)=(p', q')$.}
\vskip 0.1in
Combining with known results about the smooth classification of $E^n_{p,q}$ by
\cite{Ba}, \cite{FM}, \cite{MO}, \cite{MM}, we can prove
\vskip 0.1in
\noindent
{\bf Corollary 5.2. }{\it The stabilizing conjecture holds for $E^n_{p,q}.$}
\vskip 0.1in
Let $F_p, F_q$ be two multiple fibers and $F$ be a general fiber. Let
$A_p=[F_p], A_q=[F_q]$. It is known that $A_p=\frac{[F]}{p}, A_q=\frac{[F]}{q}
$. The primitive class  $A=[F]/pq$. Another piece of topological information is
that the canonical class $K$ is Poincare dual to
$$(n-2)F+(p-1)F_p+(q-1)F_q=((n-2)pq+(p-1)p+(q-1)q)A.
\leqno(5.1)$$
Then, the Theorem 5.2 follows from the following calculation
\vskip 0.1in
\noindent
{\bf Proposition 5.3. }{\it
$$\Psi^{E^n_{p,q}\times \P^1}_{(mA,1,1)}(\overline{\M}_{1,1}; \alpha)=
\left\{ \begin{array}{ll}
2q(A\cdot \alpha) ;& m=q (mA=A_p), \\
2p(A\cdot \alpha); &m=p (mA=A_q), \\
 0;  & m\neq p,q \mbox{ and } m<pq,
\end{array}\right. \leqno(5.2) $$
where $\alpha$ is a 4-dimensional homology class. In particular,
$$\Psi^{E^n_{p,q}\times \P^1}_{(mA,1,1)}(\overline{\M}_{1,1}; \cdot)\neq 0
\mbox{ for } m=p,q.$$}
\vskip 0.1in
\noindent
{\bf Proof: }  By the deformation theory of elliptic surfaces, we can
choose a complex structure $J_0$ on $E^n_{p,q}$
such that all singular fibers are
nodal elliptic curves. Furthermore, we can assume that the complex structures
of multiple fibers are generic,
i.e., whose $j$-invariants are not 0 or 1728.
We shall choose $\nu=0$. Therefore, there is no need to consider finite covers.
We shall drop $\mu$ from the notation.
Let $j_0$ be the standard complex structure on $\P^1$.

Let's describe $\overline{\M}_{mA}(1,1,J_0\times j_0, 0)$ for $m< pq$. For
any $f\in \overline{\M}_{mA}(1,1,J_0\times j_0, 0)$, the image $im(f)$ is a
connected effective holomorphic
curve.
Namely, $im(f)=\sum_i a_i C_i$ where $a_i>0$ and $C_i$ are irreducible
components. Note that
$$mA=\sum_i a_i [C_i].\leqno(5.3)$$
For the product complex
structure $J_0\times j_0$, $C_i=(C^1_i, C^2_i)$ where $C^1_1\subset E^n_{p,q},
C^2_i\subset \P^1$. $C^2_i$ can be realized as a holomorphic map from either
an elliptic curve or a rational curve to $\P^1$. Hence, $[C^2_i]=p_i [\P^1]$
for $p_i\geq 0$. By (5.3), $p_i=0$ and $C^2_i=\{x\}$ for some $x\in \P^1$.
Since
$im(f)$ is connected. We can write
$$im(f)=(\sum a_i C^1_i)\times \{x\},\leqno(5.4)$$
where $\sum a_i C^1_i$ is a connected effective curve. By our assumption on
singular fibers, each $C^1_i$ is
either a multi-section or a fiber. A multi-section  has positive
intersection with a general fiber. A fiber has zero intersection with general
fiber. Since $mA\cdot [F]=0$, this implies that each $C^1_i$ is a fiber,
i.e., a general fiber, a singular fiber or a multiple fiber.
Since $m<pq$ and
a singular fiber has the same homology class of a general fiber, $C^1_i$ can
only be
a multiple fiber. Since $im(f)$ is connected, therefore, $im(f)$ is either
$F_p\times \{x\}$ or $F_q\times \{x\}$. In particular,
$$\overline{\M}_{mA}(1,1,J_0, 0)=\emptyset \mbox{ for }m\neq p,q \mbox{ and }
m<pq. \leqno(5.5)$$
Obviously, $(J_0, 0)$ is $mA$-good for such $m$'s. Hence
$$\Psi^{E^n_{p,q}\times \P^1}_{(mA,1,1)}(\overline{\M}_{1,1}; \alpha)=
 0 \mbox{ for } m\neq p,q \mbox{ and } m<pq.\leqno(5.6
)$$
Fix a marked point $y_o$,
$$\overline{\M}_{A_p}(1,1,J_0\times j_0, 0)=\{f: F_p\rightarrow E^n_{p,q}
\times \P^1 ~|~ im(f)=F_p\times \{x\} \}/\Z_2$$
$$=Aut(F_p)/\Z_2\times \P^1=F_p\times \P^1,\leqno(5.7)$$
because a general element of $\overline{\M}_{1,1}$ has automorphism
group $\Z_2$.
Recall the definition of $A_p$-goodness (Definition 2.18). Since
$$\overline{\M}^r_{A_p}(1,1,J_0\times j_0)=\M_{A_p}(1,1,J_0\times j_0,0).
\leqno(5.8)$$
Definition 2.18, (2) is automatically satisfied. Unfortunately,
Definition 2.18, (1) is not satisfied.
This can also be viewed from the fact that the virtual dimension
$$2c_1(V)(mA)+(3-n)(g-1)+2=2,$$
but we have a moduli space of real dimension 4. For each $f\in \M_{A_p}(1,1,
J_0\times j_0, 0)$, the normal bundle
$$N_f (E^n_{p,q}\times \P^1)=N_{F_p}(E^n_{p,q})\otimes T_x \P^1, \leqno(5.9)$$
where $im(f)=F_p\times \{x\}$. It is known that $N_{F_p}(E^n_{p,q})$ is a
torsion element of order $p$. Hence,
$$H^1(N_f)=T_x \P^1.\leqno(5.10)$$
Furthermore, since $f$ is an embedding, we have a short exact exact sequence
$$0\rightarrow TF_{p}\rightarrow T_f E^{n}_{p,q}\times \P^1\rightarrow N_f
\rightarrow 0. \leqno(5.11)$$
It induces a long exact sequence of cohomologies
$$ H^1(F_p)\rightarrow H^1(T_f E^{n}_{p,q}\times \P^1)\rightarrow
H^1(N_f)\rightarrow 0.\leqno(5.12) $$
Hence, the obstruction space (3.11), (3.17)
$$Coker(D_f\oplus J\cdot df)=H^1(N_f)=T_x \P^1.\leqno(5.13)$$
The obstruction bundle
$$COK= \pi_2^*T\P^1,\leqno(5.14)$$
i.e., the pull-back of tangent bundle of $\P^1$.

Now we need to use
following result which is an analogue of Proposition 5.7 in \cite{R2} for genus
zero invariants. The proofs are  identical. We shall
adapt the notation of Lemma 4.6.
\vskip 0.1in
\noindent
{\bf Proposition 5.4. }{\it Suppose that $(J_0, \nu_0)$ is not $A$-good,
but satisfies the following hypotheses:
$$(1).\  K\times F\cap \Upsilon_A\times \Xi^A_{g,k}(\overline{\M}^r_{A}(g,k,
J_0,
\nu_0)-\M_A (g,k,  J_0,\nu_0))=\emptyset\leqno(5.15)$$
and hence $Z(A,g,k,J,\nu, F, K)$ {\rm (cf. (4.13))} is compact.

(2). $\dim Coker(D_f\oplus J\cdot df)$ is constant for all $f\in
Z(A,g,k,J,\nu,F,K)$ and
$Z(A,g,k,J,\nu,F,K)$ is a smooth manifold with the dimension
$ \dim Coker(D_f\oplus J\cdot df)$. For any generic $(J,\nu)$ sufficiently
close to
$(J_0,\nu_0)$, $Z(A,g,k,J,\nu,F,K)$ is oriented cobordant to the zero set of a
transverse section of the
obstruction bundle $COK$. Hence it is dual to the Euler class of the
obstruction bundle $COK$.}
\vskip 0.1in

Now we continue the the proof of
Proposition 5.3. Note that $H_4(E^n_{p,q}\times
\P^1)=H_2(E^n_{p,q})\otimes H_2(\P^1)$. Without the loss of generality, suppose
that $\alpha=\beta\otimes [\P^1]$. Choose a smooth surface  $Y\subset
E^n_{p,q}$ such that $Y$ represents $\beta$ and intersects $F_p$ transversely.
Then,
$$Z(A,1,1,J_0,0,Y,\M_{1,1})=\bigcup_{y\in F_p\cap Y} \{y\}\times \P^1.
\leqno(5.16)$$
Then, by Proposition 5.4,
$$\Psi^{E^n_{p,q}\times\P^1}_{(A_p,1,1)}(\overline{\M}_{1,1}, \alpha)=
(A_p\cdot \alpha) e(T\P^1)
=2q(A_0\cdot \alpha).\leqno(5.17)$$

The same proof yields that
$$\Psi^{E^n_{p,q}\times\P^1}_{(A_q,1,1)}(\overline{\M}_{1,1}, \alpha)=
2p(A_0\cdot \alpha).\leqno(5.18)$$
We finish the proof of Proposition 5.3.
\vskip 0.1in
\noindent
{\bf Proof of Theorem 5.1: } First of all, if $(p,q)=(p',q')$, $E^n_{p,q},
E^n_{p',q'}$ were known to be complex deformation equivalent
as K\"ahler surfaces regardless where
we perform the logarithmic transform.
It follows that $E^n_{p,q}\times \P^1$ and $ E^{n}_{p',q'}\times \P^1$ with
product symplectic structures are deformation equivalent.

Conversely, suppose that $E^n_{p,q}, E^n_{p',q'}$ are
symplectic deformation equivalent. Then, there is a diffeomorphism
$$F: E^n_{p,q}\times \P^1\rightarrow  E^n_{p',q'}\times \P^1\leqno(5.19)$$
such that
$$\Psi^{E^n_{p',q'}\times\P^1}_{(F_*(A),1,1)}(\overline{\M}_{1,1}, F_*(
\alpha))=\Psi^{E^n_{p,q}\times\P^1}_{(A,1,1)}(\overline{\M}_{1,1}, \alpha),
\leqno(5.20)$$
and
$$F^*c_i(E^n_{p',q'}\times\P^1)=c_i (E^n_{p,q}\times\P^1); F^*p_1(E^n_{p',q'}
\times\P^1)=p_1(E^n_{p,q}\times\P^1). \leqno(5.21)$$
Let  $e_0\in H^2(\P^1, \Z)$ be the positive generator. First, we claim
$$F^*(e_0)=e_0.\leqno(5.22)$$
Suppose that
$F^*(e_0)=ne_0 + \beta$ for $\beta\in H^2(E^n_{p,q}, \Z)$. Note that the first
Pontrjagan class $p_1(E^n_{p,q}\times \P^1)=p_1(E^n_{p,q})\neq 0$ and $p_1(
E^n_{p',q'}\times \P^1)=p_1(E^n_{p'q'})\neq 0$. Let $\gamma(E^n_{p,q})\in
H^4(E^n_{p,q}, \Z)$ be such that $\gamma(E^n_{p,q})[E^n_{p,q}]=1$. Define
$\gamma(E^n_{p',q'})$ in the same way. Then $p_1(E^n_{p,q})$ is a
 nonzero multiple of $\gamma(E^n_{p,q})$ and $p_1(E^n_{p',q'})$ is a non-zero
multiple of $\gamma(E^n_{p',q'})$. Thus, $F^*\gamma(E^n_{p',q'})=\gamma(
E^n_{p,q})$. Then,
$$\begin{array}{rl}
1&=(\gamma(E^n_{p',q'})\cup e_0)[E^n_{p',q'}\times \P^1]=F^*(\gamma(E^n_{p',
q'})\cup e_0)[E^n_{p,q}\times \P^1]\\
=&\gamma(E^n_{p,q})\cup (ne_0 +\beta)[
E^n_{p,q}\times \P^1]=n.\\
\end{array}
\leqno(5.23)$$
Hence $n=1$. Furthermore, $F^*(e_0^2)=0$. Then
$(e_0+\beta)^2=2e_0
 \beta +\beta^2=0$. Therefore, $2e_0\beta=0$
and $\beta^2=0$, consequently, $\beta=0$.
$$c_1(E^{n}_{p,q}\times \P^1)=c_1(E^{n}_{p,q}) +2e_0.\leqno(5.24)$$
By (5.21), (5.22),
$$F^*(c_1(E^n_{p',q'}))=c_1(E^n_{p,q}).\leqno(5.25)$$
However,
$F$ sends primitive classes to primitive classes, so $F^*(A^*)=A^*$,
where $A^*$ is the Poincare dual of $A$. Hence,
$F_*(A)=A$ and
$$\Psi^{E^n_{p',q'}\times\P^1}_{(F_*(nA),1,1)}(\overline{\M}_{1,1}, F_*(
\alpha))=\Psi^{E^n_{p',q'}\times\P^1}_{(nA,1,1)}(\overline{\M}_{1,1}, F_*(
\alpha)).\leqno(5.26)$$
Suppose that $q<p$ and $q'< p'$.
Then, $A_p(=qA)$ and $A_q(=pA)$ are  the first and the second class of
$\{nA\}$ such that
$\Psi$ is nonzero and so are $A_{p'}$ and $A_{q'}$. Hence
$$F_*(A_p)=A_{p'}, F_*(A_q)=A_{q'}.\leqno(5.27)$$
 This implies that
$$p=p', q=q'.$$
We finish the proof of Theorem 5.1.

Even though $E^n_{p,q}\times S^2$ are diffeomorphic to each other, they may
have different first Chern classes. This problem can be resolved by blowing
up $E^n_{p, q}$ at one point. Namely, if $E^n_{p,q}\#\bar{\P}^2$ is a
blow-up of $E^n_{p,q}$ at one point, $E^n_{p,q}\#\bar{\P}^2$ are diffeomorphic
to each other and have the same first Chern class up to a diffeomorphism. By
a theorem of Wall, they have the same almost complex structure up to a
 homotopy. Furthermore, we can choose the blow-up loci away
from multiple fibers. All the calculations in Theorem 5.3 and Theorem 5.1
can be carried through without change. Then, we show that
\vskip 0.1in
\noindent
{\bf Proposition 5.5. }{\it Let $X$ be the blow-up of a simple connected
elliptic surface. Then, the smooth 6-manifold $X\times S^2$ admits infinitely
many deformation classes of symplectic structures with the same tamed almost
complex structure up to a homotopy.}
\vskip 0.1in

\section {The Generalized Witten Conjecture}

In this section,
we formulate a conjecture on the structure of our invariants. This
conjecture was originated by Witten in [W2], but he used
path integrals, which are not well accepted by mathematicians.
Our only contribution here is put his arguments on a rigorous
mathematical footing. During the course of discussions, we
also use the results obtained in this paper to
derive several other equations of the generating function rigorously.
Those equations were known to physicists \cite{W2}, \cite{Ho} in the physical
context. Most arguments in this section are due to Witten.

As before, we denote by $\overline \U_{g,k}$  the universal  curve
over $\overline \M_{g,k}$. Then each marked point $x_i$ gives rise to
a section, still denoted by $x_i$, of the fibration
$\overline \U_{g,k} \mapsto
\overline \M_{g,k}$. If ${\cal K}_{\U |\M}$ denotes the cotangent
bundle to fibers of this fibration, we put ${\cal L}_{(i)} = x_i^* (
{\cal K}_{\U |\M})$. Following Witten, we put
$$\langle \tau _{d_1,\alpha _1}\tau_{d_2,\alpha _2}\cdots \tau _{d_k,\alpha _k}
\rangle _g (q) ~=~\sum _{A \in H_2(V,\Z)}
\Psi^V_{(A,g,k)}([K_{d_1,\cdots,d_k}];
\{\alpha _i\}) \, q^A \leqno (6.1)$$
where $\alpha _i \in H_*(V,\Q)$ and $[K_{d_1,\cdots,d_k}]$ is the Poincare dual
of $c_1({\cal L}_{(1)})^{d_1} \cup \cdots \cup
c_1({\cal L}_{(k)})^{d_k}$ and $q$ is an element of Novikov ring.
Symbolically, $\tau _{d,\alpha}$'s denote ``quantum field theory operators''.
For simplicity, we only consider homology classes
of even degree. Choose a
basis $\{\beta_a\}_{1\le a\le N}$ of
$H_{*,\rm even}(V, \Z)$ modulo torsions. We introduce formal variables $t_r^a$,
where $r= 0, 1, 2, \cdots$ and $1\le a \le N$. Witten's generating
function (cf. [W2]) is now simply defined to be
$$F^V(t^a_r ; q) =
\langle e^{\sum _{r,a} t^a_r \tau _{r, \beta _a}} \rangle (q)
=\sum _{n_{r,a}} \prod _{r,a} {(t^a_r)^{n_{r,a}} \over {n_{r,a}}!}
\left \langle \prod _{r,a} \tau _{r,\beta _a}^{n_{r,a}}
\right \rangle (q)
\leqno (6.2)$$
where $n_{r,a}$ are arbitrary collections of nonnegative integers,
almost all zero, labeled by $r, a$. The summation in (6.2) is over all
values of the genus $g$ and all homotopy classes $A$ of $(J,\nu)$-maps.
Sometimes, we write $F_g^V$ to be the part of $F^V$ involving only
GW-invariants of genus $g$. It is clear that $F^V$ is a generalization
of the prepotential $\Phi^V = F^V_0$ of genus 0 invariants
(cf. \cite{RT}, section 9). Indeed this generalized function contains
more information on the underlying manifold, for instance, using
Taubes' theorem \cite{T2}, one observes (cf. [T])
that for a minimal algebraic surface
$V$ of general type,
$$F^V(t^a_r;q) ~=~F^V(t^a_r; 0) +  q^{K_V} e^{\tau _{0,0}} + \cdots,\leqno (
6.3)$$
while $\Phi^V$ depends only the intersection form of $V$.

One of Witten's goals is to find out the equations which $F^V$ satisfies.
The case that $V$ is a point corresponds to
the topological gravity, where
$F^V$ is governed by a complete set of solution- the KdV hierarchy,
conjectured by Witten ([W2]) and
clarified by Kontsevich ([Ko]). In general,
it is not clear what is (or there
exists at all) the
complete set of equations which
$F^V$ satisfies, though there are partial results for $V=\C P^1$ (see \cite{EY}
). However, in [W2], Witten made a conjecture on $F^V$, which we will
describe in this section.

First of all, we have obtained several important
recursion formulas in section
2 about $\Psi$. In general, we can always rewrite a recursion formula as a
differential equation of the generating function. Assume that $\beta _1 = [V]$.
 Following Witten's
arguments in [W2], one can deduce from (2.15)
that $F^V$ satisfies the generalized string equation:
$${\partial F^V\over \partial t^1_0} = {1\over 2} \eta _{ab} t^a_0 t^b_0
+ \sum \limits_{i=0}^\infty \sum \limits _{a} t^a_{i+1}
{\partial F^V\over \partial t^a_i}
\leqno (6.4)$$

For the reader's convenience, we reproduce the arguments here.
\vskip 0.1in
\noindent
{\bf Lemma 6.1. }{\it Suppose that $(V, \omega)$ is a
semi-positive symplectic
manifold. Then, the generating function $F^V$ satisfies
the generalized string
equation
$${\partial F^V\over \partial t^1_0} = {1\over 2} \eta _{ab} t^a_0 t^b_0
+ \sum \limits_{i=0}^\infty \sum \limits _{a} t^a_{i+1}
{\partial F^V\over \partial t^a_i}$$}
{\bf Proof: } By (2.15),
$$\Psi^V_{(A,g,k+1)}([K_{d_0,d_1, \cdots, d_k}]; [V], \alpha_1, \cdots,
\alpha_k)
=\Psi^V_{(A,g,k)}([\pi(K_{d_0,d_1, \cdots, d_k})];  \alpha_1, \cdots,
\alpha_k),\leqno(6.5)$$
where for convenience, we choose to forget the first marked point instead of
the last marked point as in Proposition 2.15.
We choose $d_0=0$,
i.e., $c_1(\L_{(1)})^{d_0}=1$. Next, let us find
$[\pi(k_{d_0,d_1, \cdots, d_k})]$ for
$(g,k) \not= (0,2), (1,0)$.

Let's use $\L'_{(j)}$ to denote the
corresponding line bundle over $\overline{\M}_{g,k}$. Then, it is natural to
compare $\L_{(j)}$ with $\pi^*\L'_{(j)}$. It was known in algebraic
geometry that
$$\L_{(j)}=\pi^*\L'_{(j)}+D_j, \leqno(6.6)$$
where $D_i$ is the divisor consisting of the stable curves where $x_0, x_i$ are
in a rational component with only three special points. Hence,
$$c_1(\L_{(j)})^m=c_1(\pi^*\L'_{(j)})^m+D_j\sum_{i=1}^{m-1}c_1(\L_{(j)}
)^ic_1(\pi^*\L'_{(j)})^{m-i-1}.\leqno(6.7)$$
Furthermore,
$$c_1(\L_{(j)})\cap [D_j]=0; [D_i]\cap [D_j]=0 \mbox{ for } i\neq
j.\leqno(6.8)$$
Therefore,
$$c_1(\L_{(j)})^m=c_1(\pi^*\L'_{(j)})^m+D_jc_1(\pi^*\L'_{(j)})^{m-
1}\leqno(6.9)$$
and
$$\begin{array}{rl}
&c_1(\L_{(1)})^{d_1}\cup\cdots\cup c_1(\L_{(d_k)})^{d_k}\\
=~&(\pi^*\L'_{(1)}
)^{d_1}\cup \cdots\cup (\pi^*\L'_{(1)})^{d_k}+\sum_{j=1}^k ([D_j]\cap
\bigcup_{i=1}^n c_1(\pi^*\L_{(i)})^{d_i-\delta_{ij}}).\\
\end{array}
\leqno(6.10)$$
Note that $[\pi(D_j)]=[\overline{\M}_{g,k}].$ Then,
$$\Psi^V_{(A,g,k)}([\pi(K_{d_0,d_1, \cdots, d_k})];  \alpha_1, \cdots,
\alpha_k)\leqno(6.11)$$
$$=\Psi^V_{(A,g,k)}([K_{d_1, \cdots, d_k}];  \alpha_1, \cdots,
\alpha_k)+\sum_{j=1}^k\Psi^V_{(A,g,k)}([K_{d_1, \cdots, d_j-1, \cdots, d_k}];
\alpha_1, \cdots, \alpha_k). $$
For the dimension reason, the first term is zero. Therefore, if
$(g,k) \not= (0,2), (1,0)$, we have
$$<\tau_{0,1}\prod_{i=1}^k\tau_{d_i, \alpha_i}>=\sum_{j=1}^k<\prod_{i=1}^k
\tau_{d_i-\delta_{i,j}, \alpha_i}>, \leqno(6.12)$$
where we simply define $\tau_{r, \alpha}=0$ for $r<0$.

There are two exceptional
cases for the previous arguments, namely, $g=0, k=2$ and $g=1, k=0$.
In those special cases, forgetting one
marked point will result in a unstable curve. For the special case
$g=0, k=2$, since $\overline{\M}_{0,3}=pt$, the only non-zero term is
$$\Psi^V_{(0, 0,3)}([\overline{\M}_{0,3}]; [V], \beta_a, \beta_b).$$
Moreover, one can show that
$$\Psi^V_{(0, 0,3)}([\overline{\M}_{0,3}]; [V], \beta_a, \beta_b)=\eta_{a,b}.
\leqno(6.13)$$
In the case that $g=1, k=0$, for the dimension reason, we have that
$$\Psi^V_{(A,1,1)}([\overline{\M}_{1,1}]; [V])=0.\leqno(6.14)$$
Therefore, we have an equation
$$<\tau_{0,1}\prod_{i=1}^k\tau_{d_i, \alpha_i}>=\sum_{j=1}^k<\prod_{i=1}^k
\tau_{d_i-\delta_{i,j}, \alpha_i}>+\delta_{k,2}\delta_{d_1,0}\delta_{d_2,0}
\eta_{a_1,a_2}, \leqno(6.15)$$
The corresponding equation for
the generating function is the generalized string
equation
$${\partial F^V\over \partial t^1_0} = {1\over 2} \eta _{ab} t^a_0 t^b_0
+ \sum \limits_{i=0}^\infty \sum \limits _{a} t^a_{i+1}
{\partial F^V\over \partial t^a_i}
$$

We can choose $d_0=1$ and obtain another equation for $F^V$.
\vskip 0.1in
\noindent
{\bf Lemma 6.2. }{\it $F_g$ satisfies dilation equation
$$\frac{\partial F_g}{\partial t^1_1}=(2g-2+\sum_{i=1}^{\infty}\sum_{a}t^a_i
\frac{\partial }{\partial t^a_i})F_g+\frac{\chi(V)}{24}\delta_{g,1},
\leqno(6.16)$$
where $\chi(V)$ is the Euler characteristic of $V$.}
\vskip 0.1in
\noindent
{\bf Proof: } We  choose $d_0=1$. Repeating the analysis we just have, we get
$$c_1(\L_{(0)})\bigcup_{j=1}^k c_1(\L_{(i)})^{d_i}=c_1(\L_{(0)})\bigcup_{j=1}^k
 c_1(\L'_{(i)})^{d_i}.\leqno(6.17)$$
On the another hand, one has a natural identification
$$\alpha: \overline{\M}_{g,k+1}\cong \overline{\U}_{g,k}.\leqno(6.18)$$
Furthermore,
$$\L_{(0)}=\alpha^*(\K_{\U|\M})\otimes^n_{j=1}{\cal O}(D_j).\leqno(6.19)$$
Note that $\pi[\K_{\U|\M})]=2g-2$.
Therefore, modulo the exceptional case we have
$$<\tau_{1,1}\prod_{i=1}^k\tau_{d_i, \alpha_i}>_g=(2g-2-n)<
\prod_{i=1}^k\tau_{d_i, \alpha_i}>_g.\leqno(6.20)$$
Since $\overline{\M}_{0,3}=pt$ and $c_1(\L_{(0)})$ is a nontrivial class, the
contribution from the exceptional case $g=0, k=2$ is zero. The exceptional
case $g=1,k=0$ corresponds to
$$\Psi_{(A,1,1)}([K_1]; [V]). \leqno(6.21)$$
For the dimension reason, $A$ has to be zero. Moreover,
$$\dim \overline{\M}_{1,1}=1; [K_1]=\frac{1}{24} \{pt\}.\leqno(6.22)$$
Now we fix a generic element $(\Sigma_1, x)\in \overline{\M}_{1,1}$ and
let $J_0$ be the any almost complex structure. Then
$$K\times F \cap \Upsilon_0\times \Xi^0_{1,1}=\{f: (\Sigma, x)\rightarrow V |
Im(f)=pt\}=V.\leqno(6.23)$$
Furthermore, $Coker(D_f\oplus J\cdot df)$ can be canonically identified with
$T_{Im(f)}V$.
Therefore, $(J_0,0)$ satisfies the requirement of Proposition 5.4.
Hence, by Proposition 5.4,
$$\Psi_{(0,1,1)}([pt], [V])=e(TV)=\chi(V). \leqno(6.24)$$
Therefore, the contribution from the exceptional case is
$$\frac{1}{24}\chi(V) \leqno(6.25)$$
and
$$<\tau_{1,1}\prod_{i=1}^k\tau_{d_i, \alpha_i}>_g=(2g-2-k)<\prod_{i=1}^k
\tau_{d_i,
\alpha_i}>_g+\frac{1}{24}\chi(V)\delta_{g,1}\delta_{k,0}.\leqno(6.26)$$
In terms of the generating function, this is equivalent to
the following differential equation
$$\frac{\partial F_g}{\partial t^1_1}=(2g-2+\sum_{i=1}^{\infty}\sum_{a}t^a_i
\frac{\partial }{\partial t^a_i})F_g+\frac{\chi(V)}{24}\delta_{g,1}, \leqno(
6.27)$$
which coincides with the formula derived by Hori \cite{Ho} using a different
definition.

In the dilation equation, we have a unpleasant term $2g-2$ to prevent us to
write
it as equation of $F^V$. As pointed out by Witten, there is
a dimension constraint
$$c_1(V)(A)+(3-n)(g-1)+k=\sum_{i}(d_i+cod(\beta_{a_i})), \leqno(6.28)$$
can be rewritten as an equation
$$(\sum_i\sum_a(i-1+q_a)t^a_i\frac{\partial}{\partial
t^a_i}-c_1(A)-(3-n)(g-1))F_g=0,
\leqno(6.29)$$
where $2q_a=cod(\beta_a)$.
Combining the above equations, one can deduce
\vskip 0.1in
\noindent
{\bf Corollary 6.3. }{\it When $c_1=0$, $F^V$ satisfies dilation equation
$$\frac{\partial F^V}{\partial t^1_1}=\sum_{i=1}^{\infty}\sum_a(\frac{2}{3-n}
(i-1+q_a)+1)t^a_i \frac{\partial F^V}{\partial t^a_i}+\frac{\chi(V)}{24}.
\leqno(6.29)$$}
\vskip 0.1in
Similarly, we can also use (2.15) to derive the equations (for $d_0=0,1$) of
the generating function.

Following Witten, one can introduce
$$U = {\partial ^{2} F^V \over \partial t _{0, 1}
\partial t _{0,\sigma}},~
U' = {\partial ^{3} F^V \over \partial t _{0, 1}^2
\partial t _{0,\sigma}},~\cdots,~U^{(l)}_\sigma =
{\partial ^{l+2} F^V \over \partial t _{0, 1}^{l+1}
\partial t _{0,\sigma}},~~~~{\rm for~}~ l \ge 0
\leqno (6.30)
$$
We will regard $U^{(l)}$ to be of degree $l$. By a differential
function of degree $k$ we mean
a function $G(U, U', U'', \cdots)$ of degree in that sense. In particular,
any function of the form $G(U)$ is of degree zero, and
$(U')^2$ has degree two.
\vskip 0.1in
\noindent
{\bf Witten Conjecture.} {\it For every $g\ge 0$, there are
differential functions $G_{m,a,n,b}(U_\sigma, U_\sigma ', \cdots )$
of degree $2g$ such that
$${\partial ^{2} F_g \over \partial t _{m, a}
\partial t _{n,b}} = G_{m,a,n, b} (U_\sigma, U_\sigma ', \cdots)
\leqno(6.31)
$$
up to and including terms of genus $g$.}
\vskip 0.1in
It was pointed out by Witten (cf. [W2])
that the composition law implies
$${\partial ^3 F_0 \over \partial t_{d_1,a_1}\partial t_{d_2,a_2}
\partial t _{d_3,a_3}}= \sum _{a,b}
{\partial ^2 F_0 \over \partial t_{d_1-1,a_1}\partial t_{0,a}}
\eta ^{ab}
{\partial ^3 F_0 \over \partial t_{0,b}\partial t_{d_2,a_2}
\partial t _{d_3,a_3}}\leqno (6.32)$$
and consequently, the conjecture for $g=0$.

Recall that in the genus 0 case,  WDVV equation is a direct consequence of the
composition law. In the higher genus case, it is unclear if the composition
law is helpful at all to derive the equation and solve Witten's conjecture.
Here, we state a closely related conjecture.
\vskip 0.1in
\noindent
{\bf Conjecture 6.4. }{\it $\langle \tau _{d_1,\alpha _1}\tau_{d_2,\alpha
_2}\cdots \tau _{d_k,\alpha _k}
\rangle _g$ can be reduced to enumerative invariants $\Psi^V_{(A,g,k)}(
\overline{\M}_{g,k}; \cdots)$.}
\vskip 0.1in
\noindent
{\bf Proposition 6.5. }{\it Conjecture 6.4 holds for $g\leq 1$}
\vskip 0.1in
A special case
that $g=1$ and $V =\C P^1$ was checked in [W2].
\vskip 0.1in
\noindent
{\bf Proof:} First we assume that any $c_1({\cal L}_{(i)})$ is Poincare
dual to a divisor in ${\overline \M}_{g,k} \backslash {\M}_{g,k}$ for
$g \le 1$. Then any cycle $[K_{d_1, \cdots , d_k}]$ can be represented
by a cycle in the boundary ${\overline \M}_{g,k} \backslash {\M}_{g,k}$
so long as $d_1 + \cdots + d_k > 0$. It follows from the composition
law that $\langle \tau _{d_1,\alpha _1}\tau_{d_2,\alpha _2}\cdots \tau
_{d_k,\alpha _k}\rangle _g$ can be computed in terms of
$\langle \tau _{d_1,\alpha _1}\tau_{d_2,\alpha _2}\cdots \tau _{d_l,\alpha _l}
\rangle _{g'}$ with either $l < k$ or $g' < g$. Then the proposition
follows from a standard induction.

Next we check our assumption stated at the beginning.
Given any two points $x_1, x_2$ in $S^2$, one can construct
a canonical meromorphic section
$$s_{x_1, x_2} = {(x_1 - x_2) d z \over (z-x_1) (z-x_2)}.$$
This section has simple poles at $x_1, x_2$. Moreover, for
any $\sigma \in SL(2, \C )$,
$\sigma ^* s_{\sigma (x_1), \sigma (x_2)} = s_{x_1, x_2}$.
The moduli space $\M_{0,k}$ ($k \ge 3$) is the quotient
of $(S^2 )^k \backslash \Delta_k$ by $SL(2, \C )$, where
$SL(2, \C)$
acts on $(S^2)^k$ diagonally, and $\Delta $ denotes the set
of $(x_1, \cdots, x_k)$ with $x_i = x_j$ for some $i, j$.
Notice that the universal family $\U _{0,k}$ is
biholomorphic to $S^2 \times \M _{0, k}$. Then by the
invariance of $s_{x_1, x_2}$ under $SL(2, \C)$,
one can define a section
section $s$ of the relative cotangent bundle over $\U _{0,k}$,
such that $s(z; x_1, \cdots , x_k) = s_{x_1, x_2}(z)$.
For any $i \ge 3$, this $s$ restricts to an nonvanishing
section $s_i$ of ${\cal L}_{(i)}$
over $\M _{0,k}$, i.e., $s_i (x_1, \cdots , x_k) = s(x_i)$.
Clearly, each $s_i$ extends to be a meromorphic section
on $\M _{0,k}$. Therefore,
$c_1({\cal L}_{(i)})$ must be Poincare
dual to a divisor in ${\overline \M}_{0,k} \backslash {\M}_{0,k}$
for $i\ge 3$. Similarly, one can also show this for $i\le 3$.

Now let $g=1$. Note that each torus $T$ is a branched covering of $S^2$
of degree. There are four branched points, say $x_1$, $x_2$, $_3$, $x_4$.
Conversely, given any $x_1$, $x_2$, $_3$, $x_4$, one can construct
a branched covering $T$ with those as branched points. The resulting
torus $T$ is uniquely determined by the orbit
of $(x_1, \cdots, x_4)$ in $(S^2)^4$ by $SL(2, \C)$.
Let $\pi: T \mapsto S^2$ be the branched
covering map. Then
$\pi^*(s_{x_1, x_2} s_{x_3, x_4})$
defines a nonvanishing section $s_T$ of $K_T^2$ over $T$.
Using the invariance of $s_T$ under $SL(2, \C)$, we can easily construct
a nonvanishing section of the relative canonical bundle over
$\U_{1,1}$, which can be extended meromorphically to ${\overline \U}_{1,1}$.
It follows that any $c_1({\cal L}_{(i)})$ is Poincare
dual to a divisor in ${\overline \M}_{1,k} \backslash {\M}_{1,k}$
for any $k \ge 1$.

\end{document}